\input epsf
\documentstyle{amsppt}
\pagewidth{5.4truein}\hcorrection{0.55in}
\pageheight{7.9truein}\vcorrection{0.75in}
\TagsOnRight
\NoRunningHeads
\catcode`\@=11
\def\logo@{}
\footline={\ifnum\pageno>1 \hfil\folio\hfil\else\hfil\fi}
\topmatter
\title Perfect matchings and perfect powers
\endtitle
\author Mihai Ciucu\endauthor
\thanks Research supported in part by NSF grants DMS 9802390 and DMS 0100950.
\endthanks
\affil
  School of Mathematics, Georgia Institute of Technology\\
  Atlanta, Georgia 30332-0160
\endaffil
\abstract
In the last decade there have been many results about special families of graphs whose
number of perfect matchings is given by perfect or near perfect powers
\cite{5}\cite{11}\cite{8}. In this
paper we present an approach that allows proving them in a unified way. We use this 
approach to prove a conjecture of James Propp stating that the number of tilings of
the so-called Aztec dungeon regions is a power (or twice a power) of 13.  
We also prove a conjecture of Matt Blum stating that
the number of perfect matchings of a certain family of subgraphs of the square
lattice is a power of 3 or twice a power of 3.
In addition we obtain multi-parameter generalizations of previously known results, and 
new multi-parameter exact enumeration results. We obtain in particular a
simple combinatorial proof of Bo-Yin Yang's multivariate generalization of fortresses,
a result whose previously known proof was quite complicated, amounting to evaluation
of the Kasteleyn matrix by explicit row reduction. We also include a new multivariate exact
enumeration of Aztec diamonds, in the spirit of Stanley's multivariate version.
\endabstract
\endtopmatter
\document

\def\mysec#1{\bigskip\centerline{\bf #1}\message{ * }\nopagebreak\par\bigskip}

\def\myref#1{\item"{[{\bf #1}]}"} 
 
\def\pf{{\it Proof.\ }} 
\def\endpf{\hbox{\qed}\medskip}
\def\sendpf{\hbox{\qed}}
\def\cite#1{\relaxnext@
  \def\nextiii@##1,##2\end@{[{\bf##1},\,##2]}%
  \in@,{#1}\ifin@\def\next{\nextiii@#1\end@}\else
  \def\next{[{\bf#1}]}\fi\next}
\def\proclaimheadfont@{\smc}

\def\pf{{\it Proof.\ }}

\define\Z{{\Bbb Z}}

\define\M{\operatorname{M}}
\define\T{\operatorname{T}}
\define\wt{\operatorname{wt}}
\define\s{\operatorname{s}}

\define\de{\operatorname{d}}
\define\r{\operatorname{r}}
\define\te{\operatorname{t}}
\define\z{\operatorname{z}}
\define\y{\operatorname{y}}
\define\twoline#1#2{\line{\hfill{\smc #1}\hfill{\smc #2}\hfill}}

\def\mypic#1{\epsffile{figs/#1}}



\mysec{1. Introduction}

A lattice in the plane divides the plane into elementary regions.
A {\it tile} is the union of two elementary regions that share an edge. A {\it region}
is a connected region in the plane whose boundary consists of lattice segments. A {\it
tiling} of a region is a way to pair the elementary regions it contains into
disjoint tiles.

Tiles are allowed to carry weights, in which case the weight of a tiling is defined to
be the product of the weights of the constituent tiles. The {\it tiling generating 
function} of a region $R$, denoted $\T(R)$, is the sum of the weights of all its
tilings. 

An important motivation for studying tiling generating functions is that tilings of lattice
regions can be identified with perfect matchings of their dual graphs, and are therefore 
instances of the dimer model of statistical physics on various lattice graphs. 

In the last decade there have been many results about special families of regions on
various lattices whose number of tilings is given by perfect powers 
\cite{5}\cite{11}\cite{8}. In this
paper we present a unified approach that allows proving these cases as well as new results. 

The uniform approach in our proofs is afforded by Theorem 2.3 (the General Complementation
Theorem), a result relating the matching generating function of weighted graphs in a large class
(that includes in particular all planar bipartite graphs\footnote{Indeed, the statement of 
Theorem 2.3 clearly applies when $H$ is an arbitrary subgraph of the square grid. By performing
``vertex splittings'' (see e.g. \cite{1}) if necessary, any planar bipartite graph can be viewed
as a subgraph of the square grid, with unaffected matching generating function.}) to the matching generating function of their
``complement.'' This result is an extension of Propp's Generalized Domino-Shuffling
\cite{9}, which in turn extended the Aztec diamond case of the Complementation 
Theorem of \cite{2} to arbitrary weights by eliminating
a technical assumption in the statement of the Complementation Theorem of \cite{2}
(the assumption that the cell-factors are constant along lines of cells---see
Section 2 and \cite{2}). Our proofs employ combinatorial arguments to bring the problems to a 
form where Theorem 2.3 applies.

This more general version can in fact be
proved using the arguments employed in \cite{2} to prove the Complementation Theorem.
However, it is shorter and clearer to prove it in the spirit of \cite{9}. We
present such a proof, based on certain local graph replacement operations, in 
Section 2.

A crucial ingredient in our proofs is that certain
matrices that encode the type of each family of regions behave well under the action
of a naturally arising matrix operator $\de$ whose action is given by certain rational maps 
(see Section 3). More precisely, in all 
these examples it turns out that after applying $\de$  successively to the encoding
matrix some convenient number of times---in two of our examples, twelve and respectively
thirty(!) times---, one arrives at a matrix of the same type as 
the original one, thus obtaining simple recurrence relations that determine the tiling
generating function. 

In several of our examples we find in fact more general matrices, with entries indeterminates,
which also have this periodicity property under $\de$. These can be interpreted as parametrized
subvarieties of the algebraic variety of matrices having some fixed order under the action of 
$\de$. 

On the square lattice, Elkies, Kuperberg, Larsen and Propp considered in \cite{5}
the Aztec diamond regions and proved that the region of order $n$ has $2^{n(n+1)/2}$
tilings (for $n=4$ this region is illustrated in Figure 1.1).

\topinsert
\centerline{\mypic{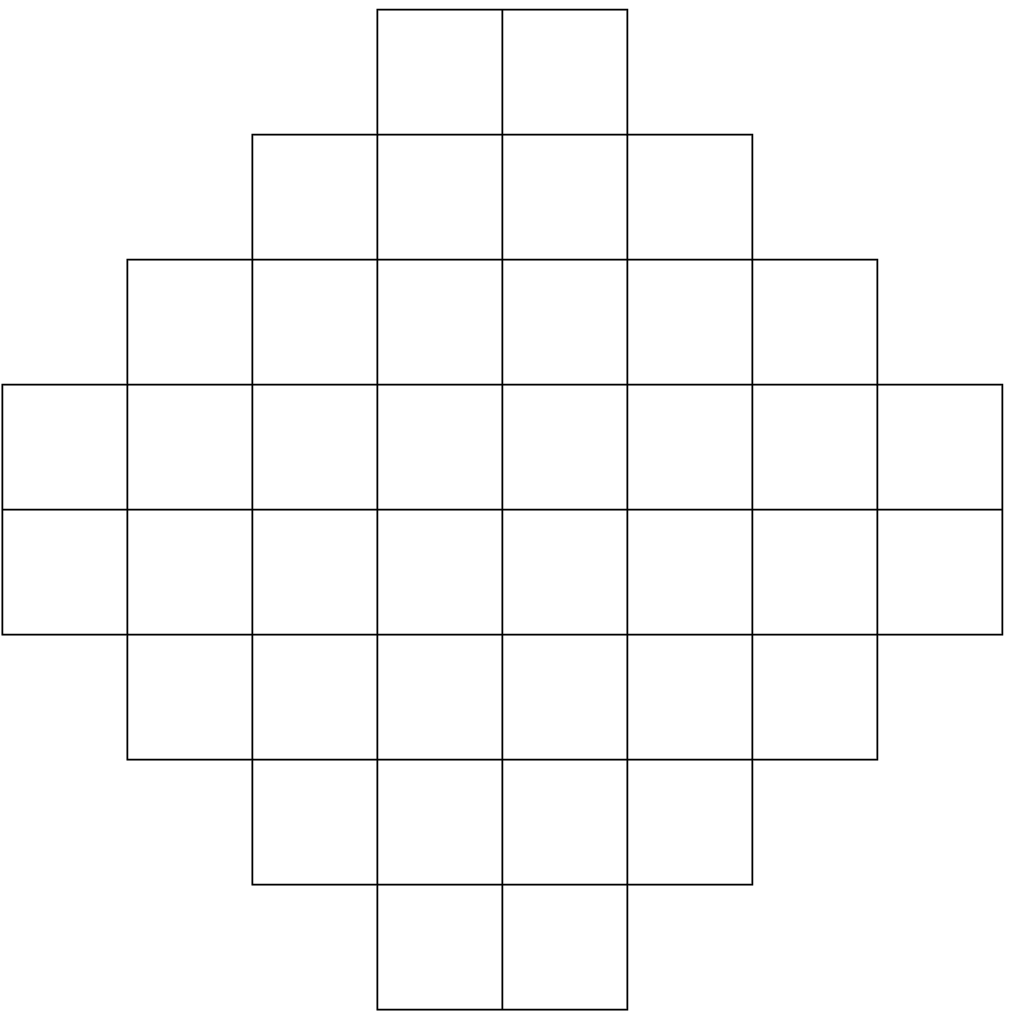}}
\centerline{{\smc Figure~1.1.} {\rm The Aztec diamond region of order 4.}}
\endinsert

For the square lattice with the diagonals drawn in---the lattice corresponding to the
affine Coxeter group $B_2$---, B-Y. Yang \cite{11} considered regions called 
fortresses, and
proved that the number of their tilings is always a power of 5 or twice a power of 5
(Figure 4.1 illustrates the case $n=4$).

On the triangular lattice, which corresponds to the affine Coxeter group $A_2$,
MacMahon \cite{7} showed that the number of tilings of the hexagon of sides $a$,
$b$, $c$, $a$, $b$, $c$ (in cyclic order) is 
$$\prod_{i=1}^a\prod_{j=1}^b\prod_{k=1}^c\frac{i+j+k-1}{i+j+k-2}.$$
It is interesting that there is no known family of regions on the triangular lattice
whose number of perfect matchings is given non-trivially by perfect powers. The
``rectangular'' region obtained from a hexagon with sides $m$, $m$, $n$, $m$, $m$, $n$
by cutting off the vertices not incident to the sides of length $n$ along zig-zag
lattice lines has $n^m$ tilings, but it can be split into $m$ hexagons that can be
tiled independently (see \cite{3}).

There is only one other affine Coxeter group of order two, $G_2$, and its
corresponding lattice is the triangular lattice with all altitudes of unit triangles
drawn in. Propp considered regions on this lattice called Aztec dungeons (see Figure~3.1 for an
example), and
conjectured that the number of their tilings is always a power of 13 or twice a
power of 13. We prove this conjecture and present two weighted
generalizations in Section 3. 

B-Y. Yang \cite{11} extended his result on enumerating tilings of fortresses by 
proving, using a quite complicated argument---explicit row reduction of the Kasteleyn
matrix---, a
conjecture of Kuperberg and Propp stating that a certain multivariate weighting of the
Aztec diamond leads to a matching generating function that factors nicely. We present
a simple combinatorial proof of this general version in Section 4, which also contains 
a weighted 
form that keeps track of the number of different types of tiles used in a
tiling of the fortress (for a simple proof of the original version of B-Y. Yang's 
result, based on the Complementation Theorem, see \cite{2}).

Section 5 contains a couterpart of Stanley's multivariate generalization (\cite{10}; see
also \cite{2}) of the Aztec diamond theorem of \cite{5}. 

In Section 6 we present a simple proof of an ex-conjecture of Propp, 
first proved by Chris Douglas \cite{4}, stating that the
number of tilings of a certain subgraph of the Aztec diamond is a power of 2. We also
obtain a natural weighted generalization that keeps track of the horizontal and
vertical edges.

In Section 7 we consider the tiling of the plane by equilateral
triangles, squares and regular hexagons. Another ex-conjecture of Propp
\cite{8}, first proved (but not yet published) by Ben Wieland,  
states that the number of tilings of a certain family of regions on this lattice,
called dragons (an example is illustrated in Figure 7.1), 
equals a power of two. We present a simple proof of this, and a
natural weighted generalization.

Section 8 presents a proof of a conjecture of Matt Blum 
stating that the number of perfect matchings of a certain family of subgraphs of the square
lattice is a power of 3 or twice a power of 3. The last section of the paper contains
some concluding remarks.

\mysec{2. The General Complementation Theorem}

A {\it perfect matching} of a graph is a collection of vertex-disjoint edges that form a
spanning subgraph. For a weighted graph $G$, the weight of a perfect matching is the
product of the weights of the constituent edges. The {\it matching generating function}
$\M(G)$ is the sum of the weights of all perfect matchings of $G$. 

We present in this section a generalization of the Complementation Theorem of \cite{2}. 
We review the necessary definitions and notations from \cite{2} below.

A {\it cellular graph} is a finite graph whose edges can be partitioned
into 4-cycles --- the {\it cells} of the graph --- such that each vertex is 
contained in at most two cells. 

Let $c_{0}$ be a cell of a graph $G$ and consider two opposite vertices $x_{0}$ and $y_{0}$
of $c_{0}$. By definition, $x_{0}$ is contained in at most one other cell
besides $c_{0}$. If such a cell $c_{1}$ exists, let~$x_{1}$ be its vertex
opposite $x_{0}$. Then $x_{1}$ in turn is contained in at most one cell 
other than $c_{1}$; if such a cell exists denote by $x_{2}$ its vertex
opposite $x_{1}$, and continue in this fashion. Repeat the procedure starting
with $y_{0}$. 

The set of cells that arise this way is said to be a {\it line}
of $G$. If the sequence $\{x_{i}\}$ (and hence the analogous sequence 
$\{y_{i}\}$) defined above is finite, the line is called a {\it path}, and
the last entry of each of the two sequences is called an {\it extremal vertex}
of $G$.

A weighted graph is a graph equipped with a weight function on its edges. 
Unless stated otherwise, the edge-weights are considered independent indeterminates.
Let $V(G)$ and $X(G)$ be the sets of vertices and extremal vertices, respectively, of $G$. 
Let $E(G)$ be the set of edges of $G$. 

The following definition is a relaxation of the corresponding one from \cite{2} 
(the original definition
involved the additional requirement that $H$ is an induced subgraph of $G$).

\proclaim{Definition 2.1}
Given a weighted graph $H$, a weighted cellular graph $G$ is said to be a {\rm cellular completion} 
of $H$ if 

$(i)$ $H$ is a subgraph of $G$, and all edges of $E(G)\setminus E(H)$ with both endpoints
in $H$ have weight 0 in $G$, and

$(ii)$ $V(G)\setminus V(H) \subseteq X(G)$.
\endproclaim

Therefore, any (not necessarily induced) subgraph of a weighted cellular graph
has a cellular completion; see Remark 2.4. A graph may have more than one cellular 
completion: Figure 2.1(c) shows a cellular completion of the graph in Figure 2.1(a);
another is obtained by shading in in Figure 2.1(b) the other class of the
chessboard indicated coloring. The cellular completion of a subgraph of the
square grid is not necessarily a subgraph of the square grid.

Let $G$ be a cellular completion of the graph $H$. The {\it complement
$H'$ of H} (with respect to $G$) is defined as follows. 
First we define the underlying graph of the complement (the weight on $H'$ is 
defined in the next paragraphs): it is the induced subgraph $H'$ of
$G$ whose vertex set is determined by the equation $V(H')\triangle V(H)=X(G)$,
where the triangle denotes symmetric difference of sets (an example is 
illustrated in Figures 2.1(a)--(d)\footnote{In Figure 2.1(d) the
graph $H'$ is indeed an induced {\it weighted} subgraph, with the missing edges
weighted by 0 (see the definition of the derived weight $\wt'$).}; the edges of $G$ not contained in $H$ 
are represented by dotted lines in Figure 2.1(c)). In other words, $V(H')$
is the set obtained from $V(H)$ after performing the following operation at 
each end of every path of $G$: if the corresponding extremal vertex belongs to 
$V(H)$, remove it; otherwise, include it. 

Starting with the weight $\wt$ on $G$, we define a new weight $\wt_1$ on $G$ as follows 
(the weight of the complement $H'$ will be defined by restriction from $\wt_1$). 
If the cell $c$ has edges weighted by $x$, $y$, $z$ and $w$ (in cyclic 
order), set $\Delta(c):=xz+yw$. We call $\Delta(c)$ the {\it cell-factor} of $c$.

\topinsert
\twoline{\mypic{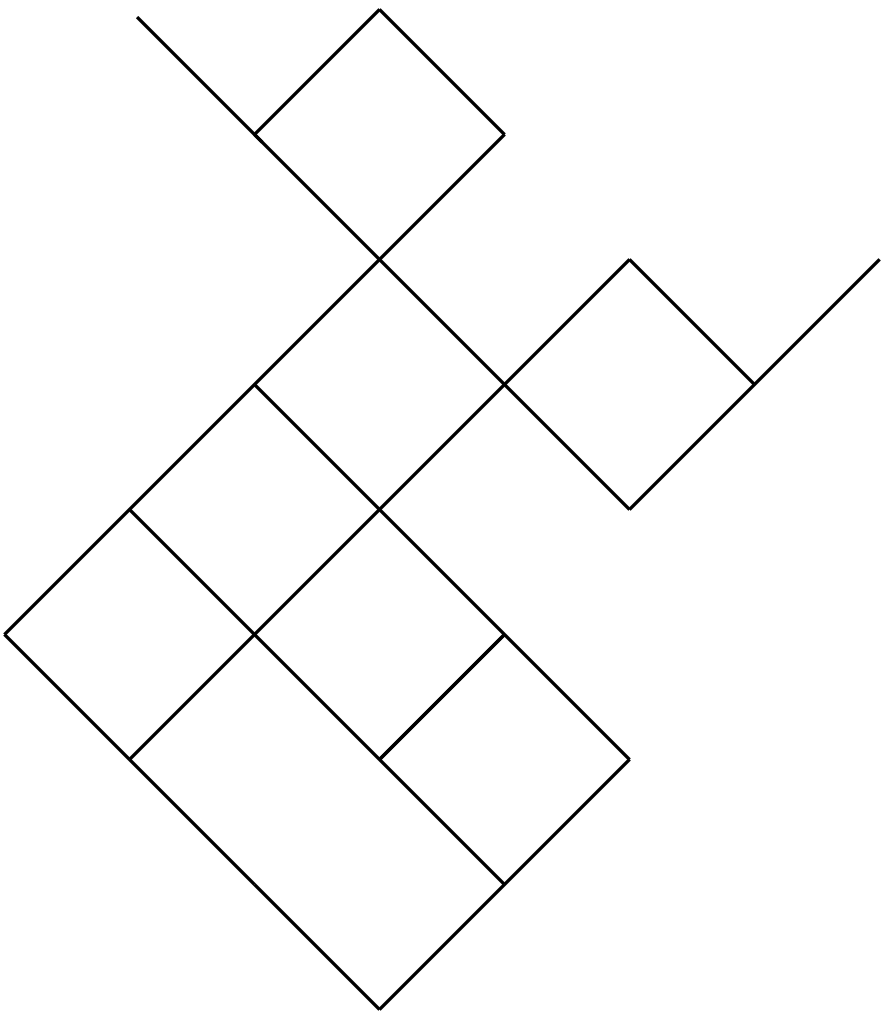}}{\mypic{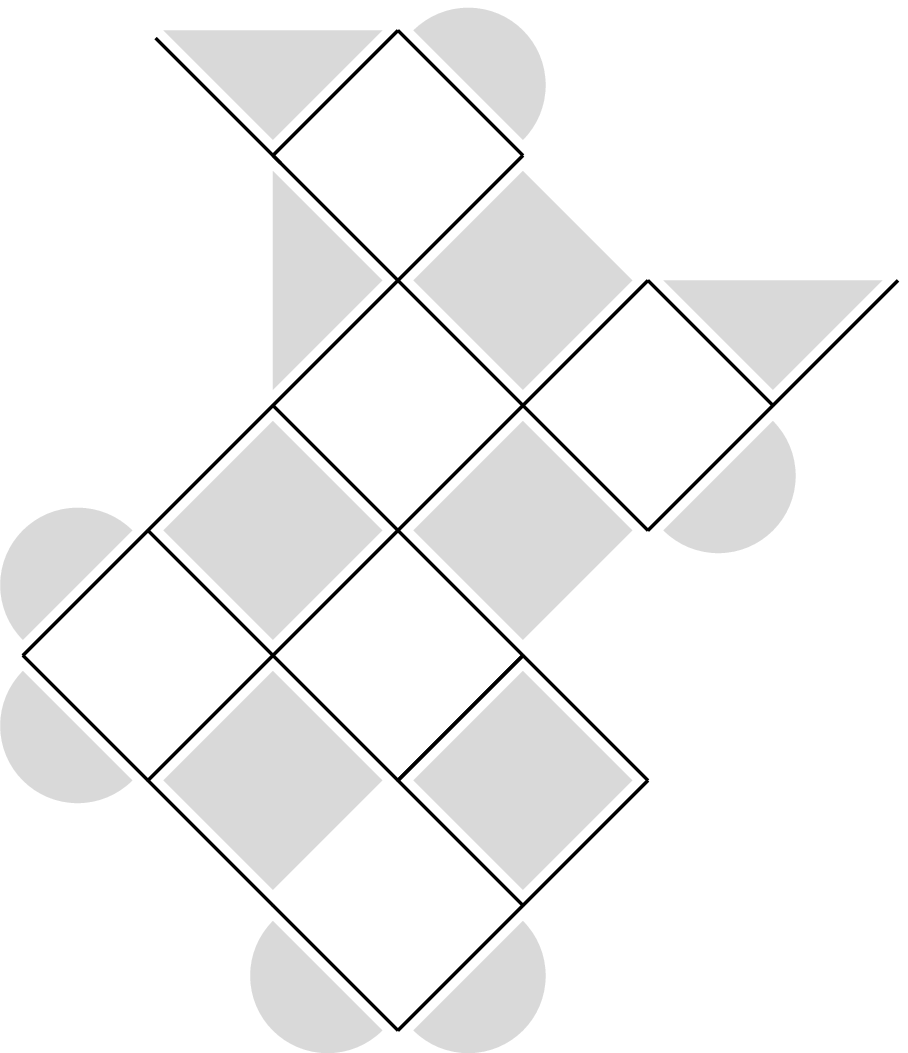}}
\twoline{Figure~2.1{\rm (a). The graph $H$.}}
{Figure~2.1{\rm (b). Whole and partial cells of $H$.}}

\medskip
\twoline{\mypic{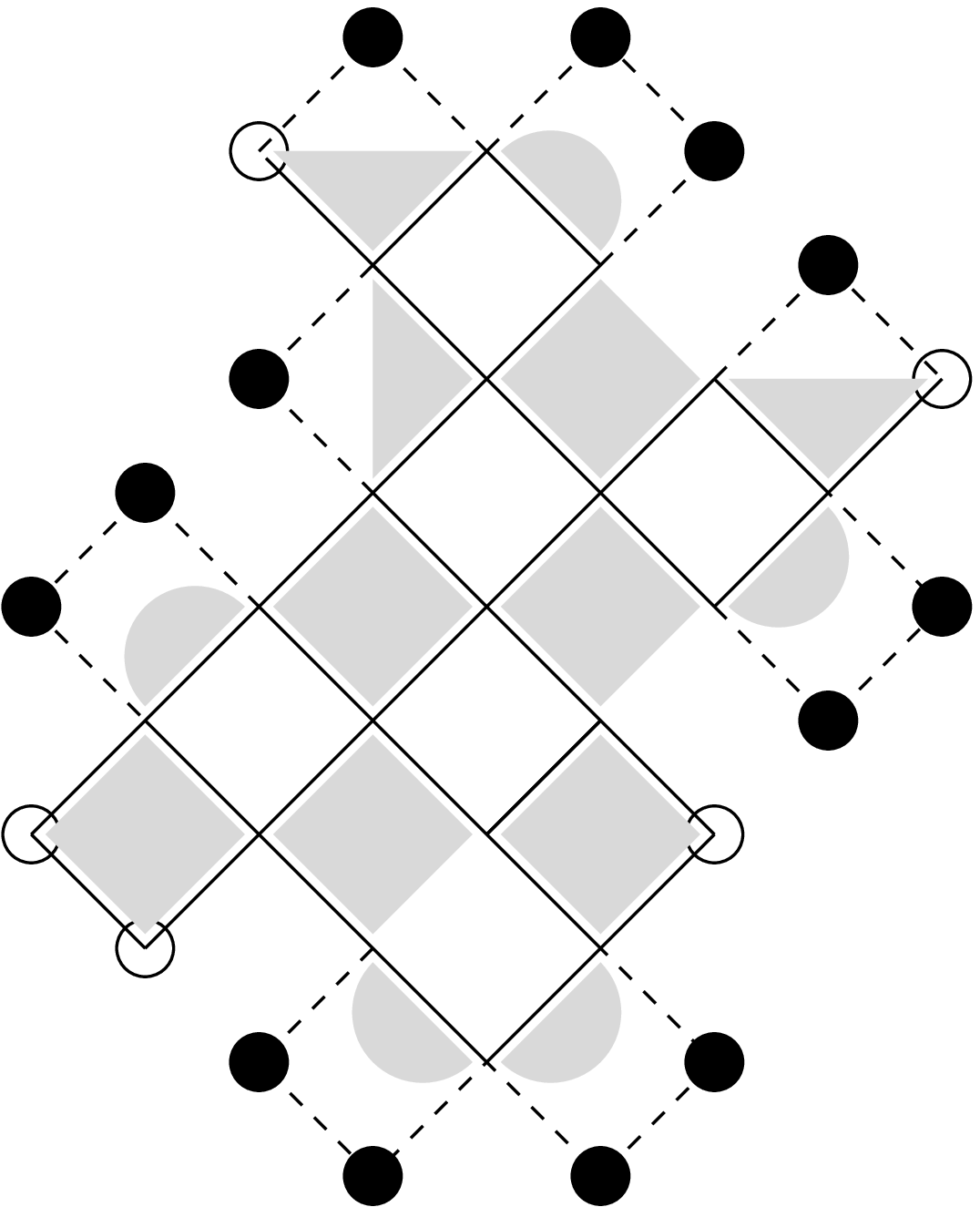}}{\mypic{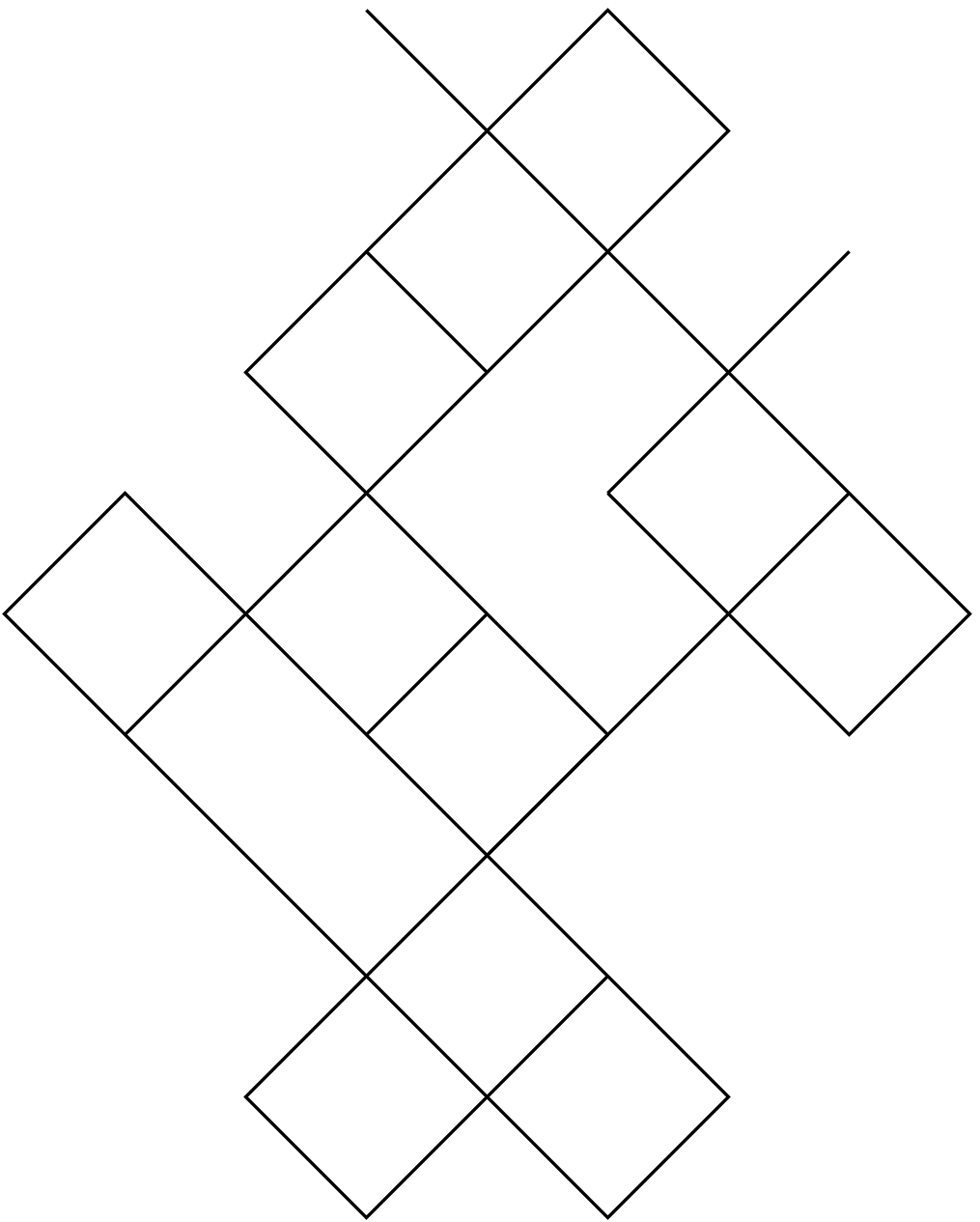}}
\twoline{Figure~2.1{\rm (c). A cellular completion of $H$.}}
{Figure~2.1{\rm (d). The graph $H'$.}}
\endinsert

Consider a cell $c$ of $G$. Assume first that the set of
endpoints of edges of $H$ contained in $c$ is the full set of four vertices of $c$---we
call such cells {\it whole}. In this
case, for each edge $e$ of $c$ set $\wt_1(e):=\wt(f)/\Delta(c)$, where $f$ is the edge 
opposite $e$ in $c$ (this is summarized in Figure 2.2, where the left part shows the
original weight on a whole cell and the right part the new weight). 
If $c$ is whole it follows that $\Delta(c)\neq0$, and this is well-defined. Note that 
this definition has the effect of sliding each 0-weight edge to the opposite edge in 
its cell, so 
by interpreting these back as missing edges, the graph changes---see Figures 2.1(a) and 
(d), where 0-weight edges are indicated as missing.

\topinsert
\centerline{\mypic{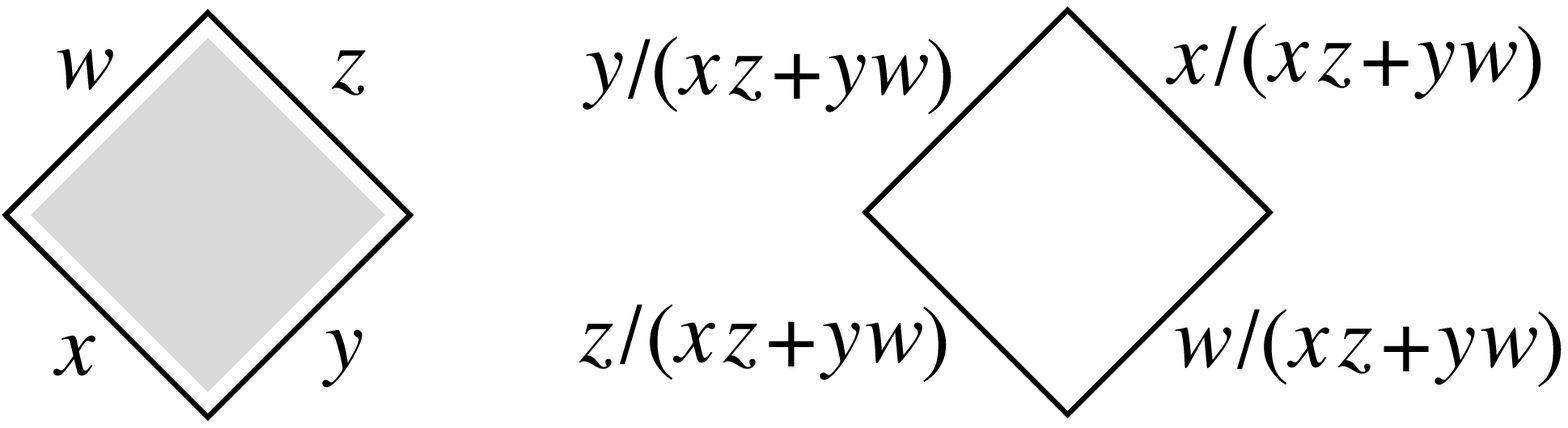}}
\centerline{{\smc Figure~2.2.} {\rm 
Weight change on whole cells.
}}

\medskip
\centerline{\mypic{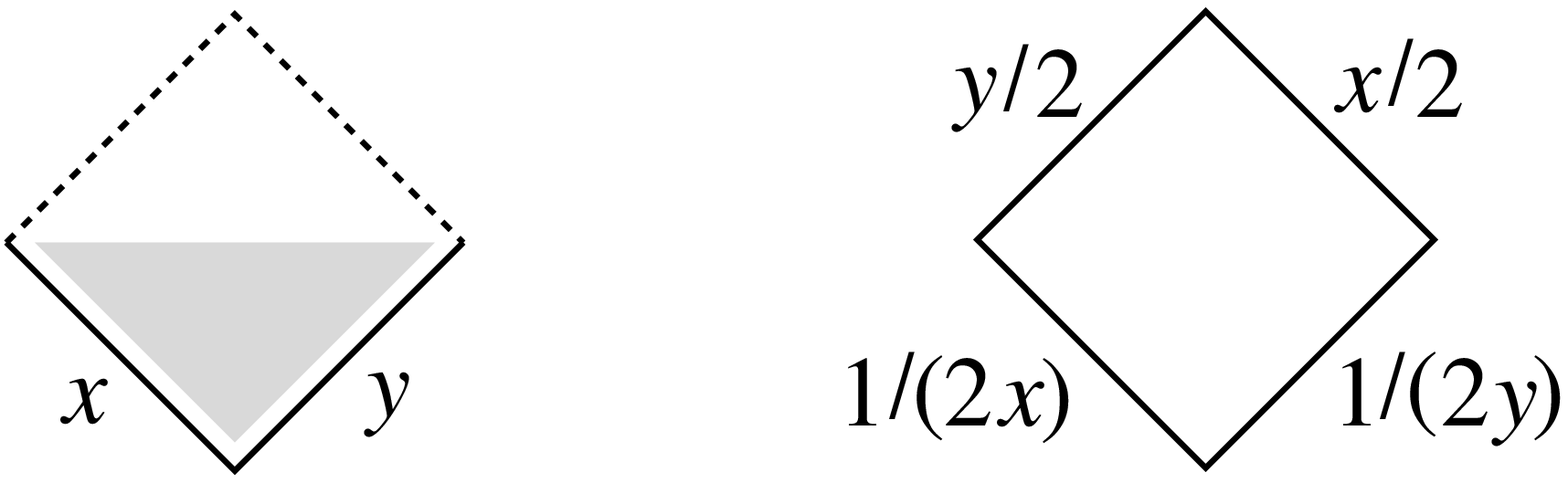}}
\centerline{{\smc Figure~2.3.} {\rm 
Weight change on partial cells with three vertices.
}}

\medskip
\centerline{\mypic{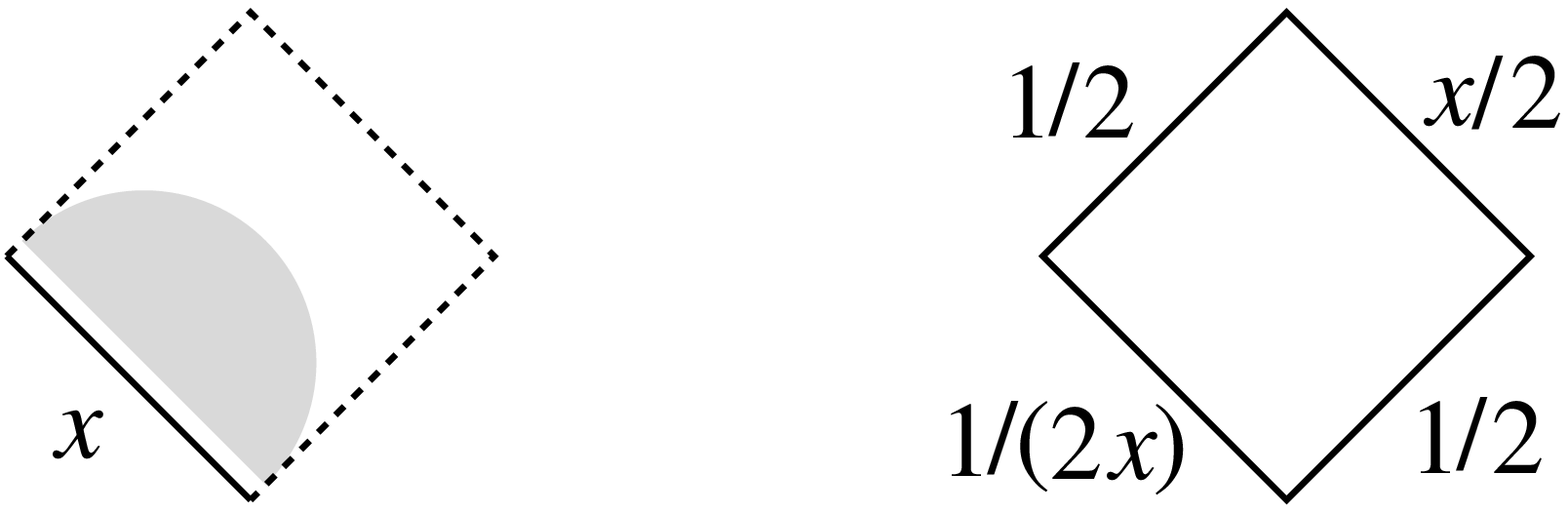}}
\centerline{{\smc Figure~2.4.} {\rm 
Weight change on partial cells with two vertices.
}}
\endinsert

If on the other hand the set of endpoints of edges of $H$ contained in $c$ has at most 
three elements---such a cell is
called {\it partial}---, define $\wt_1$ on
$c$ as indicated by Figures 2.3 and 2.4. 
More precisely, if the aforementioned set has three elements, the partial cell $c$ looks 
like the left part of Figure 2.3.
By our definition of the complement, the missing extremal point from cell
$c$ will be included in $H'$, making this a whole cell in $H'$. Assign it new weights as
indicated by the right part of Figure 2.3.

In case the set of endpoints of edges of $H$ contained in $c$ has two elements, define
$\wt_1$ similarly, using Figure 2.4. 

We note that if a partial cell contains an extremal vertex, this cell doesn't generate 
a whole cell in $H'$; in this case just use the corresponding restriction of the right side of 
Figure 2.3 or 2.4. Instances of this occur in Figures 2.1(a) and (d).

\proclaim{Definition 2.2} Let $\wt$ be a weight on the edge set of $H$. 
The {\it derived weight $\wt'$} on the complement $H'$ is the weight obtained by
restricting the above defined weight $\wt_1'$ to $H'$.
\endproclaim

\flushpar
(The above defined weight is related to, but different from the complementary weight
defined in \cite{2} and denoted there by the same symbol.) 

The following result is a generalization of Theorem 2.1 of \cite{2}. 
The main improvements are
that here the cell-factors are not required to be constant along lines, and that in the 
current set-up $H$ 
does not need to be an induced subgraph of its cellular completion $G$ (see Remark 2.4).

\proclaim{Theorem 2.3 (General Complementation Theorem)}
Let $H$ be a weighted graph and let $G$ be a cellular
completion of $H$. Let $\wt$ be a weight function on the edges of $H$. 
Then $\M(H;\wt)$ can be expressed in terms of the complement $H'$ of $H$ with respect to $G$,
weighted by the derived weight $\wt'$, as
$$
\M(H;\wt)=2^{|\Cal C_0|}\left(\prod_{c\in\Cal C_1}\Delta(c)\right)\M(H';\wt'),
$$
where $\Cal C_1$ and $\Cal C_0$ are the sets of whole and partial cells, respectively.
\endproclaim

\pf The particular case when $H$ is an Aztec diamond is due to Propp \cite{9} and
follows by repeated application of Lemma 2.5 below. For the general case we proceed as follows.
Detach the cells from one another and introduce one new vertex and two new unit-weighted 
edges between any two touching cells; 
introduce two new vertices and two new unit-weighted edges at every vertex of a cell that
doesn't belong to any other cell (the latter are precisely the extremal vertices of the cellular
completion $G$ that are contained also in $H$; for the graph in Figure 2.1(a), this procedure is
illustrated in Figures 2.5(a) and (b)). By the ``vertex splitting'' trick (see e.g. \cite{1}), the
matching generating function of the resulting graph $H_1$ is the same as that of $H$.

\topinsert
\centerline{\mypic{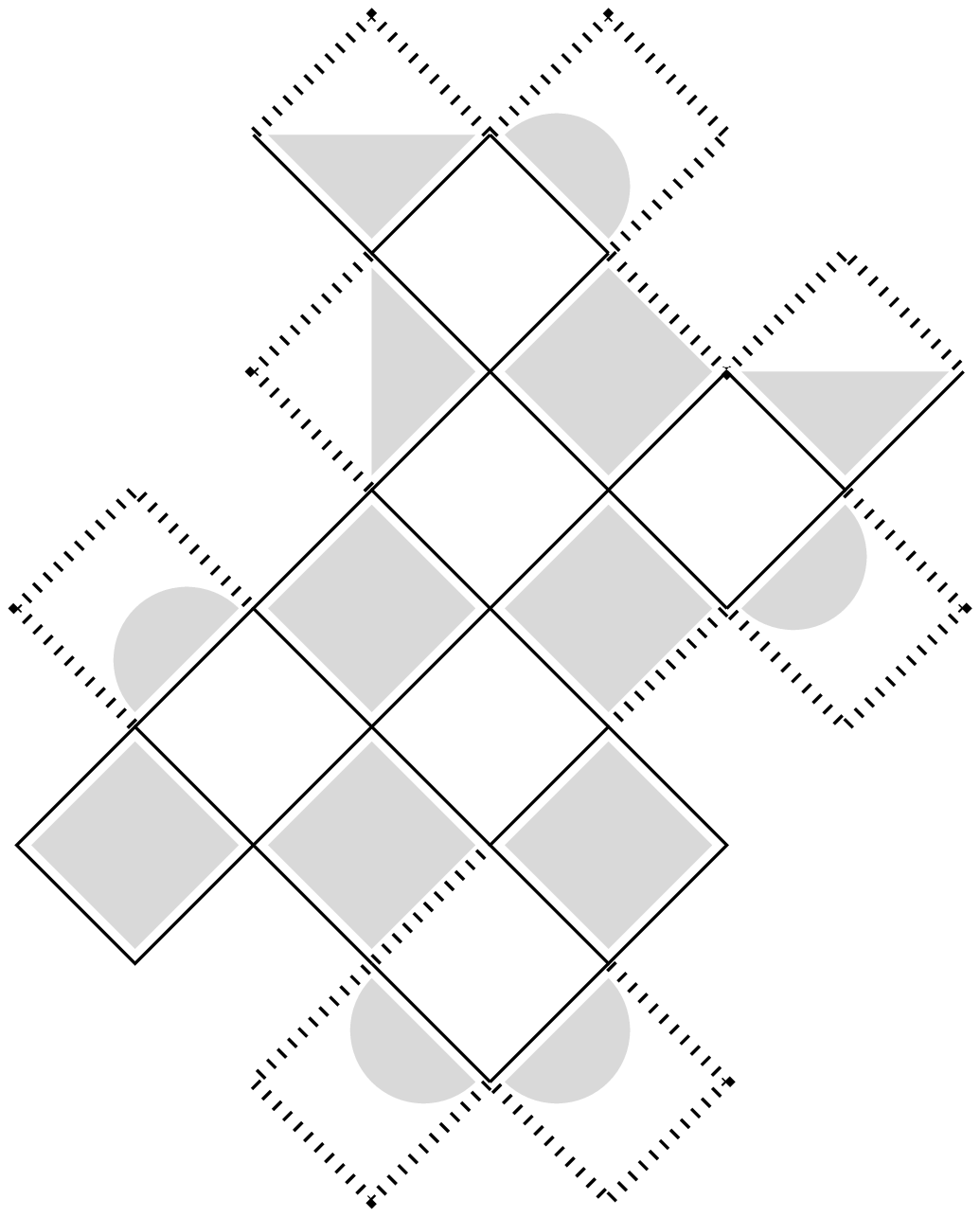}}
\centerline{{\smc Figure~2.5}{\rm (a). }}

\medskip

\twoline{\mypic{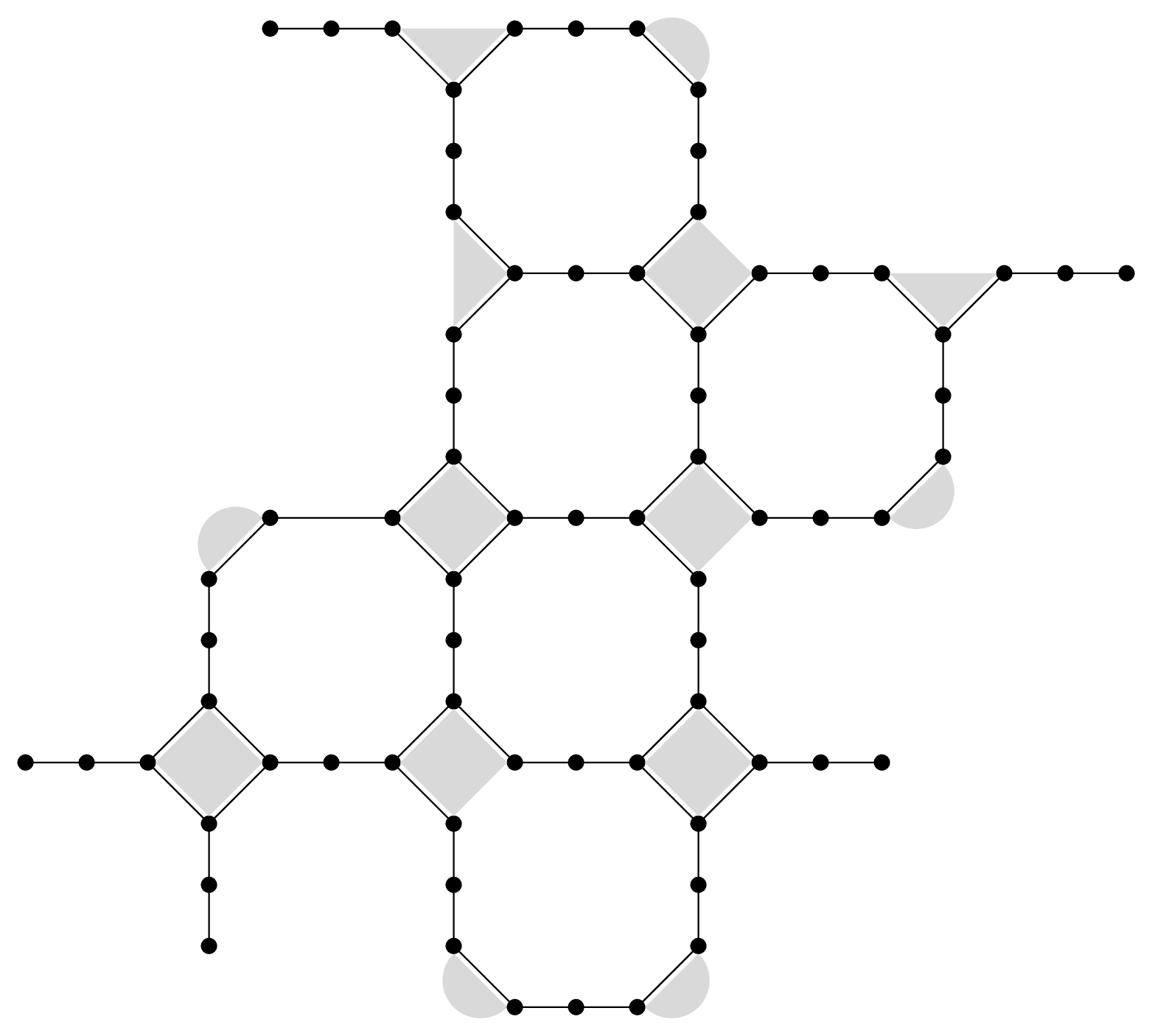}}{\mypic{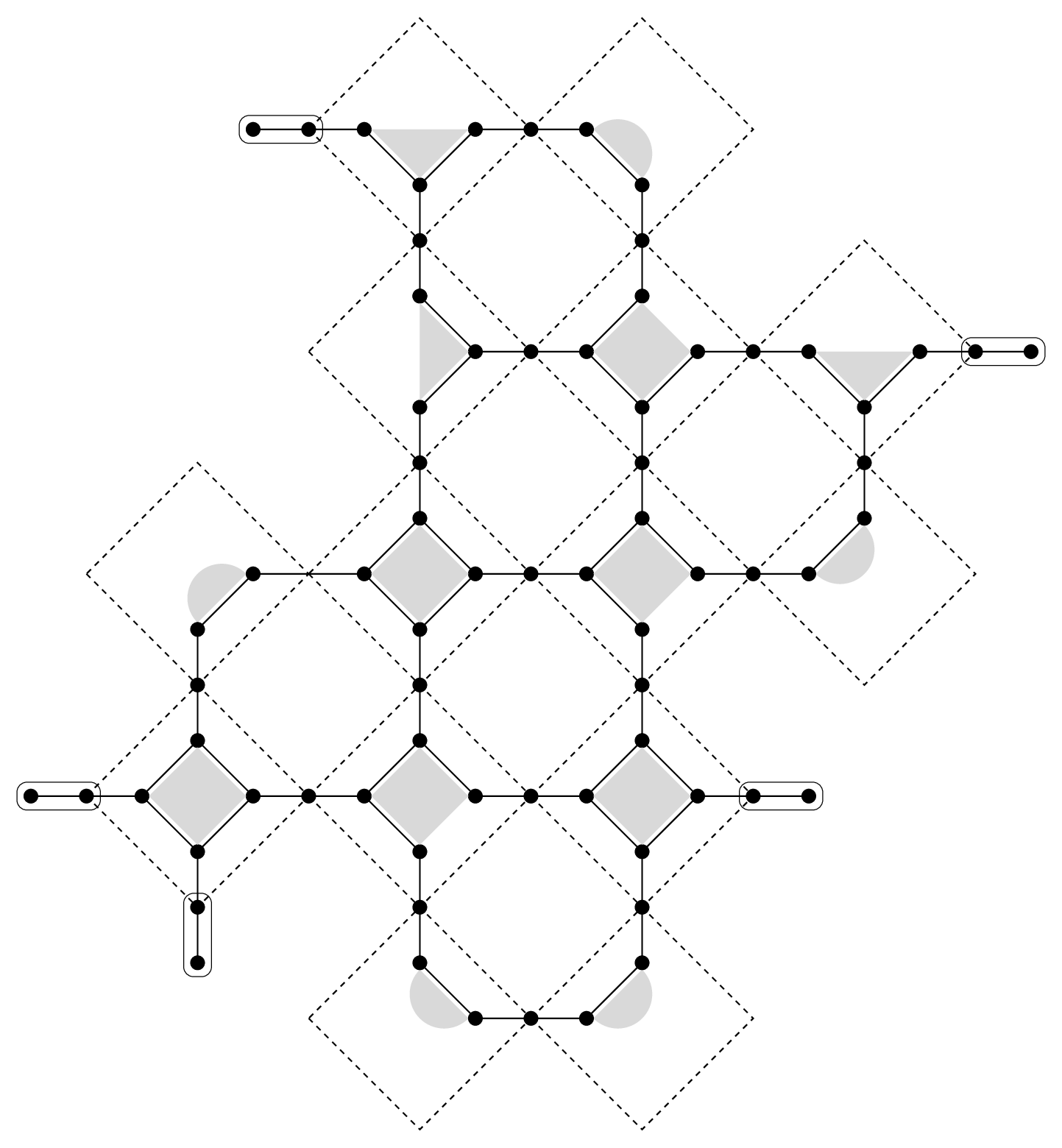}}
\twoline{Figure~2.5{\rm (b). }}
{Figure~2.5{\rm (c). }}
\endinsert

In the graph $H_1$ there is an opportunity to apply Lemma 2.5 around each whole cell, and an
opportunity to apply either part (a) or (b) of Lemma 2.6 around every partial cell. Let $H_2$ be
the graph resulting from $H_1$ after performing all local replacements prescribed by these Lemmas
(Figure 2.5(c) illustrates $H_2$ for the example of Figures 2.5(a) and (b); $H_2$ is obtained by
including the dotted 4-cycles and removing all vertices and edges they enclose).

Corresponding to each extremal vertex of $G$ contained in $H$ we have a forced edge in $H_2$
(these are circled in Figure 2.5(c)). Corresponding to each extremal vertex of $G$ {\it not}
contained in $H$, we have a new vertex in $H_2$, created by application of Lemma 2.6(a) or (b).
Therefore, the graph obtained from $H_2$ by removing all forced edges 
(which by construction have weight 1) is
precisely the complement $H'$ (compare Figures 2.5(c) and 2.1(d)). By Lemmas 2.5 and 2.6, the
resulting weight on $H'$ is precisely $\wt'$. \endpf

\flushpar
{\smc Remark 2.4.} Any subgraph $H$ of a weighted cellular graph $G$ has a 
cellular completion. Indeed, restricting $G$ to the cells containing at least one edge of 
$H$, redefining as being 0 the weight of all edges in $E(G)\setminus E(H)$ having both 
endpoints in $H$, and finally ``detaching'' the remaining cells of $G$ at vertices in 
$V(G)\setminus V(H)$ (so that they become extremal),
we obtain a cellular completion of $H$. Therefore, Theorem 2.3 applies to
any subgraph of a weighted cellular graph.

%
%

\smallpagebreak
Let $G$ be a weighted graph containing a subgraph isomorphic to the graph $K$ shown in the
left part of Figure 2.6 (the labels indicate weights; unlabeled edges have weight 1).
Suppose in addition that the four inner vertices of $K$ have no neighbors outside $K$. 
Let $\s(G)$ be the graph obtained from $G$ (by ``surgery'') by replacing $K$ by the graph 
$\bar{K}$ shown on right in Figure 2.6 (dashed lines indicate new edges, weighted as shown).


The following result is a generalization due to Propp of the ``urban renewal'' 
trick first
observed by Kuperberg (which corresponds to the case when all weights are 1). It follows by
analyzing the restrictions of matchings of $G$ and $\s(G)$ to $K$ and $\bar{K}$, respectively
(see \cite{9}).

\proclaim{Lemma 2.5} 
$\M(G)=(xz+yw)\M(\s(G)).$
\endproclaim

\topinsert
\centerline{\mypic{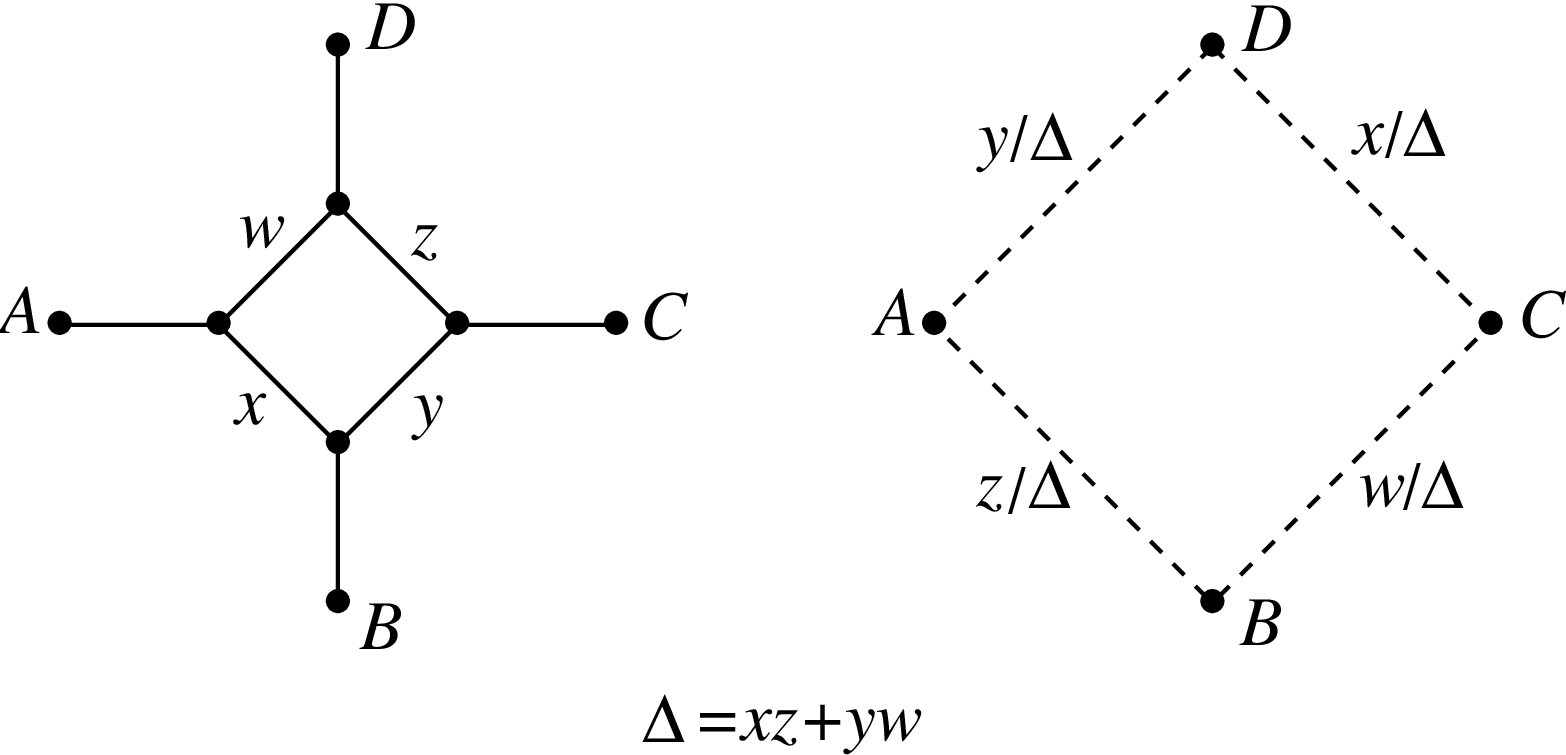}}
\centerline{{\smc Figure~2.6.} {\rm }}

\centerline{\mypic{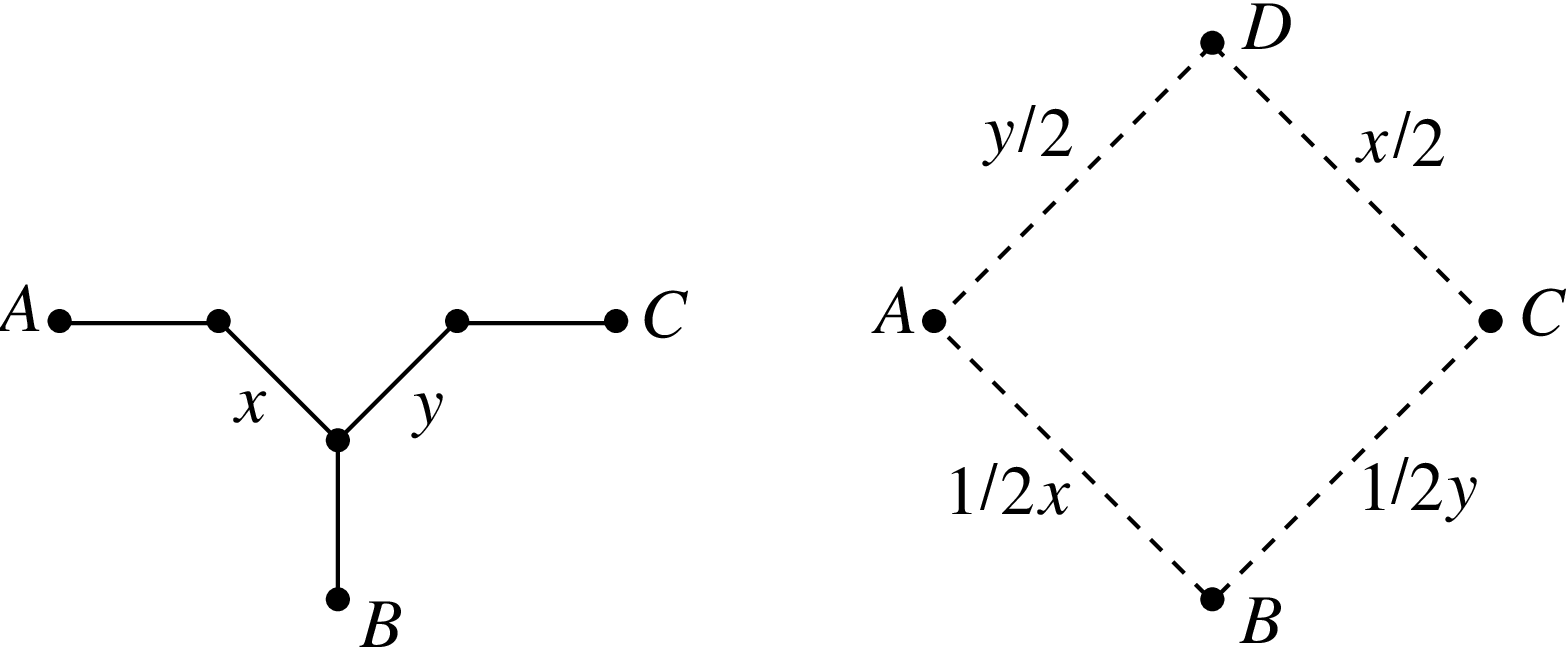}}
\centerline{{\smc Figure~2.7.} {\rm }}

\medskip
\centerline{\mypic{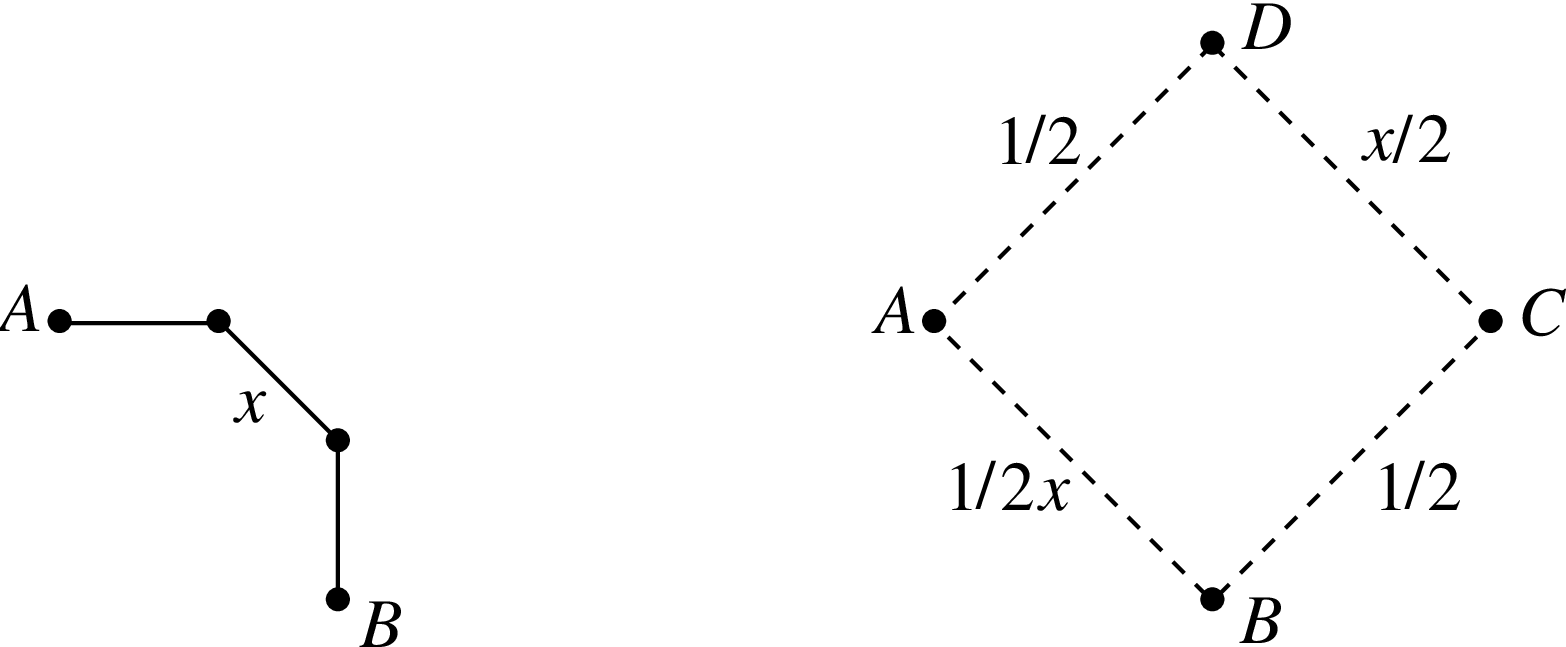}}
\centerline{{\smc Figure~2.8.} {\rm }}
\endinsert

The next result appears to be new. Together with Lemma 2.5, it is an essential ingredient in 
the proof of Theorem~2.3.

\proclaim{Lemma 2.6} $(a)$ Consider the above local replacement operation when $K$ and 
$\bar{K}$ are the graphs shown in Figure 2.7, with the indicated weights $($in particular, 
$\s(G)$ has one vertex, $D$, that was not a vertex of $G$, which is incident only to $A$ and $C)$. We have
$$
\M(G)=2\M(\s(G)).
$$

$(b)$ The statement of part $(a)$ is also true when $K$ and $\bar{K}$ are the graphs
indicated in Figure 2.8 $($in this case $\s(G)$ has two vertices, $C$ and $D$, not belonging to $G$; 
they are adjacent only to one another and to $B$ and $A$, respectively$)$.
\endproclaim

\pf To prove part (a), partition the perfect matchings of $G$ into three classes: (1) 
those containing the edge of $K$ weighted $x$ and the edge of $K$ incident to $C$, 
(2) those containing the edge of $K$ weighted $y$ and the edge of $K$ incident to $A$,
and (3) those containing the three edges of $K$ incident to $A$, $B$ or $C$. Partition the
perfect matchings of $\s(G)$ also in three classes: (1) those containing edge $CD$ but not $AB$, 
(2) those containing $AD$ but not $BC$, and (3) those containing two opposite edges of the 
4-cycle $ABCD$ (this is indeed a partition, since $D$ is incident only to $A$ and $C$). 

Given a perfect matching $\mu$ of $G$, 
construct a perfect matching $\mu'$ of $\s(G)$ by discarding 
the edges of $\mu$ with both endpoints in $K$, and including edges of $\bar{K}$ as indicated by the 
correspondence of like-numbered classes in the partitions from the previous paragraph. It is easy to 
check that $\mu\to\mu'$ is a bijection and that $\wt(\mu)=2\wt(\mu')$. This proves part~(a).

Part (b) is proved similarly. \endpf

When $H$ is an Aztec diamond we obtain the following important special case due to
Propp (called Generalized Domino-Shuffling in \cite{9}).

\proclaim{Corollary 2.7 (Reduction Theorem \cite{9})} 
$$\M(AD_n;\wt)=\M(AD_{n-1};\wt')\prod_{c\in\Cal C}\Delta(c).$$
\endproclaim

Provided all cell-factors are nonzero and keep being so after each application of the
Reduction Theorem,
Corollary 2.7 can be applied successively---with the weight changing at every
step---until we get down to the Aztec diamond of order 0, whose matching generating
function is equal to 1. This provides in particular an algorithm for finding the matching
generating function of the original weighted Aztec diamond. In certain situations, 
the weight pattern repeats after a small number of
successive applications. In these cases the matching generating function turns out to
be given by
perfect (or near-perfect) powers. Several instances of this phenomenon are presented
in the following sections.

\mysec{3. Aztec dungeons}


\topinsert
\centerline{\mypic{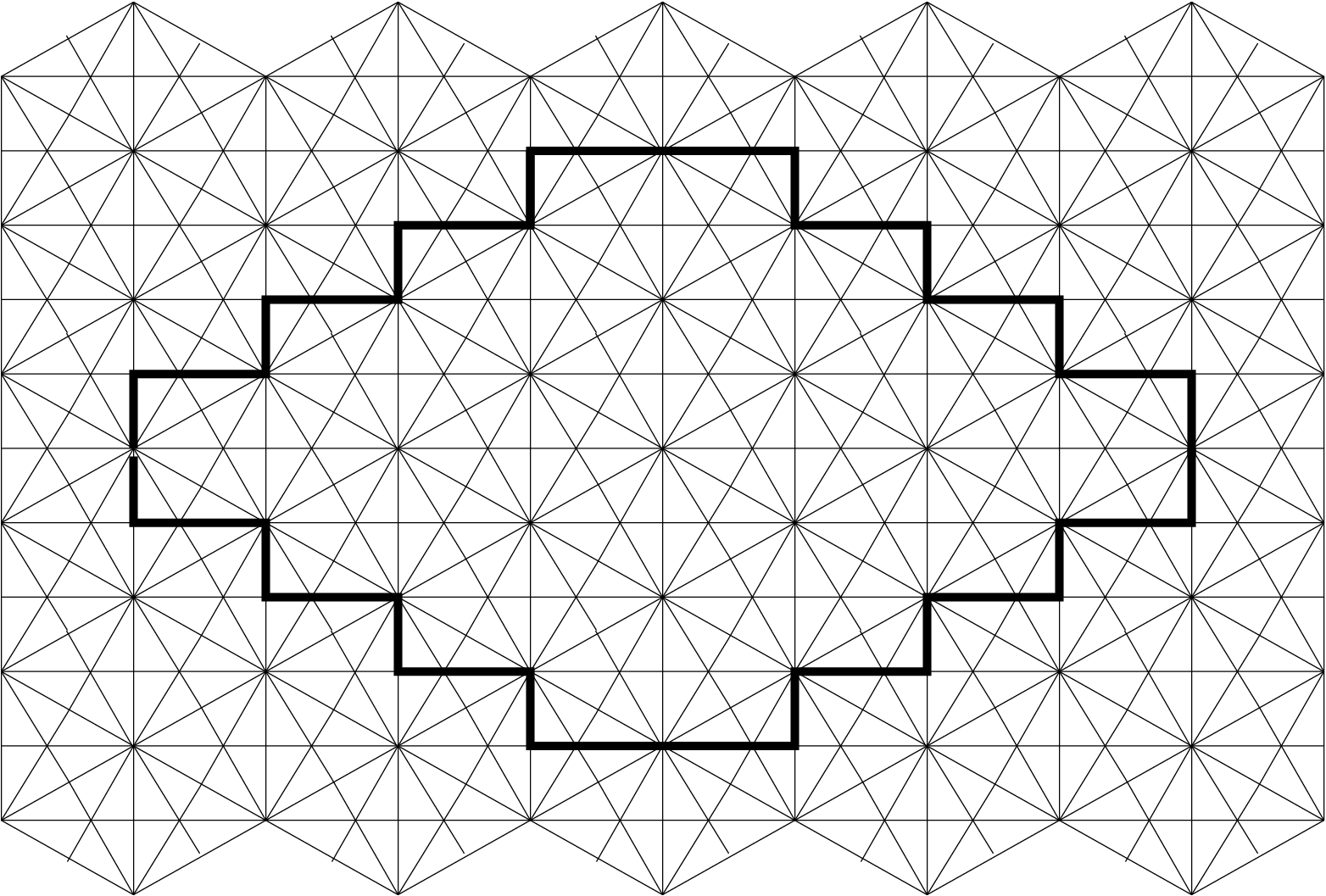}}
\centerline{{\smc Figure~3.1.} {\rm The lattice $G_2$ and the Aztec dungeon $D_4$.}}
\endinsert

The plane lattice corresponding to the affine Coxeter group $G_2$ is obtained from the
equilateral triangular lattice by drawing in all altitudes (see Figure 3.1). Consider
the family of all lattice lines parallel to a given lattice line direction $\ell$, 
and also the family of all lattice lines perpendicular to $\ell$. The union of these
two families forms a sublattice isomorphic to the square lattice. 

Tracing out what would correspond to an Aztec diamond region of order $n$
on this sublattice we obtain the {\it Aztec dungeon} $D_n$ of order $n$ (see 
Figure 3.1 for an example). Propp conjectured that the number of 
tilings of an Aztec dungeon is either a power of 13 or twice a power of 13. In
this section we prove two different generalizations of this conjecture (see Theorems
3.1 and 3.8).

There are three possible shapes of tiles on the lattice $G_2$: an equilateral
triangle, an obtuse triangle and a kite-shaped tile. It is natural to keep
track of them by assigning them distinct weights. Since the total number of tiles in a 
tiling is
determined by the order of the Aztec dungeon, we may assume without loss of
generality that one of our weights is equal to 1. Assign therefore our three
types of tiles weights $y$, $x$ and $1$, respectively. 
Denote by $\wt$ the weight obtained this way on $D_n$.

\proclaim{Theorem 3.1} Let $P$ denote the irreducible polynomial
$P:=x^6+3x^4y^2+3x^2y^4+y^6+2x^3+2xy^2+1$.
The tiling generating functions of the first six Aztec dungeons $D_n$ are given by
$$
\align
\T(D_0;\wt)&=1\\
\T(D_1;\wt)&=x^2+y^2\\
\T(D_2;\wt)&=x^2y^2P\\
\T(D_3;\wt)&=x^6y^6P^3\\
\T(D_4;\wt)&=x^{10}y^{14}(x^2+y^2)P^5\\
\T(D_5;\wt)&=x^{16}y^{24}P^8,
\endalign
$$
and for $n\geq 5$ we have the recurrence
$$\T(D_{n+1};\wt)=x^{8n-16}y^{16n-20}P^{4n-8}
\T(D_{n-5};\wt).\tag3.1$$
\endproclaim

\flushpar
Propp's original conjecture corresponds to the following special case.

\proclaim{Corollary 3.2} The number of tilings of the Aztec dungeon $D_n$ is given by the
equalities $\T(D_0)=1$, $\T(D_1)=2$, $\T(D_2)=13$, $\T(D_3)=13^3$,
$\T(D_4)=2\cdot13^5$, $\T(D_5)=13^8$, and for $n\geq5$ by the recurrence
$$\T(D_{n+1})=13^{4n-8}\T(D_{n-5}).$$
\endproclaim

\pf Set $x=y=1$ in Theorem 3.1. \endpf

The proof of Theorem 3.1 will follow from the following three preliminary results. The
first of them expresses $\M(D_n;\wt)$ in terms of a certain weighted count of
perfect matchings of the Aztec diamond graph $AD_{2n-2}$.

Let $A$ be a given $k\times l$ matrix with $k$ and $l$ even. 
The centers of the edges of the Aztec diamond 
graph $AD_n$ form a $2n\times2n$ array. Place a copy of $A$ in the upper left corner of 
this 
array and fill in the rest of the array periodically with period $A$ (i.e., translate
$A$ to the right in the array $l$ units at a time and down
in the array $k$ units at a time; if $2n$ is not a multiple of $k$ or $l$ 
some of these translates will fit only partially in the array).

\proclaim{Definition 3.3} Define the weight $\wt_A$ on the edges of $AD_n$ by
assigning each edge the corresponding entry of $A$ in the array described above.
\endproclaim

\proclaim{Lemma 3.4} Let $N$ be the matrix
$$
N=\left[ \matrix
{y\over x^2+y^2}&y&x
&{x\over x^2+y^2}\\
y&0&1&x\\
x&1&0&y\\
{x\over x^2+y^2}&x&y&{y\over x^2+y^2}
\endmatrix\right].\tag3.2
$$
We have
$$
\T(D_n;\wt)=y^{4n}(x^2+y^2)^{n^2}\M(AD_{2n-2};\wt_N).\tag3.3
$$
\endproclaim

\pf Construct the graph dual to the region $D_n$, i.e. the graph whose vertices are
the elementary regions of $D_n$ and whose edges connect precisely those elementary
regions that share an edge (this is illustrated in Figure 3.2 for $n=4$). Weight
each edge by the weight of the corresponding tile.

\topinsert
\centerline{\mypic{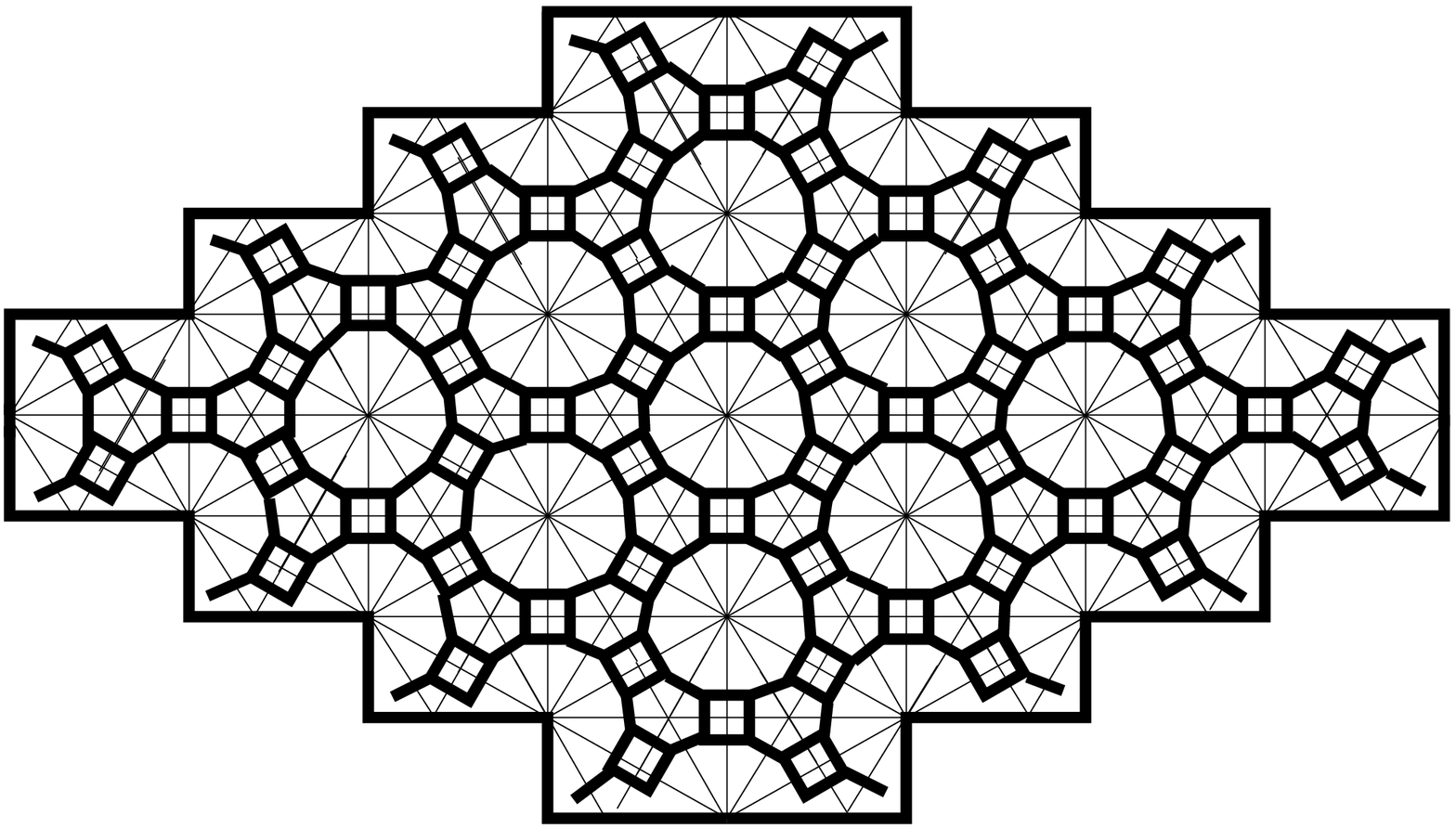}}
\centerline{{\smc Figure~3.2.} {\rm The Aztec dungeon $D_4$ and its dual graph.}}

\medskip
\centerline{\mypic{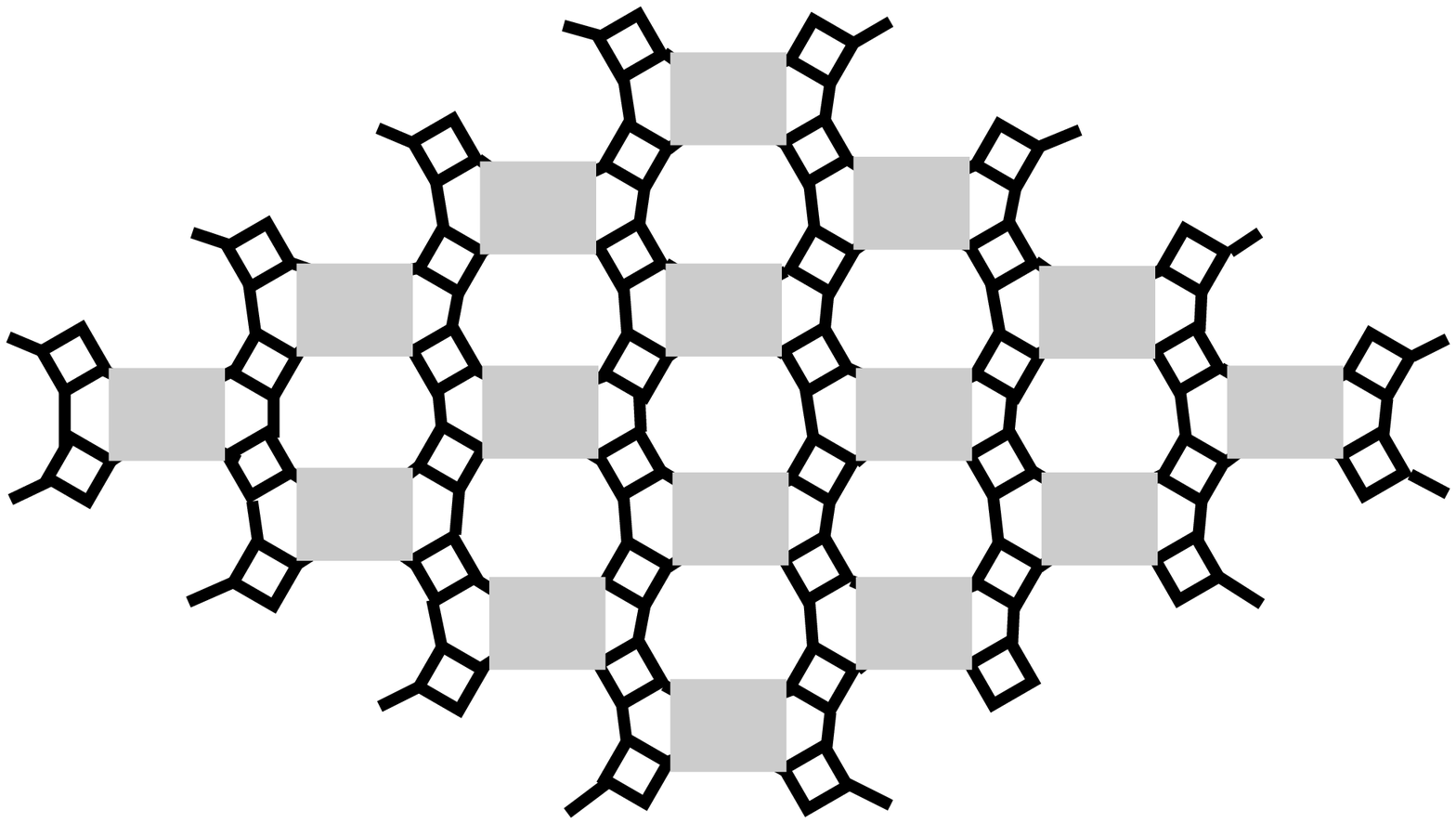}}
\centerline{{\smc Figure~3.3.} {\rm Applying Lemma 2.5 to the graph dual to $D_4$.}}
\endinsert

Clearly, the tilings of the region $D_n$ can be identified with the perfect matchings
of its dual graph. Therefore, the tiling generating function $\T(D_n;\wt)$ is the
same as the matching generating function of the dual graph of $D_n$.

This dual graph contains many local configurations like the one described at the
beginning of Section 2, providing one with as many opportunities to apply Lemma 2.5.
By their geometric orientation, these local configurations can be grouped into three
families. The largest family contains $n^2$ members (these are indicated in Figure
3.3). Apply Lemma 2.5 at each of these $n^2$ places. The resulting weighted graph is
readily seen to be isomorphic with a weighted subgraph of the square lattice
(shown in Figure 3.4 for $n=4$). 

Furthermore, because of the vertices of degree one, the resulting subgraph of the 
square lattice will have some edges around the boundary that are forced to be
contained in all of its perfect matchings (the forced edges are shown in thick
lines in Figure 3.4). It is easy to see that all these forced edges
except $4n$ of them correspond to kite-shaped tiles, and have therefore weight 1. The
remaining $4n$ edges correspond to equilateral triangle tiles, thus having weight
$y$.

Removing all vertices (together with all incident edges) connected by forced edges we
are left with a weighted spanning subgraph of the Aztec diamond of order $2n-2$
(Figure 3.5 illustrates the case $n=4$; dotted edges have weight $1/2$). 
This can be regarded as the entire graph
$AD_{2n-2}$ by weighting all missing edges by 0. By Lemma 2.5, the obtained weight
comes out to be exactly $\wt_N$, with the matrix $N$ given by (3.2). Taking into 
account the contribution of the
$4n$ forced edges of weight $y$ and applying $n^2$ times Lemma 2.5 we obtain
(3.3). \endpf

\topinsert
\twoline{\mypic{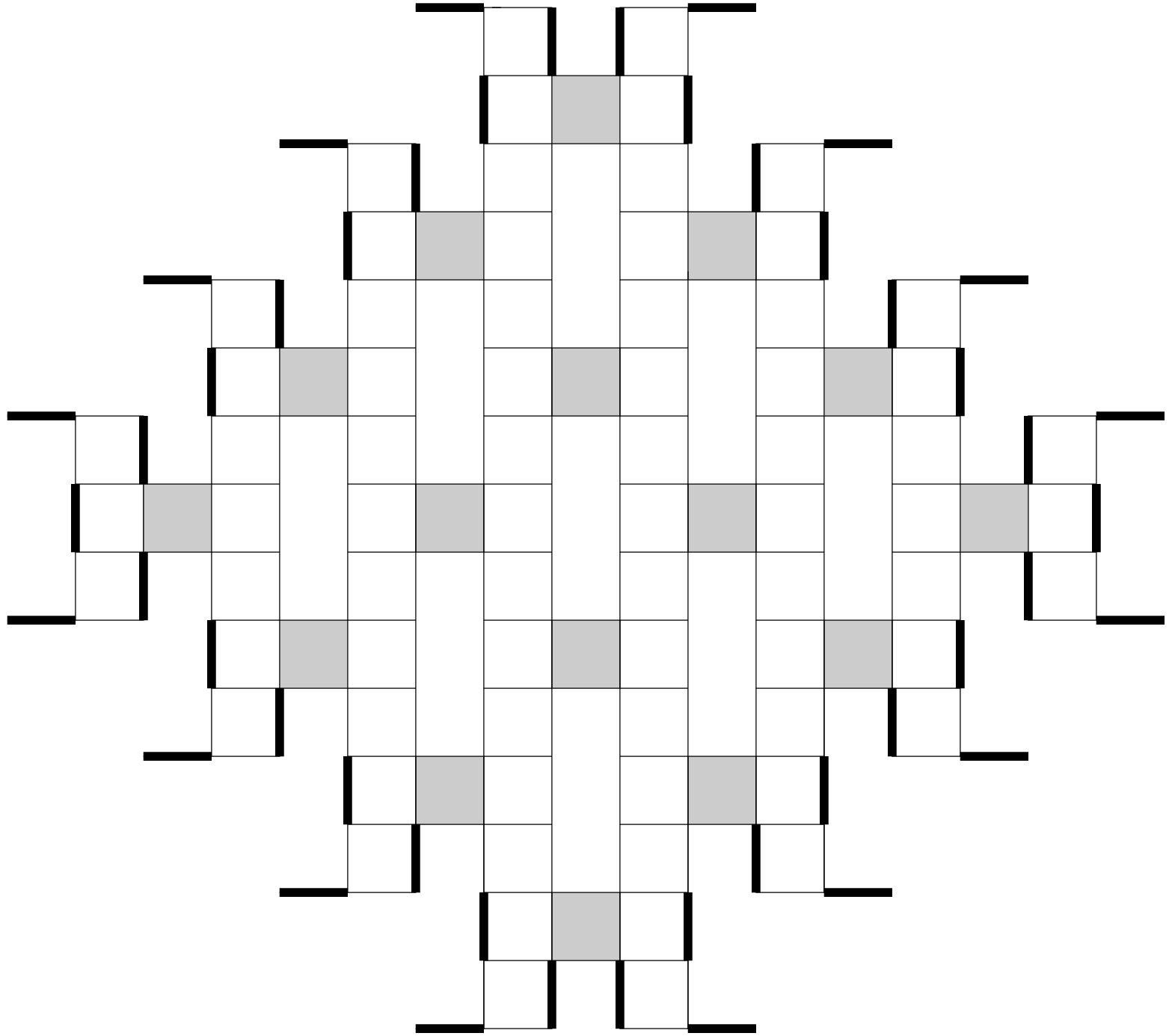}}{\mypic{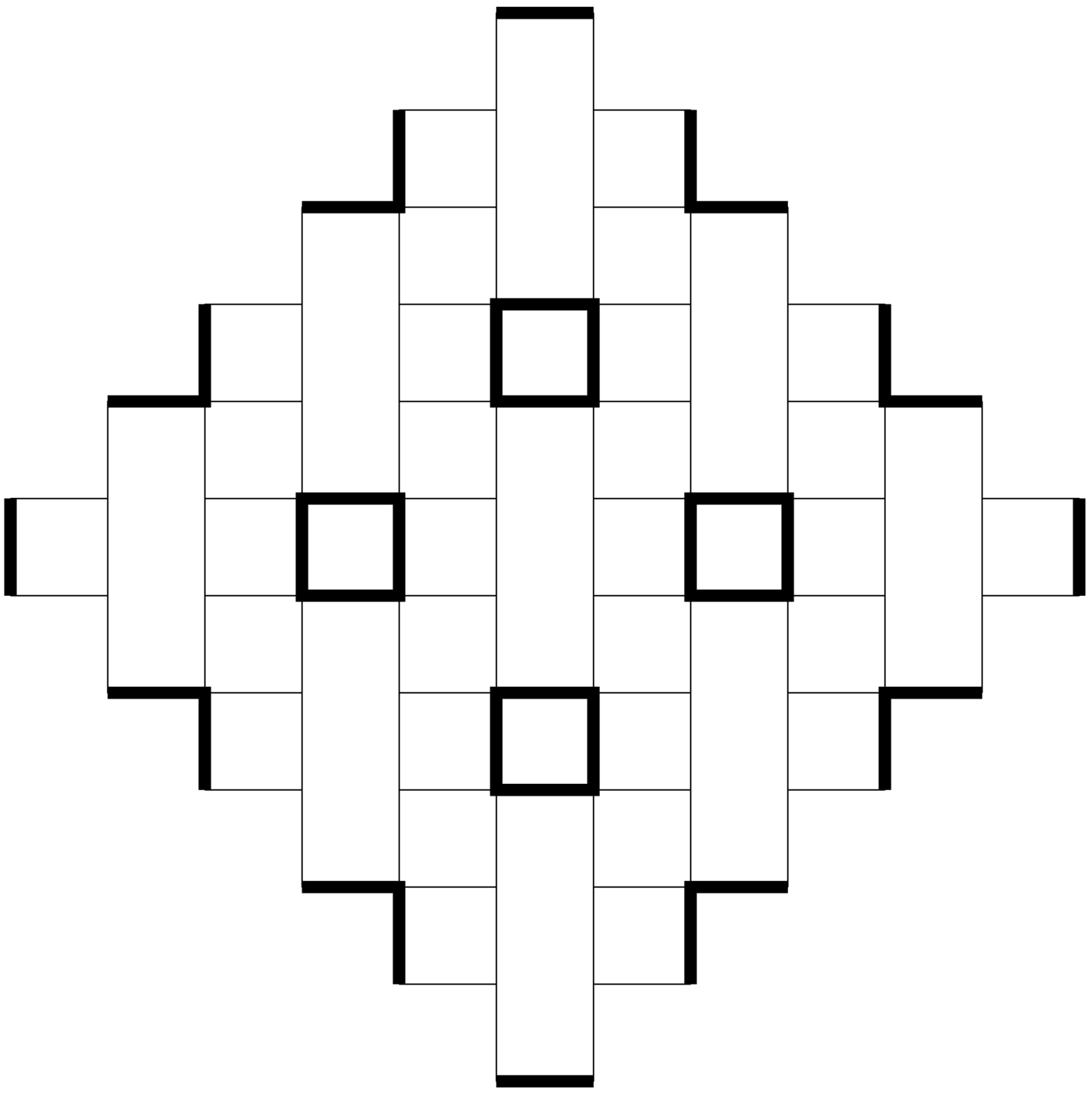}}
\twoline{Figure~3.4. {\rm A graph embedded in $\Z^2$.}}
{Figure~3.5. {\rm A weight on $AD_6$.}}
\endinsert

For a $k\times l$ matrix $A$ with $k$ and $l$ even define a new $k\times l$ matrix 
$\de(A)$ as follows. Divide matrix $A$ into $2\times2$ blocks
$$
\left[\matrix
x&w\\
y&z\\
\endmatrix\right]
$$
and assume $xz+yw\neq0$ for all such blocks. Replace each such block by
$$
\left[\matrix
z/(xz+yw)&y/(xz+yw)\\
w/(xz+yw)&x/(xz+yw)\\
\endmatrix\right]
$$
and denote the resulting $k\times l$ matrix by $B$. Define $\de(A)$ to be the 
$k\times l$ matrix obtained from $B$ by cyclically shifting its columns one unit up
and cyclically shifting the rows of the resulting matrix one unit left.

The next simple observation provides a very convenient way of keeping track of the 
evolution of periodic weights on Aztec diamonds when one applies the Reduction Theorem.

\proclaim{Lemma 3.5} Let $A$ be a $k\times l$ matrix with $k$ and $l$ even and
consider the weight $\wt_A$ it determines on $AD_n$ according to Definition 3.3. We
have
$$
\wt_A'=\wt_{\de(A)},
$$
where the weight on the left hand side is defined as in Definition 2.2.
\endproclaim

\pf This follows from Definition 2.2 and the above construction of $\de(A)$. The reason 
we need the cyclic shifts in the construction of $\de(A)$ is because the derived weight 
of Definition 2.2 is defined by viewing the Aztec diamond of order $n-1$ as embedded
concentrically into $AD_n$. \endpf

\smallpagebreak
Due to the location of the 0's in the matrix $N$ given by (3.2), all the edges in an
Aztec diamond that are weighted 0 under $\wt_N$ are parallel among themselves. By 
Definition 2.2, this property is preserved by the derived weight. In particular, in
any cell at least one pair of opposite edges is assigned nonzero weights. Therefore
all cell-factors are nonzero and the Reduction Theorem can be applied successively.

By Lemma 3.5, the successive weights that occur are the weights corresponding
via Definition 3.3 to the iterates $\de^{(i)}(N)$, for $i=1,2,3,\dotsc$. Therefore, if
one of these iterates would be the same as $N$---or indeed the same up to a
scalar multiplicative factor---then we would get a recurrence for the matching generating 
function on the 
right hand side of (3.3). By (3.3) this would then translate into a recurrence for
$\M(D_n;\wt)$ and would solve the problem of computing the latter.

The computation of the iterates $\de^{(i)}(N)$ can be done very easily with a
computer algebra package like Maple. The first few iterates don't look very promising,
but perseverance pays off: the twelfth(!) iterate turns out to be, up to a scalar 
multiple, exactly the same as $N$. More precisely, one obtains
$$
\de^{(12)}(N)=k_0N,\tag3.4
$$
with
$$
k_0=\frac{y^4(x^3+xy^2+1)^4(x^4+2x^2y^2+y^4+x)^4}
{(x^2+y^2)^4(x^6+3x^4y^2+3x^2y^4+y^6+2x^3+2xy^2+1)^4}.\tag3.5
$$
Since at each application of the Reduction Theorem the order of the resulting Aztec diamond
decreases one unit, it follows that 
$$
\M(AD_{2n};\wt_N)=k_1 \M(AD_{2n-12};\wt_N),\tag3.6
$$
where the constant $k_1$ is the product of all the cell-factors arising in the twelve
applications of the Reduction Theorem, multiplied by $k_1^{(2n-12)(2n-11)}$ (the latter factor
is due to (3.4) and to the fact that each perfect matching of $AD_n$ contains $n(n+1)$ edges).
To carry out the computation of $k_1$ by hand is a fairly strenuous enterprise, but with the 
assistance of Maple
it is quite easy. We obtain by (3.6) the following result.

\proclaim{Proposition 3.6} For $n\geq6$ we have
$$
\M(AD_{2n};\wt_N)=x^{8n-16}y^{16n-44}(x^2+y^2)^{24-12n}P^{4n-8}
\M(AD_{2n-12};\wt_N),\tag3.7
$$
where $P=x^6+3x^4y^2+3x^2y^4+y^6+2x^3+2xy^2+1$.\endpf
\endproclaim 

We are now ready to present the proof of the main result of this section.

\smallpagebreak
{\it Proof of Theorem 3.1.} By Lemma 3.4 we have
$$
\T(D_{n+1};\wt)=y^{4n+4}(x^2+y^2)^{(n+1)^2}\M(AD_{2n};\wt_N)
$$
and
$$
\T(D_{n-5};\wt)=y^{4n-20}(x^2+y^2)^{(n-5)^2}\M(AD_{2n-12};\wt_N).
$$
Together with (3.7) these two equalities imply (3.1). The stated values for the tiling
generating functions of the first six Aztec dungeons are easily checked using (3.3)
and applying repeatedly the Reduction Theorem to evaluate the weighted perfect
matching counts of the resulting Aztec diamonds. \endpf

\flushpar
{\smc Remark 3.7.} In fact, it turns out that we can further
distinguish between the two possible orientations of the equilateral triangle tile by
assigning them distinct weights in our enumeration. Indeed, keeping the weight
$y$ for the up-pointing equilateral triangle tiles and assigning weight
$z$ to the down-pointing ones, the same calculations as
above lead to the recurrence
$$
\align
\T(D_{n+1};\wt)=x^{8n-16}(yz)^{8n-10}
(x^6+3x^4yz+3x^2y^2z^2+y^3z^3+2x^3+2&xyz+1)^{4n-8}\\
\times&\T(D_{n-5};\wt),
\endalign
$$
for $n\geq5$ (where $\wt$ denotes now the new weight). Similarly, we obtain the values
$$
\align
\T(D_0;\wt)&=1\\
\T(D_1;\wt)&=x^2+yz\\
\T(D_2;\wt)&=x^2yz Q\\
\T(D_3;\wt)&=x^6(yz)^3Q^3\\
\T(D_4;\wt)&=x^{10}(yz)^7(x^2+yz)Q^5\\
\T(D_5;\wt)&=x^{16}(yz)^{12}Q^8,
\endalign
$$
where
$Q:=x^6+3x^4yz+3x^2y^2z^2+y^3z^3+2x^3+2xyz+1$.

One observes that all above polynomials have the property that the exponents of
$y$ and $z$ in each of their monomials agree. This implies that the same
property holds for each monomial in the expansion of $\T(D_{n};\wt)$, and therefore
each tiling of $D_n$ contains the same number of up-pointing and down-pointing
equilateral triangles.

However, the ultimate such
generalization that would keep track also of the orientations of the other two kinds
of tiles (six possibilities for each) does not seem to lead to a weighted count that
is expressible by a simple product formula. Indeed, as the Reduction Theorem is
applied successively larger and larger irreducible polynomials occur in the
factorization of the resulting $\Delta$-factors and their form does not seem to fit a
simple pattern.

\smallpagebreak
The key fact that allowed our proof of Theorem 3.1 is the fact that the weight function
obtained by successive applications of the Reduction Theorem repeated---up to a
scalar multiple---after a finite number of applications (in our case twelve).

As mentioned in Remark 3.7, this does not seem to happen if we weight each tile 
orientation by a different weight. In the following result we present a
weight function on Aztec diamonds that depends on eight free parameters which does
have the
property that the weight resulting at the twelfth iteration is a constant times the
original weight. The result originally conjectured by Propp corresponds, via the
argument that proved (3.3), to setting
parameters $a$ and $b$ equal to $1/2$ and the rest equal to 1. 

Consider the $4\times4$ matrix
$$
M=\left[ \matrix
a&d&d&a\\
e&0&g&e\\
f&h&0&f\\
b&c&c&b
\endmatrix\right].
$$
and the weight $\wt_M$ it determines on $AD_n$ by Definition 3.3.

\proclaim{Theorem 3.8} Let $R$ denote the irreducible polynomial
$R:=abgh+acfg+bdeh+2cdef$.
The matching generating functions of the first twelve Aztec diamonds weighted by 
$\wt_M$ are given by
$$
\align
\M(AD_{0};\wt_M)&=1\\
\M(AD_{1};\wt_M)&=de\\
\M(AD_{2};\wt_M)&=abR\\
\M(AD_{3};\wt_M)&=(de)^2abghR\\
\M(AD_{4};\wt_M)&=2(ab)^2cdefR^3\\
\M(AD_{5};\wt_M)&=2(de)^4(abgh)^2cfR^3\\
\M(AD_{6};\wt_M)&=2^3(ab)^4(cdef)^3ghR^5\\
\M(AD_{7};\wt_M)&=2^3(de)^7(abgh)^4(cf)^3R^5\\
\M(AD_{8};\wt_M)&=2^5(ab)^7(cdef)^5(gh)^3R^8\\
\M(AD_{9};\wt_M)&=2^5(de)^{10}(abgh)^7(cf)^5R^8\\
\M(AD_{10};\wt_M)&=2^8(ab)^{10}(cdef)^8(gh)^5R^{12}\\
\M(AD_{11};\wt_M)&=2^8(de)^{14}(abgh)^{10}(cf)^8R^{12},
\endalign
$$
and for higher orders we have the recurrences
$$
\M(AD_{2n};\wt_M)=2^{4n-12}(ab)^{4n-10}(cdef)^{4n-12}(gh)^{4n-16}R^{4n-8}\M(AD_{2n-12};\wt_M)
\tag3.8
$$
and
$$
\M(AD_{2n-1};\wt_M)=2^{4n-16}(de)^{4n-10}(abgh)^{4n-14}(cf)^{4n-16}R^{4n-12}
\M(AD_{2n-13};\wt_M).\tag3.9
$$
\endproclaim

\pf As we noted in the proof of Proposition 3.6, Lemma 3.5 implies that the successive weights
one obtains when applying the Reduction Theorem repeatedly are the weights corresponding to
the iterates $\de^{(i)}(M)$, for $i=1,2,3,\dotsc$. The twelfth iterate is again a constant
multiple of $M$. More precisely, we obtain
$$
\de^{(12)}(M)=k_2M,\tag3.10
$$
with
$$
k_2=\frac{1}{16}\frac{(ag+de)^2(bh+cf)^2(ag+2de)^2(bh+2cf)^2}
{(abgh+acfg+bdeh+2cdef)^4}.\tag3.11
$$
Since at each application of the Reduction Theorem the order of the resulting Aztec diamond
decreases one unit, it follows that $\M(AD_{n};\wt_M)$ differs from $\M(AD_{n-12};\wt_M)$
just by a multiplicative constant, equal to the product of all the cell-factors arising in 
the twelve applications of the Reduction Theorem, multiplied by $k_2^{(n-12)(n-11)}$ (since
$AD_n$ has $(n-12)(n-11)$ edges in each perfect matching). The explicit
value of this multiplicative constant depends on the parity of $n$. (In Proposition 3.6 we
were only interested in even orders because those are the ones that correspond to Aztec
dungeons.) 

With the assistance of Maple the multiplicative constant is easily determined. It comes
out to be the one shown in (3.8) for $n$ even, and the one in (3.9) for n odd. \endpf

\flushpar
{\smc Remark 3.9.} Recurrence (3.1) explains why there are no new prime factors (in
the polynomial ring $\Z[a,y]$) entering the factorizations of $\T(D_n;\wt)$ for $n\geq6$
besides the ones occurring for the first six values of $n$.
However, we do not have an explanation for the aesthetically pleasing and somewhat
mysterious fact that the only prime
factors occurring for $n\le5$ are the weight indeterminates themselves and the irreducible
polynomial $P$. Indeed, (3.5) illustrates two other irreducible polynomials that occur
as divisors of cell-factors arising when the Reduction Theorem is applied successively---they
specialize to 3 and 5, respectively, for $a=y=1$ (and they are the only other prime
divisors of cell-factors arising in this process). However, as indicated by the initial
values stated in Theorem 3.1, both these prime factors cancel out in the expressions giving
$\T(D_n;\wt)$ for $n=0,1,\dotsc,5$---and hence the only prime factor of the number of tilings
of $D_n$ is 13, except for the occurrence of a single factor of 2 for $n=1\,(\text{mod}\,3)$.

The same remark applies to recurrences (3.8) and (3.9) and the initial values stated in
Theorem 3.8, as (3.11) indicates four additional potential prime factors---$ag+de$, 
$bh+cf$, $ag+2de$ and $bh+2cf$ (again, they are the only other prime
divisors of cell-factors in this case)---which, however, cancel out in the expressions of 
$\M(AD_n;\wt_M)$ for $n\leq11$.

\smallpagebreak
There is another way to place the Aztec diamond contour on the sublattice of rectangles 
described at the beginning of this section. Indeed, on can place it so that the resulting region has no
forced tiles. Denote the region obtained this way by $E_n$; we will still call it an Aztec dungeon. 
Figure 3.6 illustrates the case $n=5$. 

Enumerating the tilings of $E_n$ turns out to be a distinct problem from the case of $D_n$. 
However, as indicated in the next result, the resulting number of tilings is strikingly similar 
to that corresponding to $D_n$. The enumeration of the tilings of $E_n$ 
follows by a suitable specialization of Theorem 3.8, which is thus
seen to be a common generalization of the enumeration of the two kinds of Aztec dungeons.

\topinsert
\centerline{\mypic{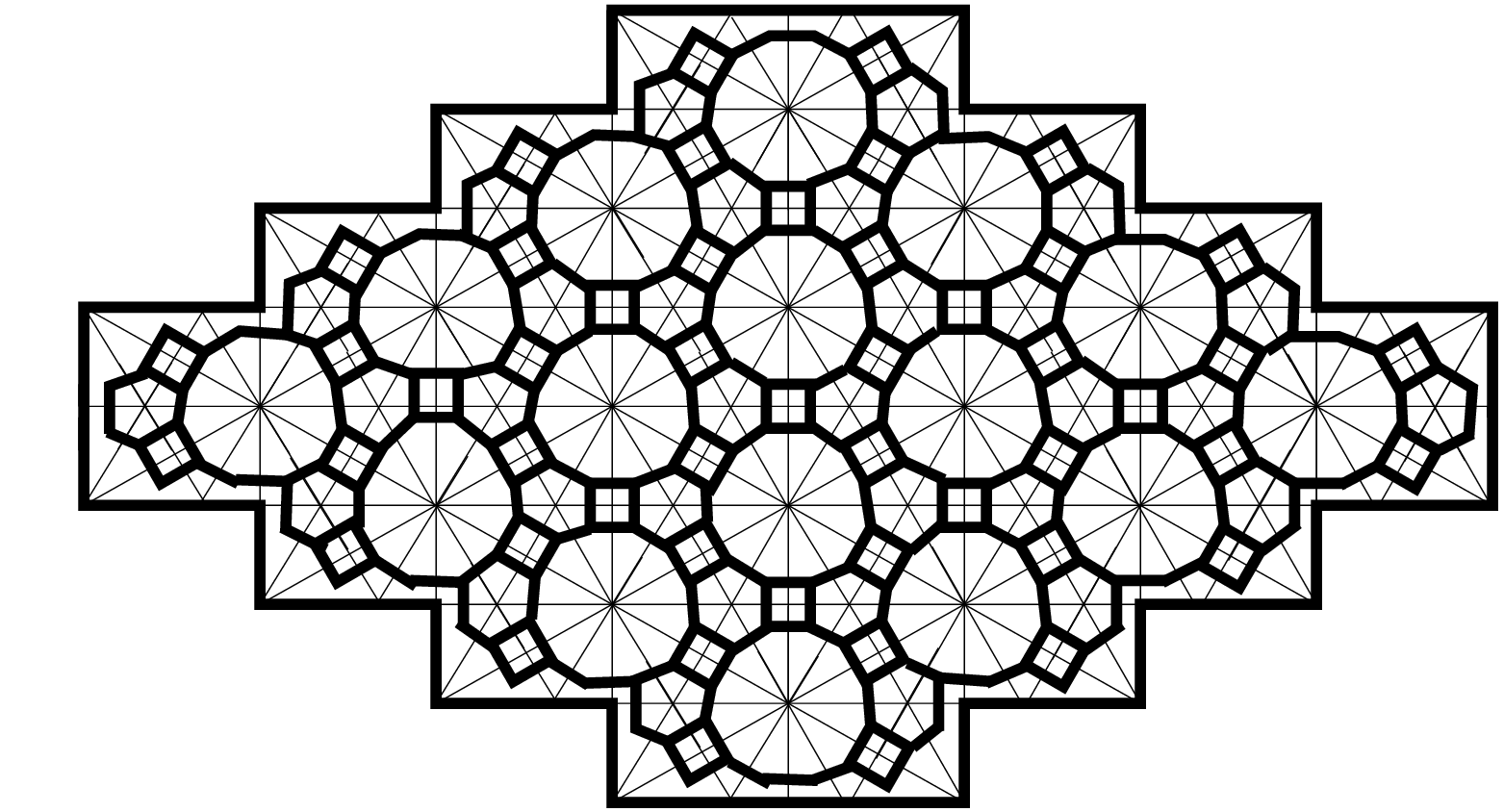}}
\centerline{{\smc Figure~3.6.} {\rm The Aztec dungeon $E_4$ and its dual graph.}}

\medskip
\centerline{\mypic{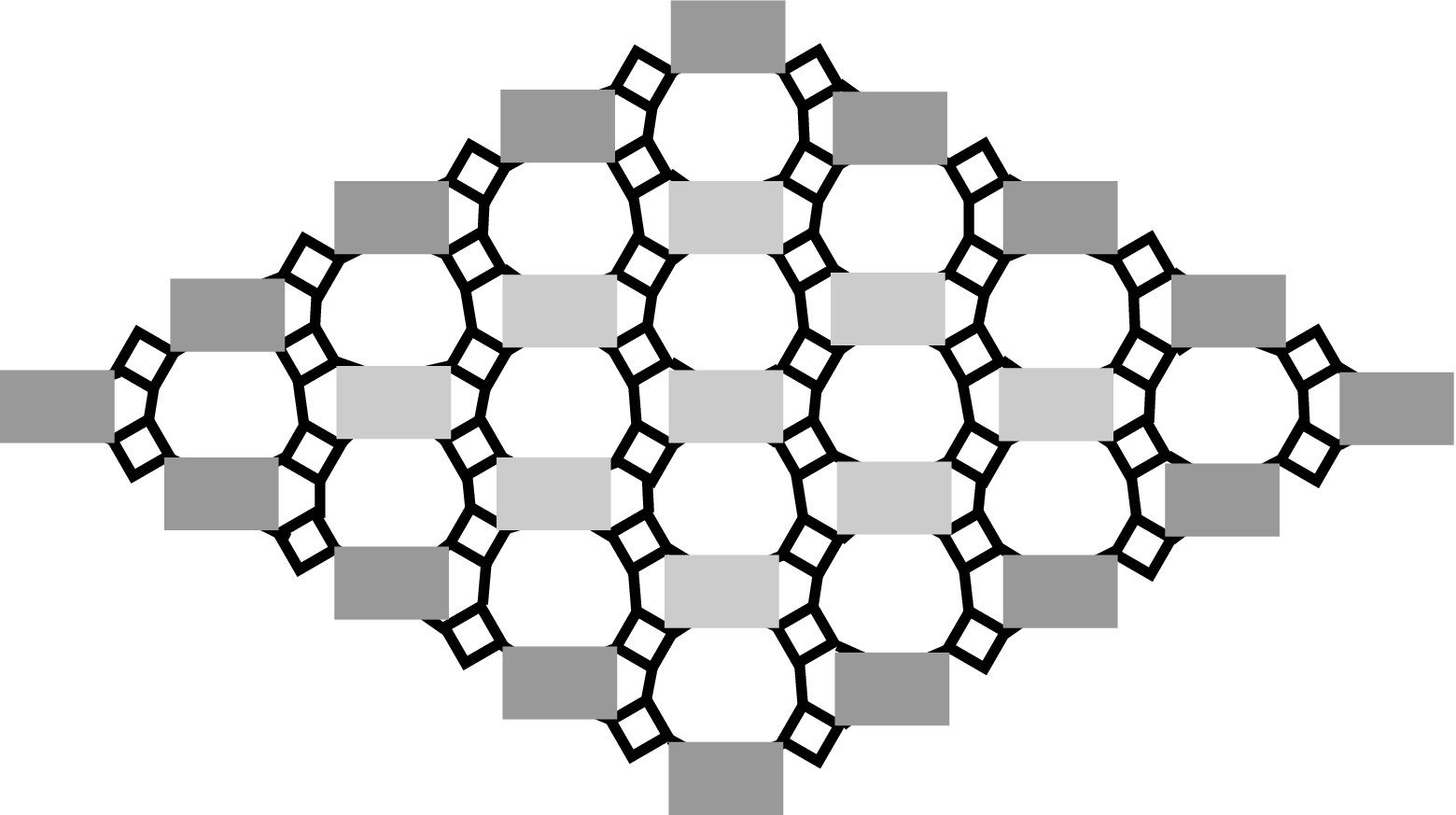}}
\centerline{{\smc Figure~3.7.} {\rm Applying Lemmas 2.5 and 2.6 to the graph dual to $E_4$.}}
\endinsert

\proclaim{Theorem 3.10} The number of tilings of the Aztec dungeon $E_n$ is given by the
equalities $\T(E_0)=1$, $\T(E_1)=2\cdot13$, $\T(E_2)=13^3$, $\T(E_3)=13^5$,
$\T(E_4)=2\cdot13^8$, $\T(E_5)=13^{12}$, and for $n\geq5$ by the recurrence
$$
\T(E_{n})=13^{4n-8}\T(E_{n-6}).
$$
\endproclaim

\pf Consider the dual graph of $E_n$ (Figure 3.6 illustrates the case $n=4$). 
Apply Lemma 2.5 at the $(n-1)^2$
places indicated by light shading in Figure 3.7, and Lemma 2.6 at the $4n$ places indicated in the same
figure by a darker shade. By these lemmas, with each such application the matching generating function 
changes by a factor of two. The resulting graph is isomorphic with $AD_{2n+1}$, and the resulting weight
is $\wt_B$ where 
$$
B=\left[ \matrix
1/2&1/2&1&1\\
1/2&1/2&1&1\\
1&1&0&1\\
1&1&1&0
\endmatrix\right].
$$
Therefore we obtain that
$$
\T(E_n)=2^{(n+1)^2}\M(AD_{2n+1};\wt_B).\tag3.12
$$
The Reduction Theorem and Lemma 3.5 yield, after combining all cell-factors of the thus weighted 
$AD_{2n+1}$, that
$$
\M(AD_{2n+1};\wt_B)=2^{n^2-1}\M(AD_{2n};\wt_{\de(B)}),\tag3.13
$$
where
$$
\de(B)=\left[ \matrix
1&1/2&1/2&1\\
1/2&0&1&1/2\\
1/2&1&0&1/2\\
1&1/2&1/2&1
\endmatrix\right].
$$
By specializing $a=b=g=h=1$ and $c=d=e=f=1/2$, (3.8) becomes
$$
M(AD_{2n};\wt_{\de(B)})=2^{60-24n}13^{4n-8}\M(AD_{2n-12};\wt_{\de(B)}).\tag3.14
$$
Using (3.12)--(3.14), $\M(E_n)$ is expressed as a constant times $\M(AD_{2n-12};\wt_{\de(B)})$. By (3.13),
the latter can be expressed in terms of $\M(AD_{2n-11};\wt_B)$, which is in turn, by (3.12), 
a constant multiple
of  $\M(E_{n-6})$. The multiplicative constant works out to be $13^{4n-8}$ and we obtain the recurrence in
the statement of the Theorem. The initial cases follow by specialization, using (3.12) and (3.13), from 
the initial cases of Theorem 3.8. $\square$

\mysec{4. Fortresses}

B-Y. Yang \cite{11} considered regions called fortresses on the lattice obtained from 
the grid lattice by drawing in all diagonals (Figure 4.1 illustrates the
fortress of order 4), and proved, using the permanent-determinant method \cite{6}, 
that the number of their tilings is either a power of 5 or 2 times a power of~5.

By applying Lemma 2.5 to the dual graph, enumerating the tilings of a fortress 
amounts to finding the matching generating function of an Aztec diamond weighted in a certain
pattern with weights 1 and 1/2. B-Y. Yang proved in fact a more general result \cite{11}, included here
as Corollary 4.3.

However, Yang's proof is quite complicated and does not explain why a nice factorization exists. 
This makes it inviting to look for a combinatorial proof. A combinatorial argument for the weight
relevant to fortress tilings can be found in \cite{2}. We present below a combinatorial proof of
the general case of Yang's theorem. Our proof also explains why such a nice factorization is obtained.

Consider the square matrix  
$$
A=\left[ \matrix
x_1&y_1&1/x_1&1/y_1&x_1&y_1&1/x_1&1/y_1&.&.&.\\
x_2&y_2&1/x_2&1/y_2&x_2&y_2&1/x_2&1/y_2&.&.&.\\
.&.&.&.&.&.&.&.\\
.&.&.&.&.&.&.&.\\
.&.&.&.&.&.&.&.\\  
x_{2n}&y_{2n}&1/x_{2n}&1/y_{2n}&x_{2n}&y_{2n}&1/x_{2n}&1/y_{2n}&.&.&.
\endmatrix\right].
$$
The weight $\wt_A$ it defines on $AD_{n}$ according to Definition 3.3. is determined by its submatrix
$$
T=\left[ \matrix
x_1&y_1\\
x_2&y_2\\
.&.\\
.&.\\
.&.\\  
x_{2n}&y_{2n}
\endmatrix\right].
$$
For simplicity of notation, rename the weight $\wt_A$ as $\y_T$.
Let $S_i:=(x_{2i-1}y_{2i-1}x_{2i}y_{2i})^{1/2}$, $i=1,\dotsc,n$. Define 
$\r$ to be the operator sending $T$ to the $(2n-2)\times2$ matrix
$$
\r(T)=\left[ \matrix
x_1/S_1&S_1/y_1\\
x_4/S_2&S_2/y_4\\
x_3/S_2&S_2/y_3\\
x_6/S_3&S_3/y_6\\
x_5/S_3&S_3/y_5\\
.&.\\
.&.\\
.&.\\  
x_{2n-2}/S_{n-1}&S_{n-1}/y_{2n-2}\\
x_{2n-3}/S_{n-1}&S_{n-1}/y_{2n-3}\\
x_{2n}/S_{n}&S_{n}/y_{2n}
\endmatrix\right].
$$

Yang's theorem (together with an explanation for why a nice factorization exists) will follow 
from our next two results.

\proclaim{Theorem 4.1} For $n\geq1$ we have
$$
\M(AD_{n};\y_T)=\M(AD_{n-1};\y_{\r(T)})\prod_{i=1}^n\Delta_i/S_i,\tag4.1
$$
for $n$ even, and
$$
\M(AD_{n};\y_T)=\M(AD_{n-1};\y_{\r(T)})\prod_{i=1}^n\Delta_i,\tag4.2
$$
for odd $n$, where $\Delta_i=x_{2i-1}y_{2i}+x_{2i}y_{2i-1}$, for $i=1,\dotsc,n$.
\endproclaim

Using the above result one can easily prove Yang's formula by induction, provided the formula has been
conjectured explicitly (it is not difficult to conjecture it based on computer calculations of the first
few cases). However, a conceptual way to obtain the formula is by using the following observation.

\proclaim{Lemma 4.2} For any $2n\times2$ matrix $T$ having rows $[x_i,y_i]$, $i=1,\dotsc,2n$, we have
$$
\r^{(4)}(T)=\left[ \matrix
x_5&y_5\\
x_6&y_6\\
x_7&y_7\\
.&.\\
.&.\\
.&.\\  
x_{2n-4}&y_{2n-4}
\endmatrix\right].
$$
\endproclaim

\pf This is checked directly from the definition of the operator $\r$ (assistance from a computer algebra
package like Maple makes this verification quite easy). \endpf

By Lemma 4.2 we are guaranteed that successive application of Theorem 4.1 will only generate a small
number of types of cell-factors and quantities $S_i$ (all types appear within four consecutive
applications of Theorem 4.1). We obtain the following result. 

\proclaim{Corollary 4.3} $(${\smc B-Y. Yang}$)$ The matching generating function $\M(AD_{n};\y_T)$
can be expressed as a simple product in the variables $x_i$ and $y_i$, $i=1,\dotsc,2n$.
\endproclaim 

In our proof of Theorem 4.1 we will employ the following simple observation.

\proclaim{Lemma 4.4} Partition the set of edges of $AD_n$ into $n+1$ classes as follows. View the set of
edges as forming a $2n\times 2n$ array. For each $i=1,\dotsc,n-1$ include the edges in rows $2i$ and $2i+1$ 
in a class. Let the edges in the first row and those in the $2n$-th row form two more classes. Then every
perfect matching of $AD_n$ contains precisely $n$ edges in each class.
\endproclaim

\pf Let $\mu$ be a perfect matching of $AD_n$. For each class, consider the set of $n$ vertices of $AD_n$
incident only to edges in that class. These guarantee there will be at least $n$ edges of $\mu$ in each
class. The statement follows from the fact that the number of edges in $\mu$ is $n(n+1)$. \endpf

{\it Proof of Theorem 4.1.} Suppose $n$ is even. For notational simplicity, we illustrate in detail the 
case $n=4$. The arguments generalize with no difficulty to general even $n$.

The original weight $y_T$ on the edges of our Aztec diamond $AD_4$ can be written out explicitely 
as $\wt_B$, where
$$
B=\left[ \matrix
x_1&y_1&1/x_1&1/y_1&x_1&y_1&1/x_1&1/y_1\\
x_2&y_2&1/x_2&1/y_2&x_2&y_2&1/x_2&1/y_2\\
x_3&y_3&1/x_3&1/y_3&x_3&y_3&1/x_3&1/y_3\\
x_4&y_4&1/x_4&1/y_4&x_4&y_4&1/x_4&1/y_4\\
x_5&y_5&1/x_5&1/y_5&x_5&y_5&1/x_5&1/y_5\\
x_6&y_6&1/x_6&1/y_6&x_6&y_6&1/x_6&1/y_6\\
x_7&y_7&1/x_7&1/y_7&x_7&y_7&1/x_7&1/y_7\\
x_8&y_8&1/x_8&1/y_8&x_8&y_8&1/x_8&1/y_8\\
\endmatrix\right].
$$
Apply the Reduction Theorem. The derived weight is $\wt_C$, with
$$
C=\left[ \matrix
x_1/\Delta_1&x_1x_2y_2/\Delta_1&x_2y_1y_2/\Delta_1&y_1/\Delta_1&x_1/\Delta_1&x_1x_2y_2/\Delta_1\\
x_4/\Delta_2&x_3x_4y_3/\Delta_2&x_3y_3y_4/\Delta_2&y_4/\Delta_2&x_4/\Delta_2&x_3x_4y_3/\Delta_2\\
x_3/\Delta_2&x_3x_4y_4/\Delta_2&x_4y_3y_4/\Delta_2&y_3/\Delta_2&x_3/\Delta_2&x_3x_4y_4/\Delta_2\\
x_6/\Delta_3&x_5x_6y_5/\Delta_3&x_5y_5y_6/\Delta_3&y_6/\Delta_3&x_6/\Delta_3&x_5x_6y_5/\Delta_3\\
x_5/\Delta_3&x_5x_6y_6/\Delta_3&x_6y_5y_6/\Delta_3&y_5/\Delta_3&x_5/\Delta_3&x_5x_6y_6/\Delta_3\\
x_8/\Delta_4&x_7x_8y_7/\Delta_4&x_7y_7y_8/\Delta_4&y_8/\Delta_4&x_8/\Delta_4&x_7x_8y_7/\Delta_4
\endmatrix\right]
$$
(where as defined in the statement of
the theorem, $\Delta_i=x_{2i-1}y_{2i}+x_{2i}y_{2i-1}$, for $i=1,\dotsc,4$). 
Taking into account the cell-factors,
the Reduction Theorem yields
$$\M(AD_4;\wt_B)=\prod_{i=1}^4(x_{2i-1}y_{2i}+x_{2i}y_{2i-1})^2
\left(\frac{1}{x_{2i-1}y_{2i}}+\frac{1}{x_{2i}y_{2i-1}}\right)^2\M(AD_{3};\wt_C).
\tag4.3
$$
Group the rows of $C$ into 4 blocks corresponding to the partition of edges of $AD_3$ from Lemma 4.4.
By Lemma 4.4 we can scale the blocks to obtain
$$
\M(AD_{3};\wt_C)=\prod_{i=1}^4\Delta_i^{-3}\M(AD_{3};\wt_D),\tag4.4
$$
where
$$
D=\left[ \matrix
x_1&x_1x_2y_2&x_2y_1y_2&y_1&x_1&x_1x_2y_2\\
x_4&x_3x_4y_3&x_3y_3y_4&y_4&x_4&x_3x_4y_3\\
x_3&x_3x_4y_4&x_4y_3y_4&y_3&x_3&x_3x_4y_4\\
x_6&x_5x_6y_5&x_5y_5y_6&y_6&x_6&x_5x_6y_5\\
x_5&x_5x_6y_6&x_6y_5y_6&y_5&x_5&x_5x_6y_6\\
x_8&x_7x_8y_7&x_7y_7y_8&y_8&x_8&x_7x_8y_7
\endmatrix\right].
$$
This matrix does not yet have the form of $B$, but can be brought to that form by 
further scaling of the blocks defined above. Indeed, what characterizes a matrix of
type $B=(b_{ij})$ is that the product $b_{ij}b_{i,j+2}$ is equal to 1 for all choices of 
indices for which both entries are defined. Matrix $D=(d_{ij})$ does not have this property,
but it does have the property that the products $d_{ij}d_{i,j+2}$ are constant---equal
to $x_{2i-1}y_{2i-1}x_{2i}y_{2i}$---within
each block. Therefore, Lemma 4.4 can be used to bring our matrix to the form of matrix
$B$. With $S_i=(x_{2i-1}y_{2i-1}x_{2i}y_{2i})^{1/2}$, $i=1,\dotsc,4$, we obtain by Lemma 4.4 
$$
\M(AD_{3};\wt_D)=\prod_{i=1}^4S_i^{3}\M(AD_{3};\wt_E),\tag4.5
$$
where 
$$
E=\left[ \matrix
x_1/S_1&S_1/y_1&S_1/x_1&y_1/S_1&x_1/S_1&S_1/y_1\\
x_4/S_2&S_2/y_4&S_2/x_4&y_4/S_2&x_4/S_2&S_2/y_4\\
x_3/S_2&S_2/y_3&S_2/x_3&y_3/S_2&x_3/S_2&S_2/y_3\\
x_6/S_3&S_3/y_6&S_3/x_6&y_6/S_3&x_6/S_3&S_3/y_6\\
x_5/S_3&S_3/y_5&S_3/x_5&y_5/S_3&x_5/S_3&S_3/y_5\\
x_8/S_4&S_4/y_8&S_4/x_8&y_8/S_4&x_8/S_4&S_4/y_8
\endmatrix\right].
$$
By (4.3)--(4.5) we obtain that
$$
\M(AD_4;\wt_B)=\prod_{i=1}^4\frac{\Delta_i}{S_i}\M(AD_{3};\wt_E).
$$
For general even $n$, the same argument proves (4.1).

The case of odd $n$ is treated similarly. The only difference is that now the 
exponents in the product in (4.3) are $(n+1)/2$ and $(n-1)/2$, as opposed to both 
being $n/2$ when $n$ is even. Relations (4.4) and (4.5) hold with no change for
odd $n$. One obtains the recurrence (4.2). \endpf

On the square lattice with all square diagonals drawn in there are two types of tiles: a
square tile and a triangular tile. The latter has four different possible orientations.

\topinsert
\twoline{\mypic{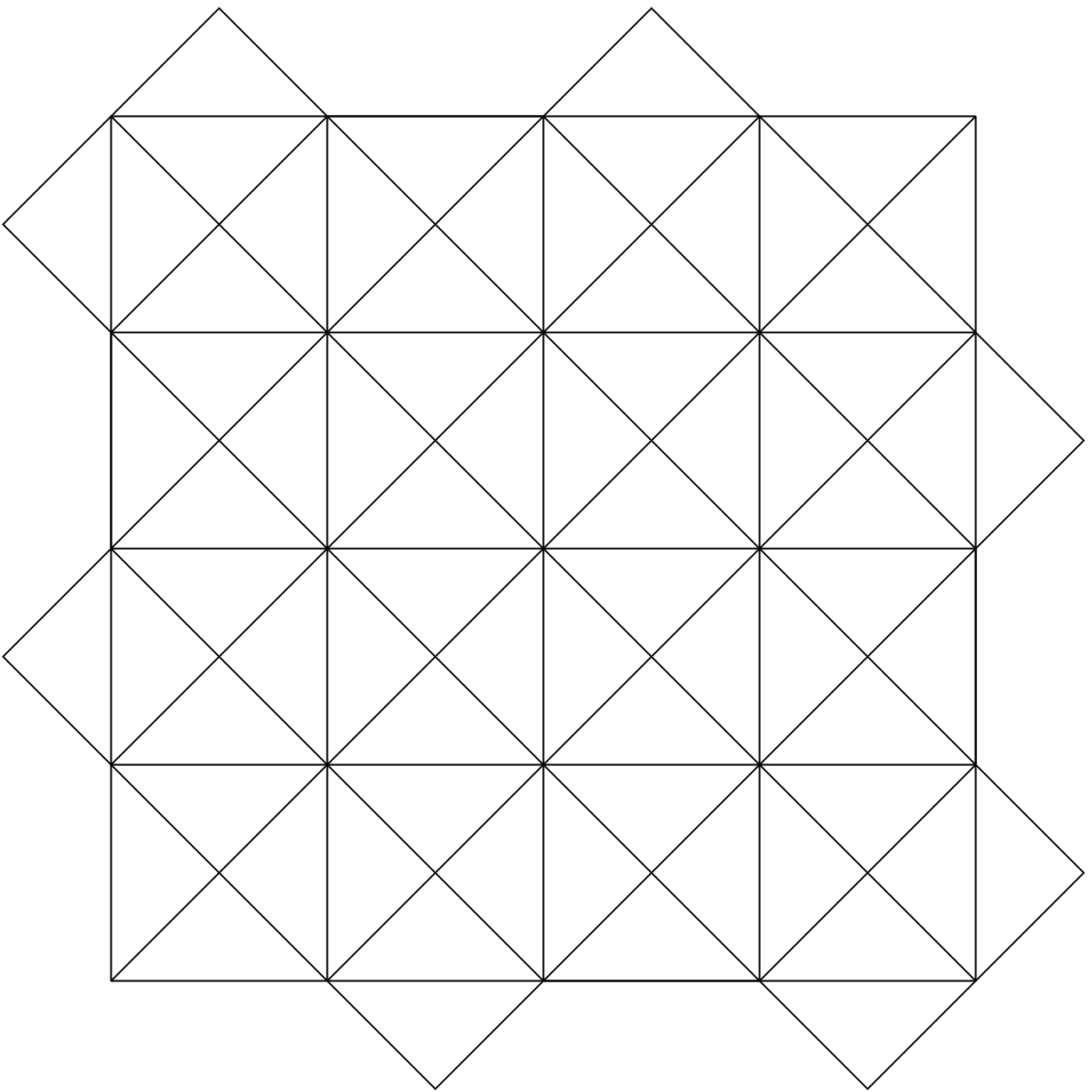}}{\mypic{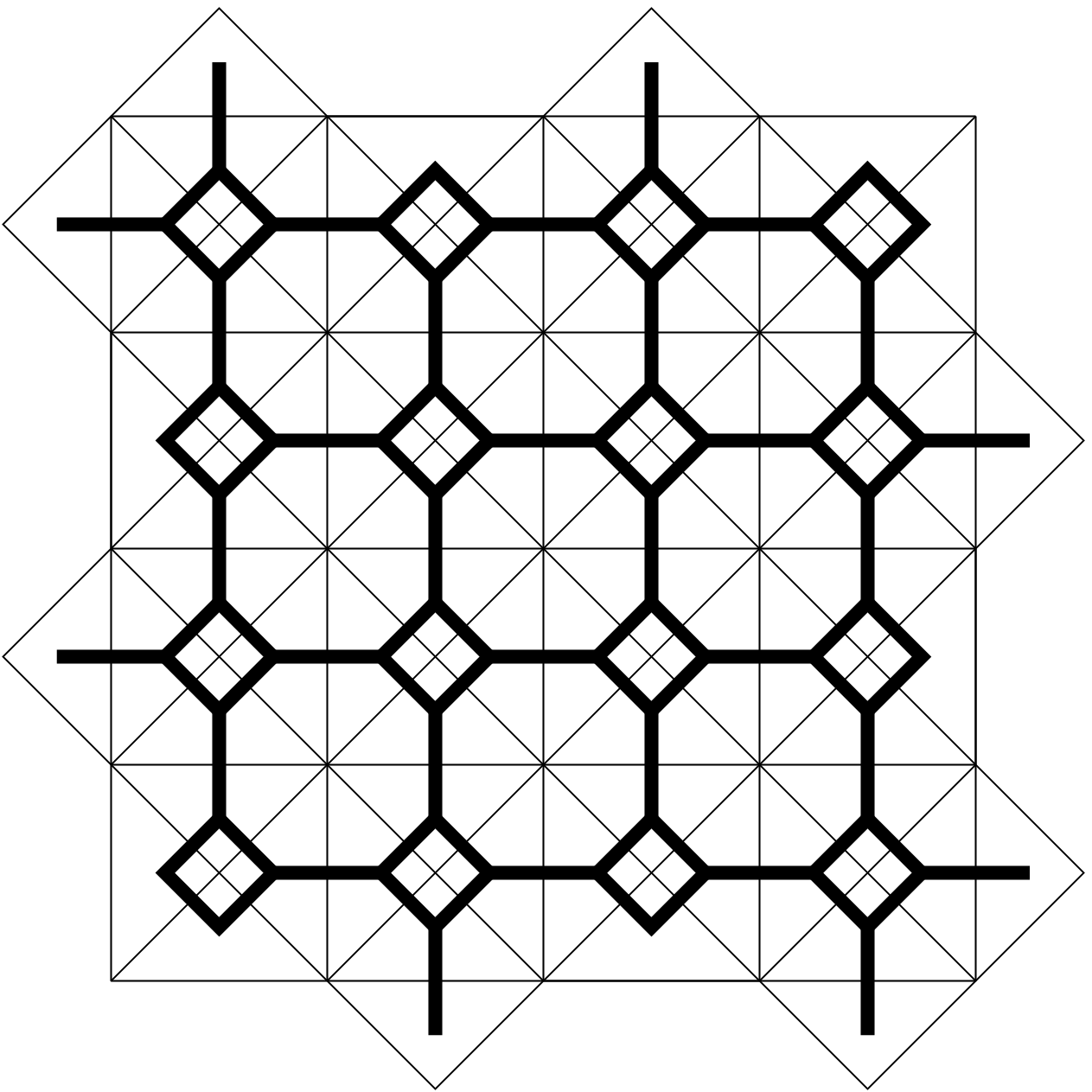}}
\twoline{Figure~4.1. {\rm The fortress of order 4.}}
{Figure~4.2. {\rm The dual graph of a fortress.}}
\endinsert

\topinsert
\centerline{\mypic{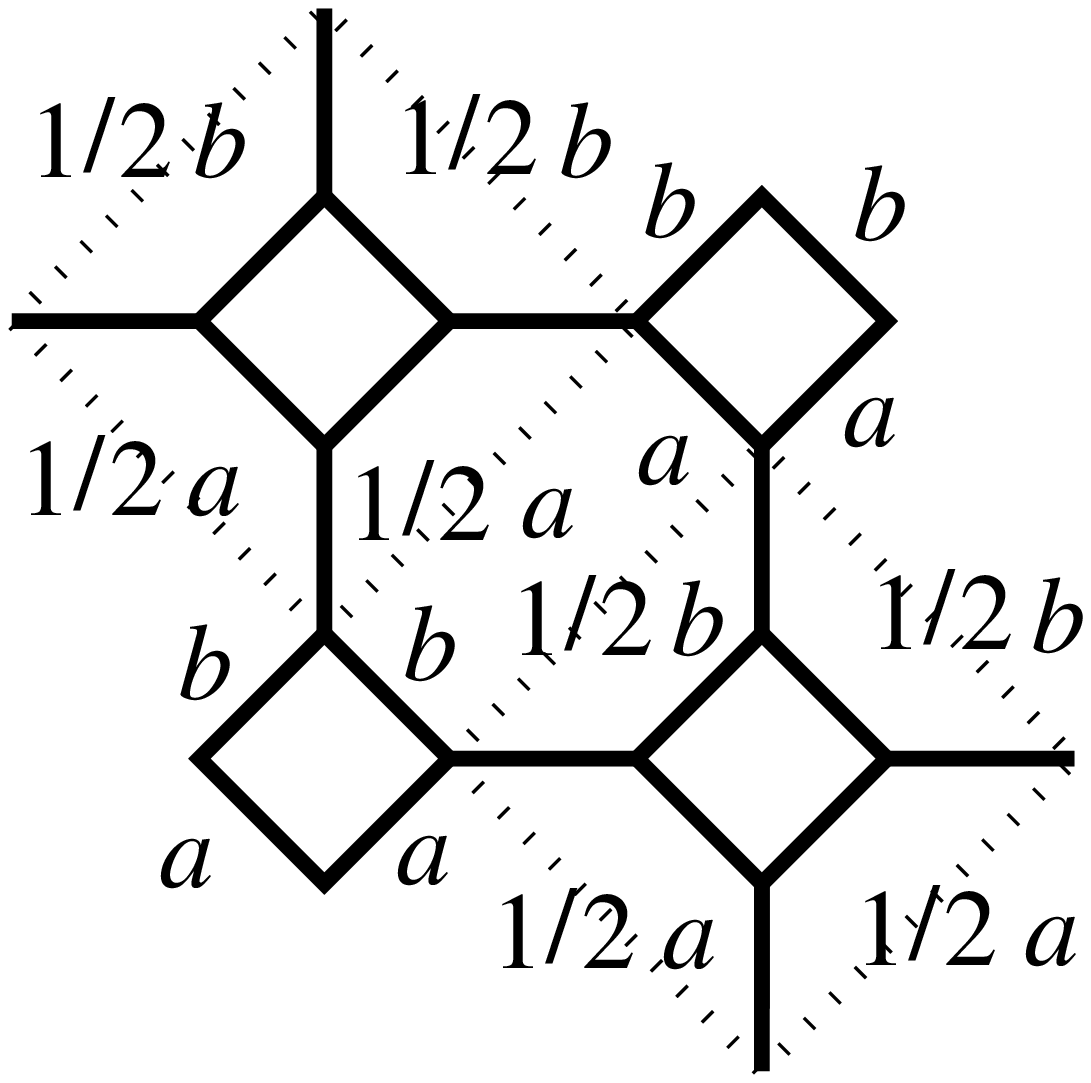}}
\centerline{{\smc Figure~4.3.} {\rm The natural weight of a fortress.}}
\endinsert

As it was the case for Aztec dungeons, if we assign different weights to the five
oriented tiles the resulting tiling generating function does not seem to be expressible
as a simple product of low-degree polynomials in the weights. However, if the two
triangular tiles whose hypothenuses face northeast and northwest have the same weight $a$,
and also the two remaining types have the same weight $b$, the weighted count is a simple
product. 

As noted before the statement of Theorem 3.1, there is no loss of generality in 
assuming that the weight of the square tiles is 1.

\proclaim{Corollary 4.5} The tiling generating function of a fortress with tiles weighted
as described above is a simple product of quadratic polynomials in $a$ and $b$.
\endproclaim

\pf The described tiling generating function is clearly the same as the matching
generating function of the dual graph of the fortress, with weights on edges induced by
the corresponding tiles (Figures 4.1 and 4.2 show the fortess of order 4 and its dual
graph). By choosing convenient local graph replacements of the type described in Lemma
2.5, computing it amounts to
finding the matching generating function of an Aztec diamond weighted with the periodic 
weight indicated in Figure 4.3. By scaling this weight matrix by $\sqrt{2}$ (with the effect of
multiplying the matching generating function by a power of $\sqrt{2}$), it becomes of
the type addressed in Theorem 4.1. Our statement follows then from Corollary 4.3.
\sendpf

\mysec{5. Aztec diamonds with yet another weight}

We have seen in the previous section a periodic weighting of the Aztec diamond that
factors into small factors. Another such factorization has been found by Stanley
(\cite{10}; see also \cite{2}) for the weight pattern
$$
\left[ \matrix
x_1&y_1&x_1&y_1&x_1&y_1&x_1&y_1&.&.&.\\
x_2&y_2&x_2&y_2&x_2&y_2&x_2&y_2&.&.&.\\
.&.&.&.&.&.&.&.\\
.&.&.&.&.&.&.&.\\
.&.&.&.&.&.&.&.\\  
x_{2n}&y_{2n}&x_{2n}&y_{2n}&x_{2n}&y_{2n}&x_{2n}&y_{2n}&.&.&.
\endmatrix\right].
$$
We show in this section that a third weight pattern, namely the one given by the matrix
$$
N=\left[ \matrix
x_1&y_1&y_1&x_1&x_1&y_1&y_1&x_1&.&.&.\\
x_2&y_2&y_2&x_2&x_2&y_2&y_2&x_2&.&.&.\\
.&.&.&.&.&.&.&.\\
.&.&.&.&.&.&.&.\\
.&.&.&.&.&.&.&.\\  
x_{2n}&y_{2n}&y_{2n}&x_{2n}&x_{2n}&y_{2n}&y_{2n}&x_{2n}&.&.&.
\endmatrix\right],
$$
also leads to a simple product formula. A matrix having the form above is determined
by its submatrix  
$$
T=\left[ \matrix
x_1&y_1\\
x_2&y_2\\
.&.\\
.&.\\
.&.\\  
x_{2n}&y_{2n}
\endmatrix\right].
$$
consisting of its first two columns. For notational simplicity,
denote the weight $\wt_N$ by $\z_T$.

Let $\te$ be the operator sending $T$ to the $(2n-4)\times2$ matrix
$$
\te(T)=\left[ \matrix
1/x_4&1/y_4\\
1/x_3&1/y_3\\
1/x_6&1/y_6\\
1/x_5&1/y_5\\
.&.\\
.&.\\
.&.\\  
1/x_{2n-2}&1/y_{2n-2}\\
1/x_{2n-3}&1/y_{2n-3}
\endmatrix\right].
$$
Let $\Delta_i=x_{2i-1}y_{2i}+x_{2i}y_{2i-1}$, for $i=1,\dotsc,n$

\proclaim{Theorem 5.1} For $n\geq1$ we have
$$
\M(AD_{n};\z_T)=2^{n-1}\prod_{i=1}^n\Delta_i
\prod_{i=1 \atop i\neq2,2n-1}^{2n}x_i^{n/2}y_i^{n/2-1}\M(AD_{n-2};\z_{\te(T)}),\tag5.1
$$
for $n$ even, and
$$
\M(AD_{n};\y_T)=2^{n-1}\prod_{i=1}^n\Delta_i
\prod_{i=1 \atop i\neq2,2n-1}^{2n}x_i^{(n-1)/2}y_i^{(n-1)/2}
\M(AD_{n-2};\z_{\te(T)}),\tag5.2
$$
for odd $n$.
\endproclaim

\pf Suppose first that $n$ is even. 
As for Theorem 4.1, we illustrate the proof for $n=4$. Our arguments generalize
with no difficulty to general even $n$. The weights on the edges of $AD_4$ form the
pattern
$$
B=\left[ \matrix
x_1&y_1&y_1&x_1&x_1&y_1&y_1&x_1\\
x_2&y_2&y_2&x_2&x_2&y_2&y_2&x_2\\
x_3&y_3&y_3&x_3&x_3&y_3&y_3&x_3\\
x_4&y_4&y_4&x_4&x_4&y_4&y_4&x_4\\
x_5&y_5&y_5&x_1&x_5&y_5&y_5&x_5\\
x_6&y_6&y_6&x_1&x_6&y_6&y_6&x_6\\
x_7&y_7&y_7&x_1&x_7&y_7&y_7&x_7\\
x_8&y_8&y_8&x_1&x_8&y_8&y_8&x_8\\
\endmatrix\right].
$$
Apply the Reduction Theorem. The derived weight is $\wt_C$, with
$$
C=\left[ \matrix
x_1/\Delta_1&x_1/\Delta_1&y_1/\Delta_1&y_1/\Delta_1&x_1/\Delta_1&x_1/\Delta_1\\
x_4/\Delta_2&x_4/\Delta_2&y_4/\Delta_2&y_4/\Delta_2&x_4/\Delta_2&x_4/\Delta_2\\
x_3/\Delta_2&x_3/\Delta_2&y_3/\Delta_2&y_3/\Delta_2&x_3/\Delta_2&x_3/\Delta_2\\
x_6/\Delta_3&x_6/\Delta_3&y_6/\Delta_3&y_6/\Delta_3&x_6/\Delta_3&x_6/\Delta_3\\
x_5/\Delta_3&x_5/\Delta_3&y_5/\Delta_3&y_5/\Delta_3&x_5/\Delta_3&x_5/\Delta_3\\
x_8/\Delta_4&x_8/\Delta_4&y_8/\Delta_4&y_8/\Delta_4&x_8/\Delta_4&x_8/\Delta_4\\
\endmatrix\right].
$$
The Reduction Theorem gives
$$
\M(AD_4;\wt_B)=\prod_{i=1}^4\Delta_i^4
\M(AD_{3};\wt_C).
\tag5.3
$$
By Lemma 4.4, we can factor out the denominators along rows to obtain
$$
\M(AD_{3};\wt_C)=\prod_{i=1}^4\Delta_i^{-3}\M(AD_{3};\wt_D),\tag5.4
$$
where
$$
D=\left[ \matrix
x_1&x_1&y_1&y_1&x_1&x_1\\
x_4&x_4&y_4&y_4&x_4&x_4\\
x_3&x_3&y_3&y_3&x_3&x_3\\
x_6&x_6&y_6&y_6&x_6&x_6\\
x_5&x_5&y_5&y_5&x_5&x_5\\
x_8&x_8&y_8&y_8&x_8&x_8
\endmatrix\right].
$$
Apply the Reduction Theorem again. All new cell-factors are monomials in the variables. We obtain
$$
\M(AD_{3};\wt_D)=2^{3^2}\prod_{i=1\atop i\neq2,7}^4x_i^{2}y_i
\M(AD_{3};\wt_E),\tag5.5
$$
with
$$
E=\left[ \matrix
1/(2x_4)&1/(2y_4)&1/(2y_4)&1/(2x_4)\\
1/(2x_3)&1/(2y_3)&1/(2y_3)&1/(2x_3)\\
1/(2x_6)&1/(2y_6)&1/(2y_6)&1/(2x_6)\\
1/(2x_5)&1/(2y_5)&1/(2y_5)&1/(2x_5)
\endmatrix\right].
$$
By Lemma 4.4, the factors of $1/2$ can be factored out, with an overall multiplicative contribution
of $1/2^6$. By (5.3)--(5.5) we obtain 
$$
\M(AD_4;\z_T)=2^3\prod_{i=1}^4\Delta_i
\prod_{i=1 \atop i\neq2,7}^{8}x_i^{2}y_i\M(AD_{2};\z_{\te(T)}),
$$
where $T$ is the submatrix of $B$ consisting of its first two columns. For general even $n$, the
same arguments prove (5.1).

The case when $n$ is odd can be treated similarly. The only difference occurs in (5.5), where the
exponents of $x_i$ and $y_i$ are now both $(n-1)/2$, as opposed to $n/2$ and $n/2-1$,
respectively, in the case of even $n$. Therefore, (5.3) and (5.4) imply (5.2). \endpf

\proclaim{Lemma 5.2} For any $2n\times2$ matrix $T$ having rows $[x_i,y_i]$, $i=1,\dotsc,2n$, 
we have
$$
\te^{(2)}(T)=\left[ \matrix
x_5&y_5\\
x_6&y_6\\
x_7&y_7\\
.&.\\
.&.\\
.&.\\  
x_{2n-4}&y_{2n-4}
\endmatrix\right].
$$
\endproclaim

\pf This is readily checked by the definition of the operator $\te$. \endpf

\proclaim{Corollary 5.3 (Yang \cite{11})} The matching generating function $\M(AD_{n};\z_T)$
can be expressed as a simple product in the variables $x_i$ and $y_i$, $i=1,\dotsc,2n$.
\endproclaim 

\pf By Lemma 5.2 and application of Theorem 5.1 two successive times we obtain that, for
$n\geq4$, 
$\M(AD_{n};\z_T)$ is a simple product times $\M(AD_{n-4};\z_{\bar{T}})$, where
$\bar{T}$ is the matrix obtained from $T$ by discarding the first four and last four rows.
Repeated application of this proves the statement of the Corollary. \sendpf

\mysec{6. Squares and hexagons}

Consider the lattice of squares and hexagons shown in Figure 6.1, and regard it as an
infinite graph $H$. This graph can
naturally be embedded in the square grid, as indicated in Figure 6.2. Draw an Aztec
diamond on the grid (see Figure 6.2; dotted lines indicate portions of the boundary
lying along missing edges).

\topinsert
\centerline{\mypic{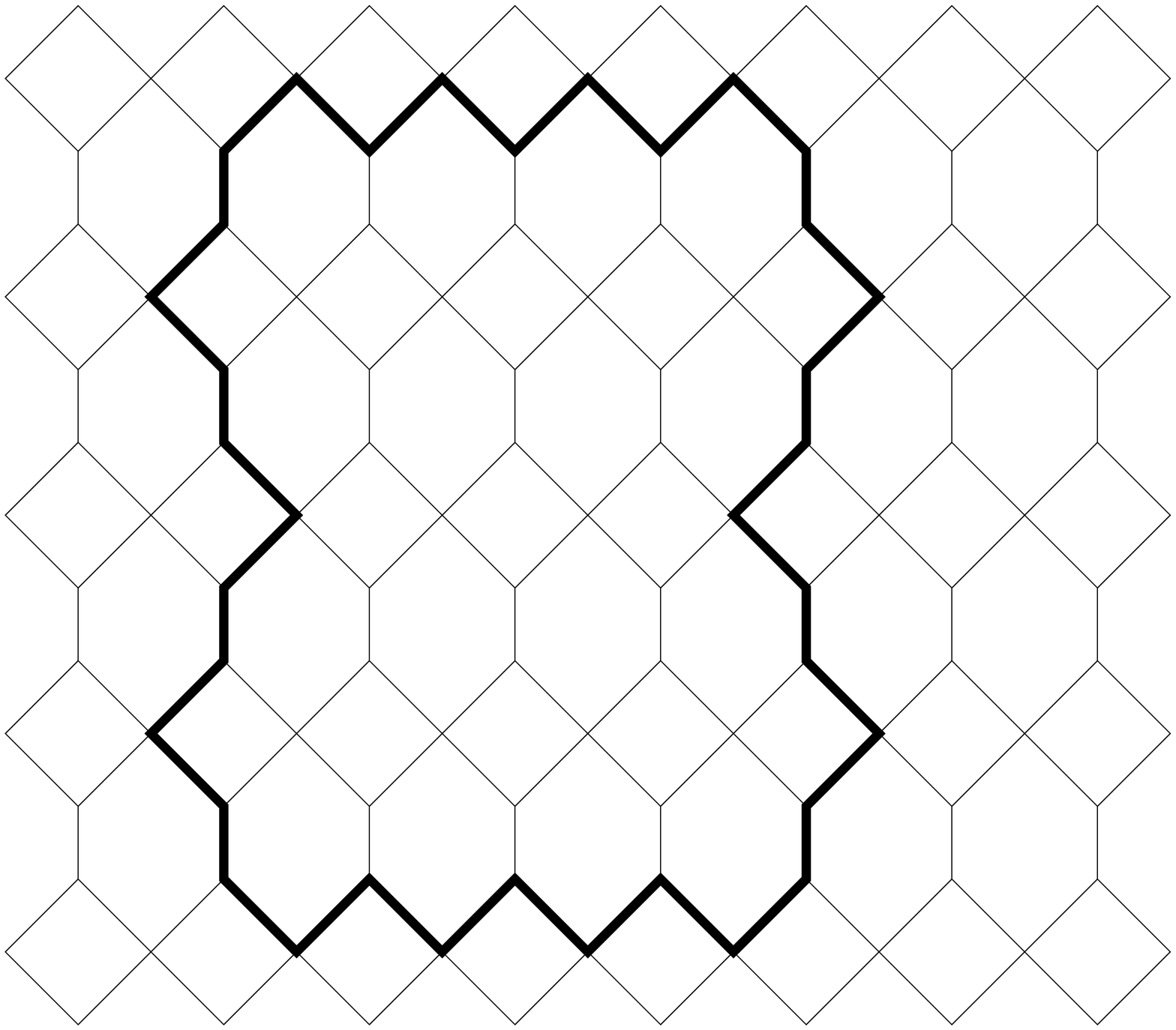}}
\centerline{{\smc Figure~6.1.} {\rm }}
\endinsert

\topinsert
\centerline{\mypic{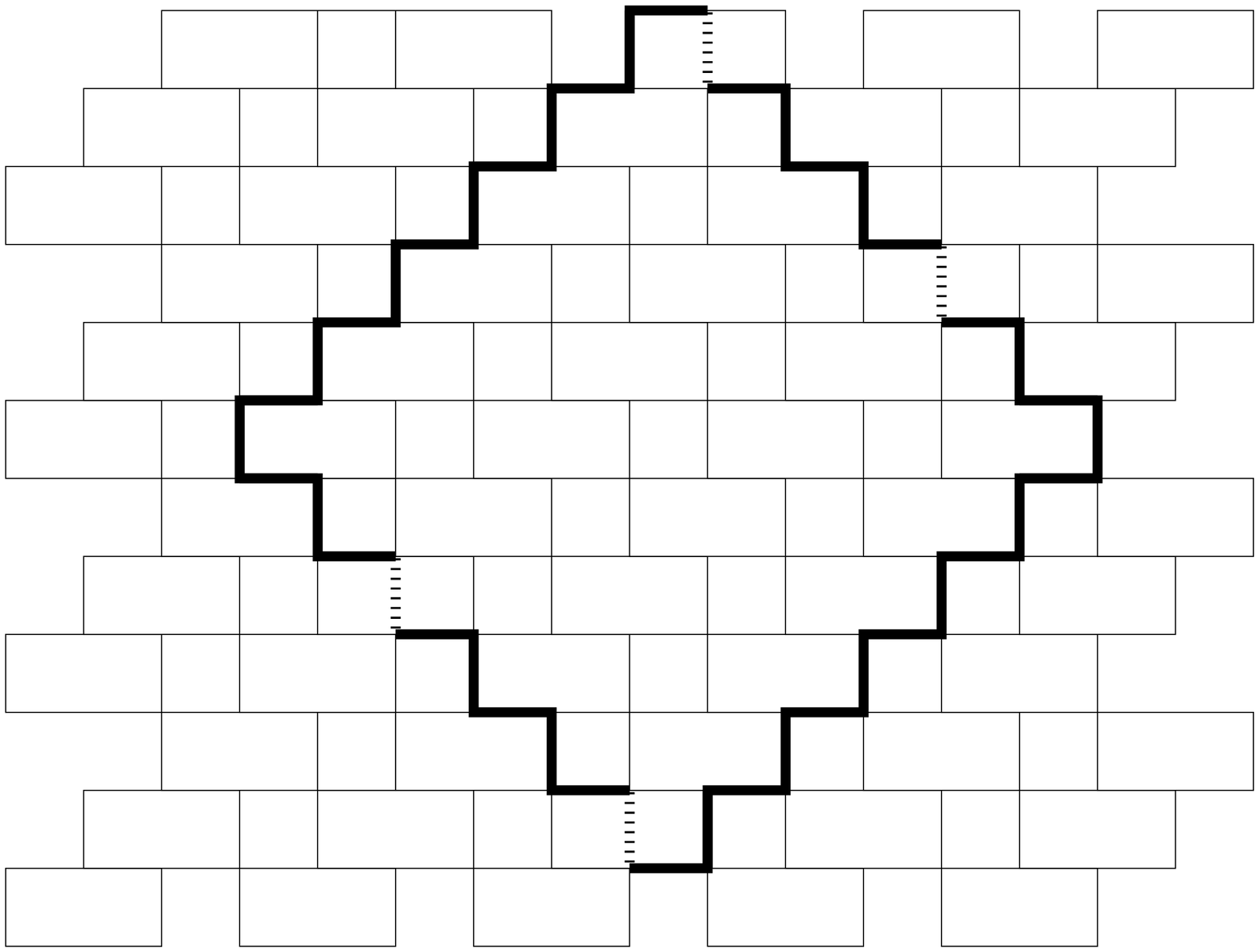}}
\centerline{{\smc Figure~6.2.} {\rm }}
\endinsert

A natural question to look at is to find the number
of perfect matchings of the portion of the embedding of $H$ that is a subgraph of our 
Aztec diamond (these subgraphs of $H$ have a natural definition also on the original 
squares and hexagons lattice; the example shown in Figure 6.2 is isomorphic, after
removing forced edges, to the subgraph of $H$ indicated in Figure 6.1). Propp
conjectured (in an e-mail post to the ``domino'' mailing list in February 1996)
that this number is always a power of 2. This was first
proved by Chris Douglas \cite{4}. 

We present in this section a weighted generalization 
of this result that keeps track of the orientation of the edges in the perfect matching.
As in the previous sections, the proof will follow by repeated applications of the
Reduction Theorem.

For definiteness, let us draw the Aztec diamond contour so that its westernmost edge 
belongs to a 4-cycle of $H$ lying outside the contour (this is the way the contour was 
drawn in Figure 6.2). Assign weight 1 to the horizontal edges, and weight $a$ to the 
vertical edges. It is easy to see that the weight on $AD_n$ arising this way is
precisely $\wt_A$, with period
$$
A=\left[ \matrix
1&0&1&a&1&a\\
a&1&a&1&0&1
\endmatrix\right].
$$
\proclaim{Theorem 6.1} For all $n\geq1$, $\M(AD_{2n};\wt_A)=\M(AD_{2n+1};\wt_A)=(1+a^2)^{n(n+1)}$.
\endproclaim

In our proof we will employ the following counterpart of Lemma 4.4.

\proclaim{Lemma 6.2} Partition the set of edges of $AD_n$ into $n$ classes as follows. View the set of
edges as forming a $2n\times 2n$ array. For each $i=1,\dotsc,n$ include the edges in columns $2i-1$ and 
$2i$ in a class. Then every perfect matching of $AD_n$ contains precisely $n+1$ edges in each class.
\endproclaim

\pf Let $\mu$ be a perfect matching of $AD_n$. For each class, consider the set of $n+1$ vertices of $AD_n$
incident only to edges in that class. These guarantee there will be at least $n+1$ edges of $\mu$ in each
class. Since the number of edges in $\mu$ is $n(n+1)$, the statement follows. \endpf

{\it Proof of Theorem 6.1.} Apply the Reduction Theorem. 
Since all the cell-factors are either 1 or $1+a^2$, we
obtain
$$
\M(AD_n;\wt_A)=(1+a^2)^{k_1}\M(AD_{n-1};\wt_{\de(A)}),\tag6.1
$$
where $k_1$ is a positive integer and
$$
\de(A)=\left[ \matrix
1&a/(1+a^2)&1/(1+a^2)&a&1&0\\
a&1/(1+a^2)&a/(1+a^2)&1&0&1
\endmatrix\right].
$$
Apply the Reduction Theorem again. All the cell-factors are now equal to 1, so we
obtain
$$
\M(AD_{n-1};\wt_{\de(A)})=\M(AD_{n-2};\wt_{\de^{(2)}(A)}),\tag6.2
$$
where 
$$
\de^{(2)}(A)=\left[ \matrix
1&a&1/(1+a^2)&0&1&a/(1+a^2)\\
a&1&a/(1+a^2)&1&0&1/(1+a^2)
\endmatrix\right].
$$
The cell-factors for $\wt_{\de^{(2)}(A)}$ are either $1+a^2$ or $(1+a^2)^{-1}$. One more
application of the Reduction Theorem yields
$$
\M(AD_{n-2};\wt_{\de^{(2)}(A)})=(1+a^2)^{k_2}\M(AD_{n-3};\wt_{\de^{(3)}(A)}),\tag6.3
$$
where $k_2$ is an integer and 
$$
\de^{(3)}(A)=\left[ \matrix
1/(1+a^2)&0&1&a&1+a^2&a/(1+a^2)\\
a/(1+a^2)&1+a^2&a&1&0&1/(1+a^2)
\endmatrix\right].
$$
Lemma 4.4 clearly holds also when partitioning the edges of an Aztec diamond into 
vertical blocks. Apply this version of Lemma 4.4 to the weight pattern $\wt_{A_3}$ 
we have on $AD_{n-3}$. Factor out from every third vertical block,
starting with the leftmost one, a factor of $(1+a^2)^{-2}$; leave all other blocks unchanged. 
By Lemma 4.4 we obtain
$$
\M(AD_{n-3};\wt_{\de^{(3)}(A)})=(1+a^2)^{k_3}\M(AD_{n-3};\wt_{B}),\tag6.4
$$
where $k_3$ is an integer and
$$
B=\left[ \matrix
1+a^2&0&1&a&1+a^2&a(1+a^2)\\
a(1+a^2)&1+a^2&a&1&0&1+a^2
\endmatrix\right].
$$
View $B$ as consisting of three $2\times2$ blocks. By Lemma 6.2, scaling the entries
within a block has a simple multiplicative effect on the matching generating function.
Factoring out the constants $(1+a^2)^{-1}$, $1$ and $(1+a^2)^{-1}$,
respectively, the resulting scaled matrix is precisely our original matrix $A$.
We obtain by Lemma 6.2 that
$$
\M(AD_{n-3};\wt_{B})=(1+a^2)^{k_4}\M(AD_{n-3};\wt_{A}),\tag6.5
$$
where $k_4$ is an integer.

It is easy to see that, because of forced edges, $\M(AD_{2n+1};\wt_A)=\M(AD_{2n};\wt_A)$ 
for all $n\geq1$. Therefore, by (6.1)--(6.5) we obtain, after working out the explicit values of 
$k_1,\dotsc,k_4$, that
$$
\M(AD_{2n+1};\wt_A)=(1+a^2)^{2n}\M(AD_{2n-2};\wt_A)=(1+a^2)^{2n}\M(AD_{2n-1};\wt_A).
$$
Repeated application of this proves the Theorem. \sendpf

\mysec{7. The Aztec dragon}

Consider the lattice of hexagons, squares and equilateral triangles illustrated in
Figure~7.1. Propp \cite{8} considered a family of regions $R_n$ called {\it Aztec
dragons} on this lattice ($R_6$ is outlined in Figure 7.1) and conjectured that the 
number of their tilings is $2^{n(n+1)}$. This was first proved, in work not yet published,
by Ben Wieland (as announced in \cite{8}). In this section we use the Reduction Theorem 
to obtain a weighted version of this result.

\topinsert
\centerline{\mypic{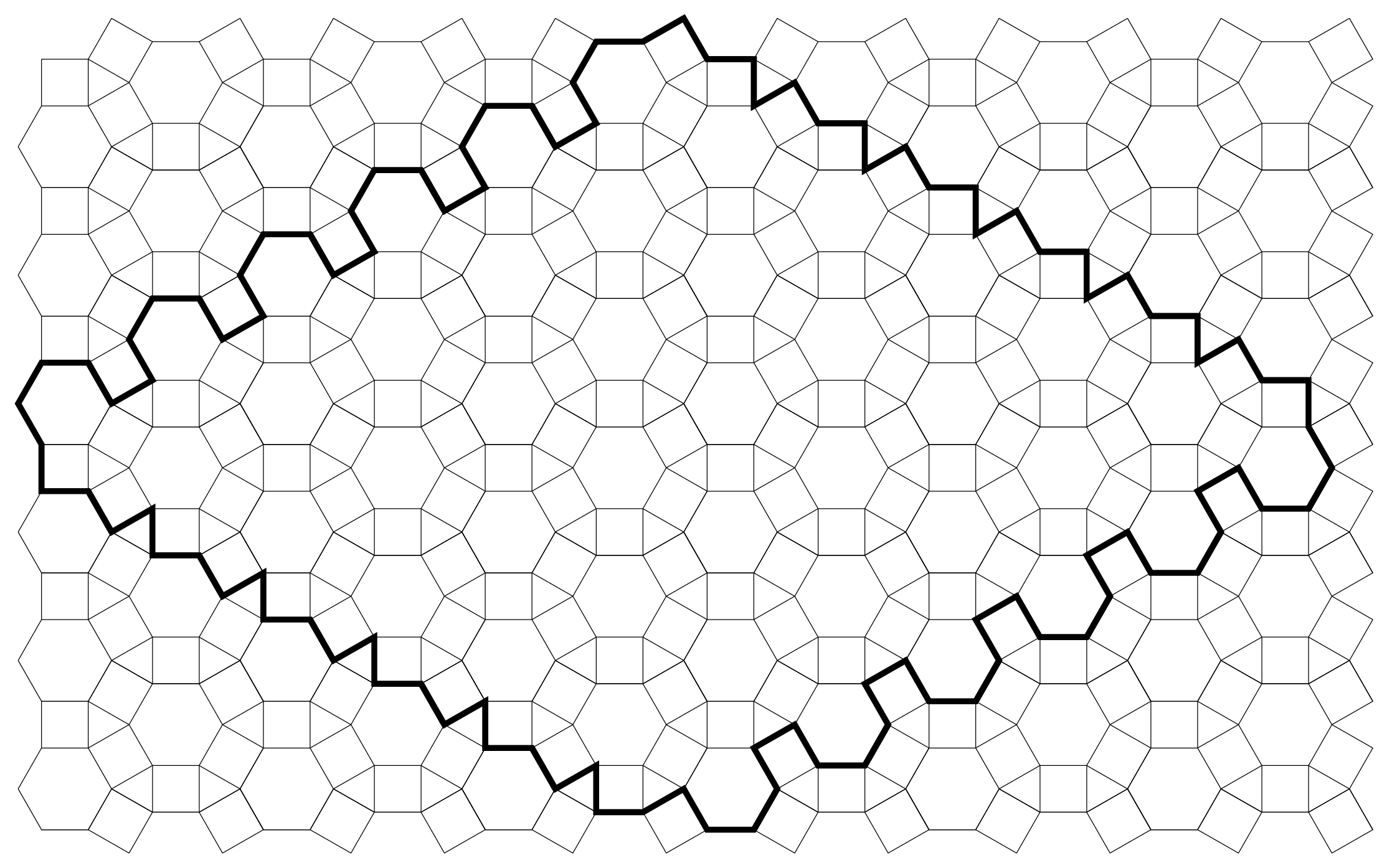}}
\centerline{{\smc Figure~7.1.} {\rm The Aztec dragon $R_6$.}}
\endinsert

There are two types of tiles. Since each hexagon must be paired with a square, and
each triangle also with a square, the number of tiles of each type is constant across
all tilings of a given dragon. Thus, our weighted version (Theorem 7.1 below) 
will not concern the
natural weight in the original set-up, but a certain other weight that is natural once 
we rephrase the problem.

Let $R_n'$ be the region obtained from $R_n$ by including $n$ square-triangle tiles
along its northwestern boundary, as indicated in Figure 7.2. Clearly, these tiles are
forced to be part of any tiling of the enlarged region. Enumerating the tilings of
$R_n$ is equivalent therefore to enumerating tilings of $R_n'$. 

As in the previous sections, consider the dual graph of $R_n'$ (for $n=6$ this is
pictured in Figure 7.3). It is easy to see that this graph can be deformed
isomorphically to the one illustrated in Figure 7.4.

\topinsert
\centerline{\mypic{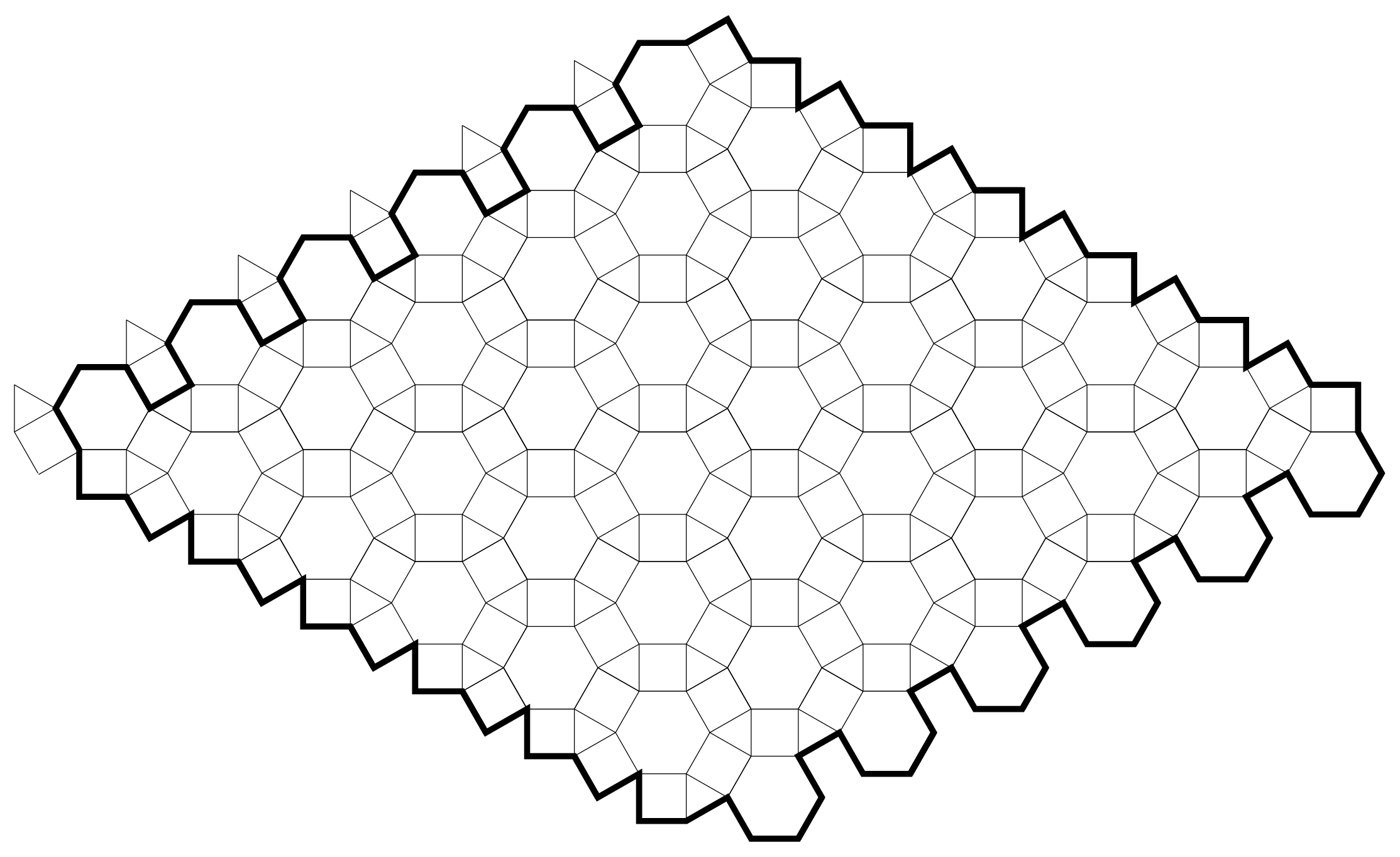}}
\centerline{{\smc Figure~7.2.} {\rm The region $R_6'$.}}

\medskip

\centerline{\mypic{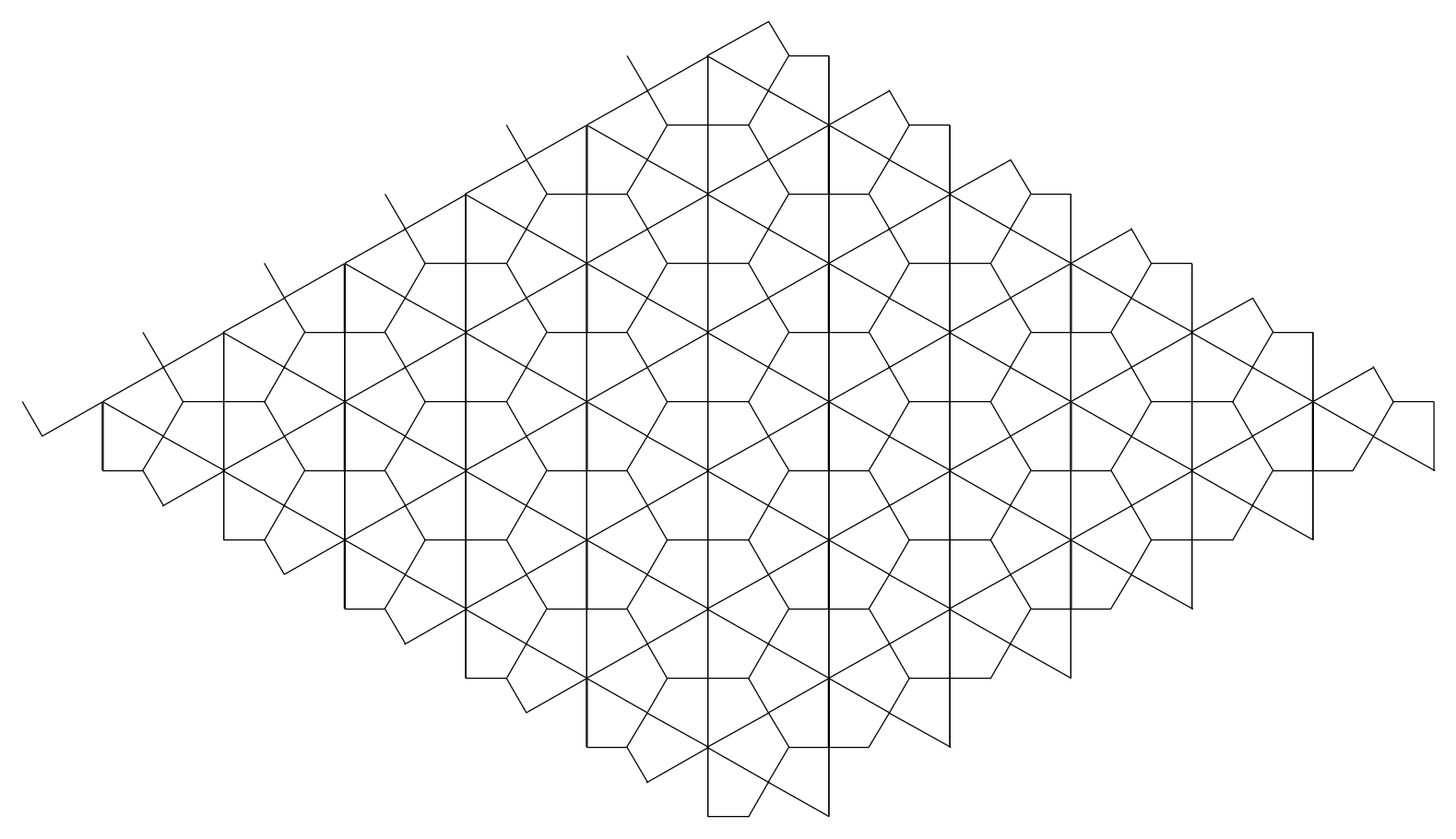}}
\centerline{{\smc Figure~7.3.} {\rm The dual graph of $R_6'$.}}
\endinsert

Furthermore, by employing the
``vertex-splitting'' trick (see e.g. Lemma 1.3 of \cite{1}), the latter graph is
readily seen to have its perfect matchings identified with those of the graph in
Figure 7.5 (this figure should be regarded as a subgraph of the grid graph; in
particular, it has vertices at the midpoints of the segments of length 2).

\topinsert
\twoline{\mypic{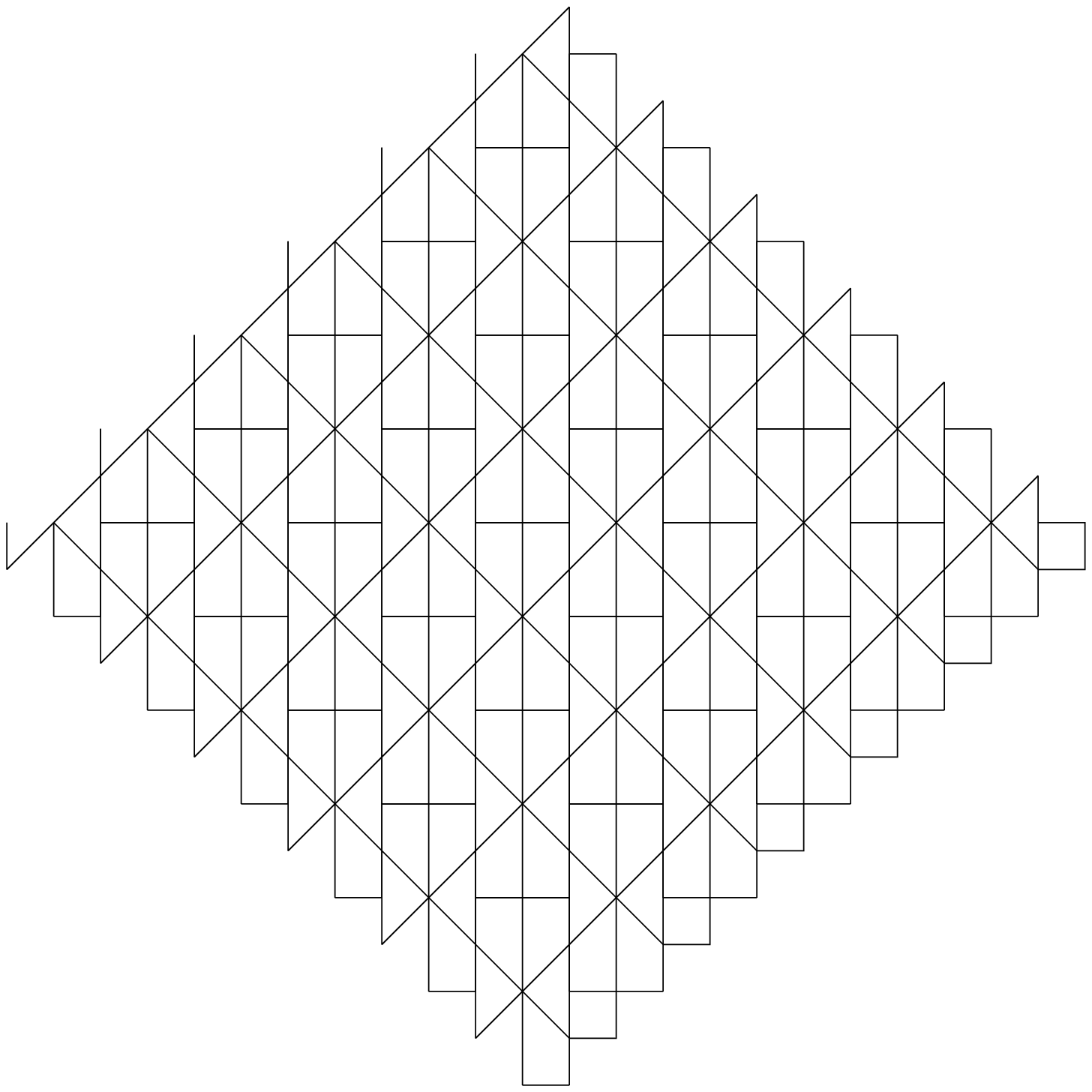}}{\mypic{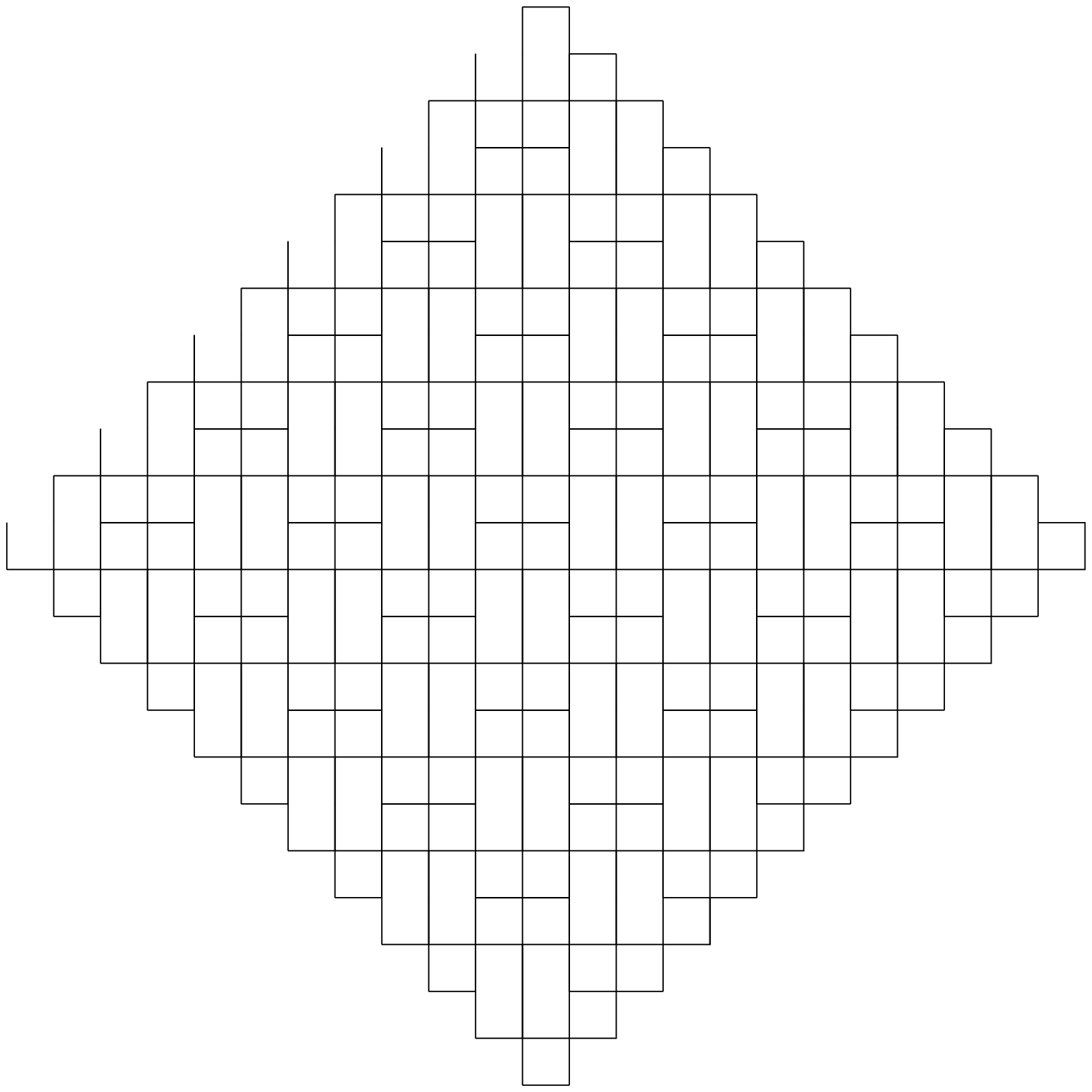}}
\twoline{Figure~7.4. {\rm An isomorphic deformation.}}
{Figure~7.5. {\rm A subgraph of $AD_{12}$.}}
\endinsert

However, the graph we obtained this way is just the periodic weighting $\wt_A$ of the 
Aztec diamond $AD_{2n}$, with
$$
A=\left[ \matrix
1&1&1&1\\
1&0&1&1\\
0&1&1&1\\
1&1&1&1\\
\endmatrix\right].
$$
Therefore, we have
$$
\M(R_n)=\M(AD_{2n};\wt_A).\tag7.1
$$
One can also regard the obtained graph (shown in Figure 7.5 for $n=6$) as a subgraph 
of $AD_{2n}$. Weight the vertical edges in this subgraph by 1, and the horizontals by
$a$, i.e. consider the weight $\wt_B$ on $AD_{2n}$, where
$$
B=\left[ \matrix
a&1&a&1\\
1&0&1&a\\
0&1&a&1\\
1&a&1&a\\
\endmatrix\right].
$$

\proclaim{Theorem 7.1} $\M(AD_{2n};\wt_B)=(1+a^2)^{n(n+1)}.$
\endproclaim

\proclaim{Corollary 7.2} $\T(R_n)=2^{n(n+1)}.$
\endproclaim

\pf This follows from (7.1) and the case $a=1$ of Theorem 7.1. \endpf

{\it Proof of Theorem 7.1.} By applying the Reduction Theorem twice and collecting
cell-factors one obtains
$$
\M(AD_{2n};\wt_B)=\Delta^{3n-1}\M(AD_{2n-2};\wt_{C}),\tag7.2
$$
where $\Delta=1+a^2$ and $C=\de^{(2)}(B)$ is the matrix
$$
C=\left[ \matrix
a&1&a/\Delta&1/\Delta\\
\Delta&0&1&a\\
0&\Delta&a&1\\
1&a&1/\Delta&a/\Delta\\
\endmatrix\right].
$$
By Lemma 4.4, we can factor out $1/\Delta$ from the first and last rows of $C$ to obtain
$$
\M(AD_{2n-2};\wt_C)=\left((1/\Delta)^{2n-2}\right)^n\M(AD_{2n-2};\wt_{D}),\tag7.3
$$
where
$$
D=\left[ \matrix
a\Delta&\Delta&a&1\\
\Delta&0&1&a\\
0&\Delta&a&1\\
\Delta&a\Delta&1&a\\
\endmatrix\right].
$$
In turn, by Lemma 6.2 we can factor out $\Delta$ from the first two columns of $D$ to get
$$
\M(AD_{2n-2};\wt_D)=\left(\Delta^{2n-1}\right)^{n-1}\M(AD_{2n-2};\wt_{B}).\tag7.4
$$
By (7.2)--(7.4) we obtain the recurrence
$$
\M(AD_{2n};\wt_B)=\Delta^{2n}\M(AD_{2n-2};\wt_{B}).
$$
Repeated application of this proves the statement of the Theorem. \endpf

\flushpar
{\smc Remark 7.3.} The number $2^{n(n+1)}$ of perfect matchings of the above subgraph of $AD_{2n}$
is approximately equal to the square root of the total number $2^{2n(2n+1)}$ of perfect matchings
of the full $AD_{2n}$.

\mysec{8. Powers of 3}

\topinsert
\centerline{\mypic{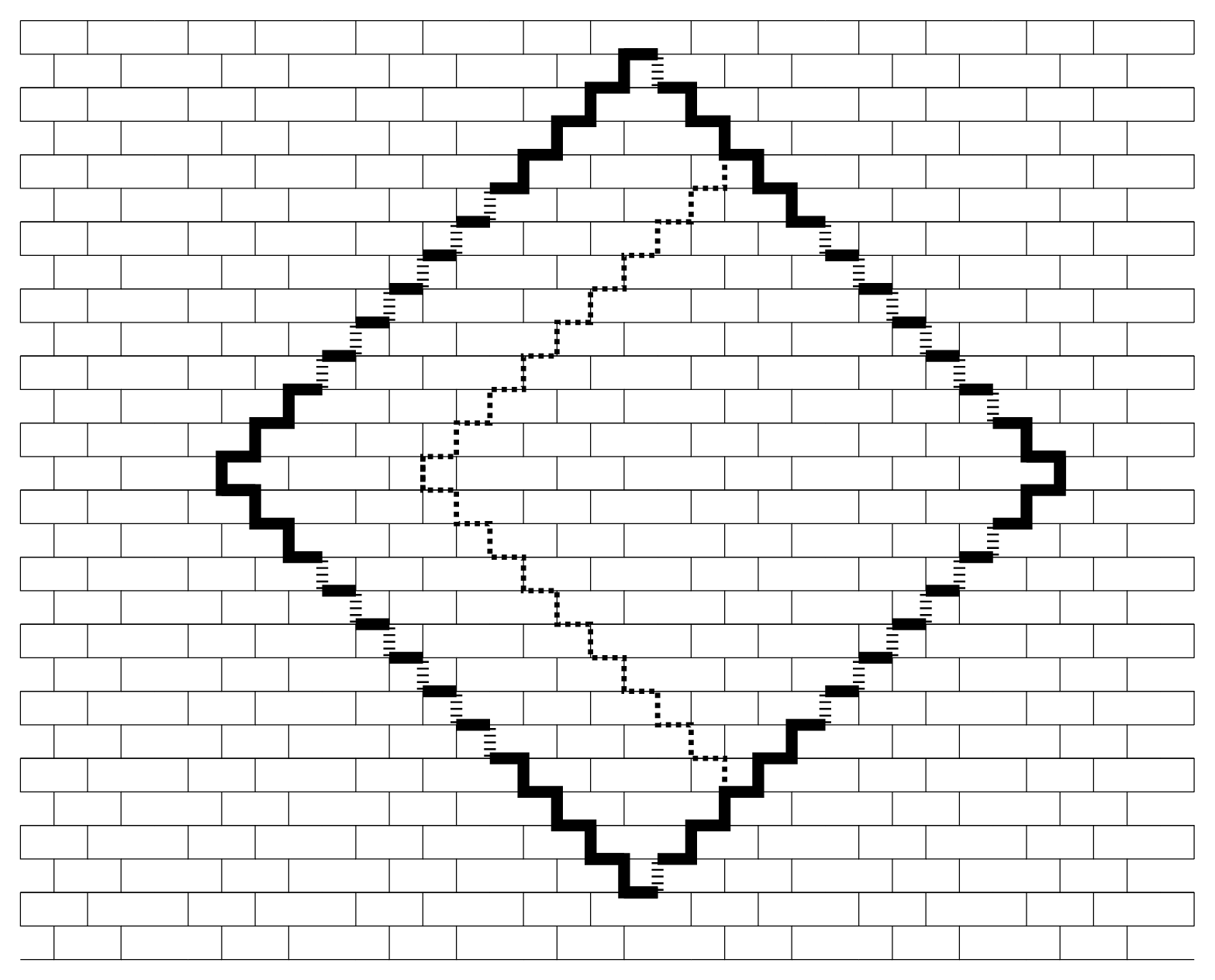}}
\centerline{{\smc Figure~8.1.} {\rm The region $B_{13}$.}}

\medskip

\twoline{\mypic{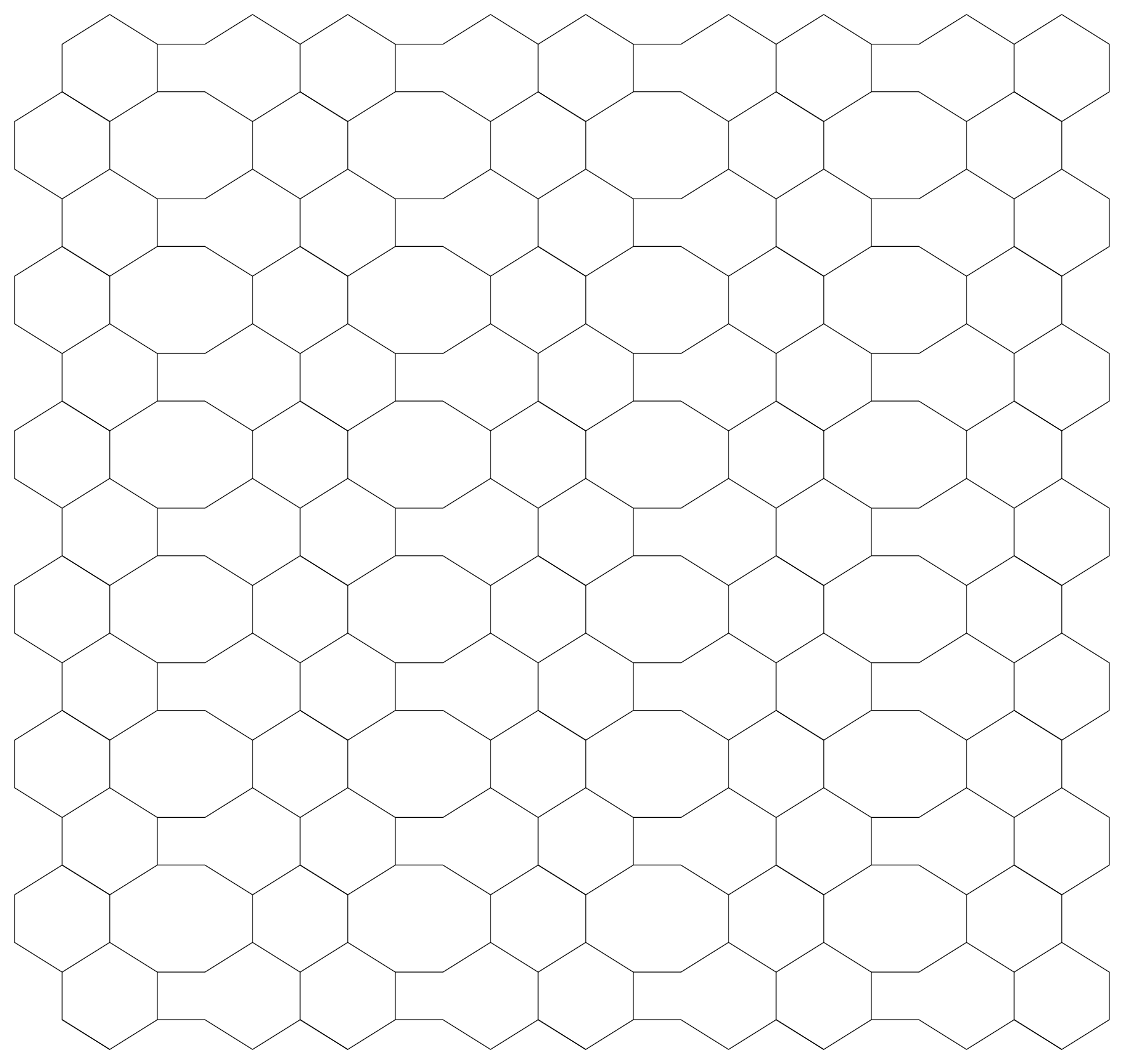}}{\mypic{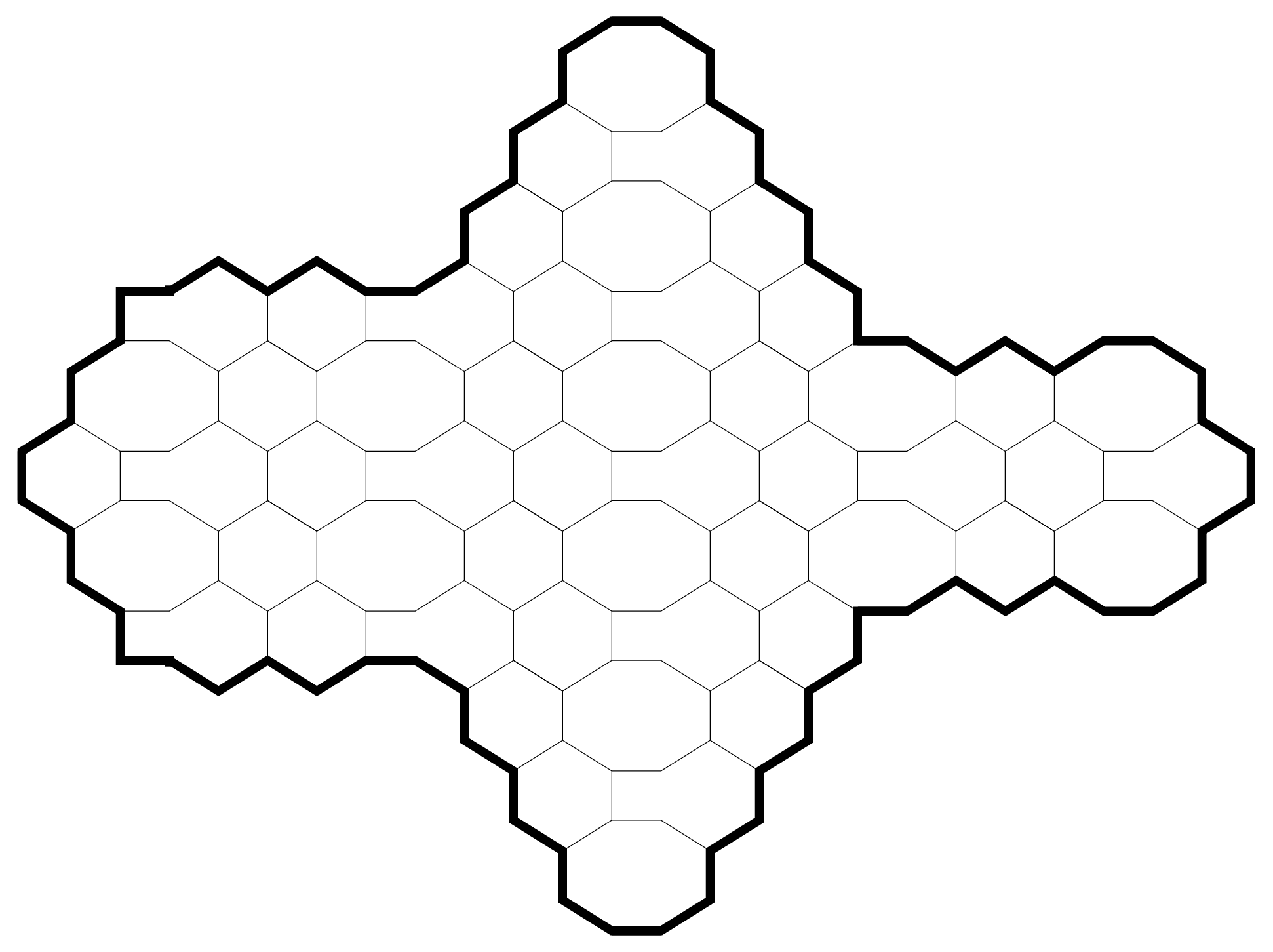}}
\twoline{Figure~8.2. {\rm A deformation of the lattice.}}
{Figure~8.3. {\rm $B_{13}$ with forced edges removed.}}
\endinsert

We have seen in the previous sections families of (unweighted) graphs whose number of
perfect matchings are perfect or near-perfect powers of the base 2, 5 or 13. What
about~3? Matt Blum 
has considered (as described in an e-mail message to James Propp dated December 1997)
a family of subgraphs of the square grid
for which he noticed (from data on concrete cases) that the number of perfect
matchings seems to be always a power of 3 or twice a power of 3. We prove in this
section that this is indeed the case.

Consider, following Blum, the sublattice of the square grid lattice showed in Figure 
8.1 and view it as an infinite graph $S$. Draw the
boundary of an Aztec diamond of order $n$ on this lattice in such a way that the
easternmost edge has an embedded hexagon east of it. Let $B_n$ be the induced subgraph
of $S$ spanned by the vertices lying inside or on the Aztec diamond boundary (the
graph inside the larger contour in Figure 8.1 is $B_{13}$). An alternative way to view
$B_n$ is as a subgraph of the lattice showed in Figure 8.2, which is just a
deformation of the lattice considered by Blum (the graph obtained from $B_{13}$ after
removing the forced edges is shown in this guise in Figure 8.3).

\proclaim{Theorem 8.1} For $n\geq31$ we have 
$$
\M(B_n)=3^{4x_n}\M(B_{n-30}),\tag8.1
$$
where $x_{5k+1}=4k-12$, $x_{5k+2}=4k-10$ and
$x_{5k+3}=x_{5k+4}=x_{5k+5}=4k-8$, $k\geq6$.
\endproclaim

The values of $\M(B_n)$ for $n=1,\dotsc,30$ are given in Table 8.1. By 
(8.1) one obtains that $\M(B_n)$ is either a perfect power of 3 or twice a perfect
power of 3, for all $n$.

\smallpagebreak
\pf The graph $B_n$ is the periodic weighting of $AD_n$ with period indicated by the
dashed contour in Figure 4.1. Rotating the indicated subgraph by $180^\circ$ (due to
our choice of the ``base'' edge being easternmost), we read off the period to be the
$20\times20$ matrix

\topinsert \eightpoint
$$\hbox{\vbox{\offinterlineskip
  \def\spacing{\omit\vrule height4pt&&&&&&\cr}
  \halign{\vrule\strut#&&\ \hfil$#$\hfil\ &\vrule#\cr
  \noalign{\hrule}\spacing
  &n&&\T(B_n)&&\text{factorization}&\cr
  \spacing\noalign{\hrule}\spacing
  &1&&1&&1&\cr\spacing
  &2&&2&&2&\cr\spacing
  &3&&6&&2\cdot3&\cr\spacing
  &4&&6&&2\cdot3&\cr\spacing
  &5&&6&&2\cdot3&\cr\spacing
  &6&&6&&2\cdot3&\cr\spacing
  &7&&27&&3^3&\cr\spacing
  &8&&486&&2\cdot3^5&\cr\spacing
  &9&&486&&2\cdot3^5&\cr\spacing
  &10&&486&&2\cdot3^5&\cr\spacing
  &11&&486&&2\cdot3^5&\cr\spacing
  &12&&6561&&3^8&\cr\spacing
  &13&&531441&&3^{12}&\cr\spacing
  &14&&531441&&3^{12}&\cr\spacing
  &15&&531441&&3^{12}&\cr\spacing  
  \noalign{\hrule}}}}
\qquad
\hbox{\vbox{\offinterlineskip
  \def\spacing{\omit\vrule height4pt&&&&&&\cr}
  \halign{\vrule\strut#&&\ \hfil$#$\hfil\ &\vrule#\cr
  \noalign{\hrule}\spacing
  &n&&\T(B_n)&&\text{factorization}&\cr
  \spacing\noalign{\hrule}\spacing
  &16&&531441&&3^{12}&\cr\spacing  
  &17&&43046721&&3^{16}&\cr\spacing  
  &18&&10460353203&&3^{21}&\cr\spacing  
  &19&&10460353203&&3^{21}&\cr\spacing  
  &20&&10460353203&&3^{21}&\cr\spacing  
  &21&&10460353203&&3^{21}&\cr\spacing  
  &22&&7625597484987&&3^{27}&\cr\spacing  
  &23&&5559060566555523&&3^{33}&\cr\spacing  
  &24&&5559060566555523&&3^{33}&\cr\spacing  
  &25&&5559060566555523&&3^{33}&\cr\spacing  
  &26&&5559060566555523&&3^{33}&\cr\spacing  
  &27&&24315330918113857602&&2\cdot3^{40}&\cr\spacing  
  &28&&79766443076872509863361&&3^{48}&\cr\spacing  
  &29&&79766443076872509863361&&3^{48}&\cr\spacing  
  &30&&79766443076872509863361&&3^{48}&\cr\spacing  
  \noalign{\hrule}}}}
$$
\tenpoint\centerline{{\smc Table~8.1.} The first thirty values of $\T(B_n)$.}
\endinsert

$$
A=\left[ \matrix 
\eightpoint1&\eightpoint1&\eightpoint1&\eightpoint1&\eightpoint0&\eightpoint1&\eightpoint0&\eightpoint1&\eightpoint0&\eightpoint1&\eightpoint0&\eightpoint1&\eightpoint0&\eightpoint1&\eightpoint0&\eightpoint1&\eightpoint1&\eightpoint1&\eightpoint1&\eightpoint1\\
\eightpoint1&\eightpoint0&\eightpoint1&\eightpoint0&\eightpoint1&\eightpoint1&\eightpoint1&\eightpoint1&\eightpoint1&\eightpoint1&\eightpoint1&\eightpoint1&\eightpoint1&\eightpoint0&\eightpoint1&\eightpoint0&\eightpoint1&\eightpoint0&\eightpoint1&\eightpoint0\\
\eightpoint1&\eightpoint1&\eightpoint0&\eightpoint1&\eightpoint0&\eightpoint1&\eightpoint0&\eightpoint1&\eightpoint0&\eightpoint1&\eightpoint0&\eightpoint1&\eightpoint0&\eightpoint1&\eightpoint1&\eightpoint1&\eightpoint1&\eightpoint1&\eightpoint1&\eightpoint1\\
\eightpoint1&\eightpoint0&\eightpoint1&\eightpoint1&\eightpoint1&\eightpoint1&\eightpoint1&\eightpoint1&\eightpoint1&\eightpoint1&\eightpoint1&\eightpoint0&\eightpoint1&\eightpoint0&\eightpoint1&\eightpoint0&\eightpoint1&\eightpoint0&\eightpoint1&\eightpoint0\\
\eightpoint0&\eightpoint1&\eightpoint0&\eightpoint1&\eightpoint0&\eightpoint1&\eightpoint0&\eightpoint1&\eightpoint0&\eightpoint1&\eightpoint0&\eightpoint1&\eightpoint1&\eightpoint1&\eightpoint1&\eightpoint1&\eightpoint1&\eightpoint1&\eightpoint1&\eightpoint1\\
\eightpoint1&\eightpoint1&\eightpoint1&\eightpoint1&\eightpoint1&\eightpoint1&\eightpoint1&\eightpoint1&\eightpoint1&\eightpoint0&\eightpoint1&\eightpoint0&\eightpoint1&\eightpoint0&\eightpoint1&\eightpoint0&\eightpoint1&\eightpoint0&\eightpoint1&\eightpoint0\\
\eightpoint0&\eightpoint1&\eightpoint0&\eightpoint1&\eightpoint0&\eightpoint1&\eightpoint0&\eightpoint1&\eightpoint0&\eightpoint1&\eightpoint1&\eightpoint1&\eightpoint1&\eightpoint1&\eightpoint1&\eightpoint1&\eightpoint1&\eightpoint1&\eightpoint0&\eightpoint1\\
\eightpoint1&\eightpoint1&\eightpoint1&\eightpoint1&\eightpoint1&\eightpoint1&\eightpoint1&\eightpoint0&\eightpoint1&\eightpoint0&\eightpoint1&\eightpoint0&\eightpoint1&\eightpoint0&\eightpoint1&\eightpoint0&\eightpoint1&\eightpoint0&\eightpoint1&\eightpoint1\\
\eightpoint0&\eightpoint1&\eightpoint0&\eightpoint1&\eightpoint0&\eightpoint1&\eightpoint0&\eightpoint1&\eightpoint1&\eightpoint1&\eightpoint1&\eightpoint1&\eightpoint1&\eightpoint1&\eightpoint1&\eightpoint1&\eightpoint0&\eightpoint1&\eightpoint0&\eightpoint1\\
\eightpoint1&\eightpoint1&\eightpoint1&\eightpoint1&\eightpoint1&\eightpoint0&\eightpoint1&\eightpoint0&\eightpoint1&\eightpoint0&\eightpoint1&\eightpoint0&\eightpoint1&\eightpoint0&\eightpoint1&\eightpoint0&\eightpoint1&\eightpoint1&\eightpoint1&\eightpoint1\\
\eightpoint0&\eightpoint1&\eightpoint0&\eightpoint1&\eightpoint0&\eightpoint1&\eightpoint1&\eightpoint1&\eightpoint1&\eightpoint1&\eightpoint1&\eightpoint1&\eightpoint1&\eightpoint1&\eightpoint0&\eightpoint1&\eightpoint0&\eightpoint1&\eightpoint0&\eightpoint1\\
\eightpoint1&\eightpoint1&\eightpoint1&\eightpoint0&\eightpoint1&\eightpoint0&\eightpoint1&\eightpoint0&\eightpoint1&\eightpoint0&\eightpoint1&\eightpoint0&\eightpoint1&\eightpoint0&\eightpoint1&\eightpoint1&\eightpoint1&\eightpoint1&\eightpoint1&\eightpoint1\\
\eightpoint0&\eightpoint1&\eightpoint0&\eightpoint1&\eightpoint1&\eightpoint1&\eightpoint1&\eightpoint1&\eightpoint1&\eightpoint1&\eightpoint1&\eightpoint1&\eightpoint0&\eightpoint1&\eightpoint0&\eightpoint1&\eightpoint0&\eightpoint1&\eightpoint0&\eightpoint1\\
\eightpoint1&\eightpoint0&\eightpoint1&\eightpoint0&\eightpoint1&\eightpoint0&\eightpoint1&\eightpoint0&\eightpoint1&\eightpoint0&\eightpoint1&\eightpoint0&\eightpoint1&\eightpoint1&\eightpoint1&\eightpoint1&\eightpoint1&\eightpoint1&\eightpoint1&\eightpoint1\\
\eightpoint0&\eightpoint1&\eightpoint1&\eightpoint1&\eightpoint1&\eightpoint1&\eightpoint1&\eightpoint1&\eightpoint1&\eightpoint1&\eightpoint0&\eightpoint1&\eightpoint0&\eightpoint1&\eightpoint0&\eightpoint1&\eightpoint0&\eightpoint1&\eightpoint0&\eightpoint1\\
\eightpoint1&\eightpoint0&\eightpoint1&\eightpoint0&\eightpoint1&\eightpoint0&\eightpoint1&\eightpoint0&\eightpoint1&\eightpoint0&\eightpoint1&\eightpoint1&\eightpoint1&\eightpoint1&\eightpoint1&\eightpoint1&\eightpoint1&\eightpoint1&\eightpoint1&\eightpoint0\\
\eightpoint1&\eightpoint1&\eightpoint1&\eightpoint1&\eightpoint1&\eightpoint1&\eightpoint1&\eightpoint1&\eightpoint0&\eightpoint1&\eightpoint0&\eightpoint1&\eightpoint0&\eightpoint1&\eightpoint0&\eightpoint1&\eightpoint0&\eightpoint1&\eightpoint0&\eightpoint1\\
\eightpoint1&\eightpoint0&\eightpoint1&\eightpoint0&\eightpoint1&\eightpoint0&\eightpoint1&\eightpoint0&\eightpoint1&\eightpoint1&\eightpoint1&\eightpoint1&\eightpoint1&\eightpoint1&\eightpoint1&\eightpoint1&\eightpoint1&\eightpoint0&\eightpoint1&\eightpoint0\\
\eightpoint1&\eightpoint1&\eightpoint1&\eightpoint1&\eightpoint1&\eightpoint1&\eightpoint0&\eightpoint1&\eightpoint0&\eightpoint1&\eightpoint0&\eightpoint1&\eightpoint0&\eightpoint1&\eightpoint0&\eightpoint1&\eightpoint0&\eightpoint1&\eightpoint1&\eightpoint1\\
\eightpoint1&\eightpoint0&\eightpoint1&\eightpoint0&\eightpoint1&\eightpoint0&\eightpoint1&\eightpoint1&\eightpoint1&\eightpoint1&\eightpoint1&\eightpoint1&\eightpoint1&\eightpoint1&\eightpoint1&\eightpoint0&\eightpoint1&\eightpoint0&\eightpoint1&\eightpoint0
\endmatrix\right].
$$

By Lemma 3.5, the periodic weightings arising by successive applications of
the Reduction Theorem have periods given by the orbit of $A$ under the operator $\de$
defined in Section~3. Using a computer algebra program (for instance Maple) one can
easily compute $\de^{(i)}(A)$ for successive values of $i$. One finds that all entries
of $\de^{(30)}(A)$ are either 0 or integer powers of 9, with the position of the
zeroes matching perfectly the positions of the zeros of $A$. The exponents of 9 form
the pattern
$$
\left[ 
\matrix
{\scriptscriptstyle{4}}\!\!\!&{\scriptscriptstyle{1}}\!\!\!&{\scriptscriptstyle{-2}}\!\!\!&{\scriptscriptstyle{-1}}\!\!\!&{\scriptscriptstyle{-\infty}}\!\!\!&{\scriptscriptstyle{-1}}\!\!\!&{\scriptscriptstyle{-\infty}}\!\!\!&{\scriptscriptstyle{1}}\!\!\!&{\scriptscriptstyle{-\infty}}\!\!\!&{\scriptscriptstyle{0}}\!\!\!&{\scriptscriptstyle{-\infty}}\!\!\!&{\scriptscriptstyle{-1}}\!\!\!&{\scriptscriptstyle{-\infty}}\!\!\!&{\scriptscriptstyle{1}}\!\!\!&{\scriptscriptstyle{-\infty}}\!\!\!&{\scriptscriptstyle{1}}\!\!\!&{\scriptscriptstyle{2}}\!\!\!&{\scriptscriptstyle{-1}}\!\!\!&{\scriptscriptstyle{-4}}\!\!\!&{\scriptscriptstyle{0}}\\
{\scriptscriptstyle{1}}\!\!\!&{\scriptscriptstyle{-\infty}}\!\!\!&{\scriptscriptstyle{-1}}\!\!\!&{\scriptscriptstyle{-\infty}}\!\!\!&{\scriptscriptstyle{-1}}\!\!\!&{\scriptscriptstyle{-2}}\!\!\!&{\scriptscriptstyle{1}}\!\!\!&{\scriptscriptstyle{4}}\!\!\!&{\scriptscriptstyle{0}}\!\!\!&{\scriptscriptstyle{-4}}\!\!\!&{\scriptscriptstyle{-1}}\!\!\!&{\scriptscriptstyle{2}}\!\!\!&{\scriptscriptstyle{1}}\!\!\!&{\scriptscriptstyle{-\infty}}\!\!\!&{\scriptscriptstyle{1}}\!\!\!&{\scriptscriptstyle{-\infty}}\!\!\!&{\scriptscriptstyle{-1}}\!\!\!&{\scriptscriptstyle{-\infty}}\!\!\!&{\scriptscriptstyle{0}}\!\!\!&{\scriptscriptstyle{-\infty}}\\
{\scriptscriptstyle{-2}}\!\!\!&{\scriptscriptstyle{-1}}\!\!\!&{\scriptscriptstyle{-\infty}}\!\!\!&{\scriptscriptstyle{-1}}\!\!\!&{\scriptscriptstyle{-\infty}}\!\!\!&{\scriptscriptstyle{1}}\!\!\!&{\scriptscriptstyle{-\infty}}\!\!\!&{\scriptscriptstyle{0}}\!\!\!&{\scriptscriptstyle{-\infty}}\!\!\!&{\scriptscriptstyle{-1}}\!\!\!&{\scriptscriptstyle{-\infty}}\!\!\!&{\scriptscriptstyle{1}}\!\!\!&{\scriptscriptstyle{-\infty}}\!\!\!&{\scriptscriptstyle{1}}\!\!\!&{\scriptscriptstyle{2}}\!\!\!&{\scriptscriptstyle{-1}}\!\!\!&{\scriptscriptstyle{-4}}\!\!\!&{\scriptscriptstyle{0}}\!\!\!&{\scriptscriptstyle{4}}\!\!\!&{\scriptscriptstyle{1}}\\ 
{\scriptscriptstyle{-1}}\!\!\!&{\scriptscriptstyle{-\infty}}\!\!\!&{\scriptscriptstyle{-1}}\!\!\!&{\scriptscriptstyle{-2}}\!\!\!&{\scriptscriptstyle{1}}\!\!\!&{\scriptscriptstyle{4}}\!\!\!&{\scriptscriptstyle{0}}\!\!\!&{\scriptscriptstyle{-4}}\!\!\!&{\scriptscriptstyle{-1}}\!\!\!&{\scriptscriptstyle{2}}\!\!\!&{\scriptscriptstyle{1}}\!\!\!&{\scriptscriptstyle{-\infty}}\!\!\!&{\scriptscriptstyle{1}}\!\!\!&{\scriptscriptstyle{-\infty}}\!\!\!&{\scriptscriptstyle{-1}}\!\!\!&{\scriptscriptstyle{-\infty}}\!\!\!&{\scriptscriptstyle{0}}\!\!\!&{\scriptscriptstyle{-\infty}}\!\!\!&{\scriptscriptstyle{1}}\!\!\!&{\scriptscriptstyle{-\infty}}\\ 
{\scriptscriptstyle{-\infty}}\!\!\!&{\scriptscriptstyle{-1}}\!\!\!&{\scriptscriptstyle{-\infty}}\!\!\!&{\scriptscriptstyle{1}}\!\!\!&{\scriptscriptstyle{-\infty}}\!\!\!&{\scriptscriptstyle{0}}\!\!\!&{\scriptscriptstyle{-\infty}}\!\!\!&{\scriptscriptstyle{-1}}\!\!\!&{\scriptscriptstyle{-\infty}}\!\!\!&{\scriptscriptstyle{1}}\!\!\!&{\scriptscriptstyle{-\infty}}\!\!\!&{\scriptscriptstyle{1}}\!\!\!&{\scriptscriptstyle{2}}\!\!\!&{\scriptscriptstyle{-1}}\!\!\!&{\scriptscriptstyle{-4}}\!\!\!&{\scriptscriptstyle{0}}\!\!\!&{\scriptscriptstyle{4}}\!\!\!&{\scriptscriptstyle{1}}\!\!\!&{\scriptscriptstyle{-2}}\!\!\!&{\scriptscriptstyle{-1}}\\
{\scriptscriptstyle{-1}}\!\!\!&{\scriptscriptstyle{-2}}\!\!\!&{\scriptscriptstyle{1}}\!\!\!&{\scriptscriptstyle{4}}\!\!\!&{\scriptscriptstyle{0}}\!\!\!&{\scriptscriptstyle{-4}}\!\!\!&{\scriptscriptstyle{-1}}\!\!\!&{\scriptscriptstyle{2}}\!\!\!&{\scriptscriptstyle{1}}\!\!\!&{\scriptscriptstyle{-\infty}}\!\!\!&{\scriptscriptstyle{1}}\!\!\!&{\scriptscriptstyle{-\infty}}\!\!\!&{\scriptscriptstyle{-1}}\!\!\!&{\scriptscriptstyle{-\infty}}\!\!\!&{\scriptscriptstyle{0}}\!\!\!&{\scriptscriptstyle{-\infty}}\!\!\!&{\scriptscriptstyle{1}}\!\!\!&{\scriptscriptstyle{-\infty}}\!\!\!&{\scriptscriptstyle{-1}}\!\!\!&{\scriptscriptstyle{-\infty}}\\ 
{\scriptscriptstyle{-\infty}}\!\!\!&{\scriptscriptstyle{1}}\!\!\!&{\scriptscriptstyle{-\infty}}\!\!\!&{\scriptscriptstyle{0}}\!\!\!&{\scriptscriptstyle{-\infty}}\!\!\!&{\scriptscriptstyle{-1}}\!\!\!&{\scriptscriptstyle{-\infty}}\!\!\!&{\scriptscriptstyle{1}}\!\!\!&{\scriptscriptstyle{-\infty}}\!\!\!&{\scriptscriptstyle{1}}\!\!\!&{\scriptscriptstyle{2}}\!\!\!&{\scriptscriptstyle{-1}}\!\!\!&{\scriptscriptstyle{-4}}\!\!\!&{\scriptscriptstyle{0}}\!\!\!&{\scriptscriptstyle{4}}\!\!\!&{\scriptscriptstyle{1}}\!\!\!&{\scriptscriptstyle{-2}}\!\!\!&{\scriptscriptstyle{-1}}\!\!\!&{\scriptscriptstyle{-\infty}}\!\!\!&{\scriptscriptstyle{-1}}\\
{\scriptscriptstyle{1}}\!\!\!&{\scriptscriptstyle{4}}\!\!\!&{\scriptscriptstyle{0}}\!\!\!&{\scriptscriptstyle{-4}}\!\!\!&{\scriptscriptstyle{-1}}\!\!\!&{\scriptscriptstyle{2}}\!\!\!&{\scriptscriptstyle{1}}\!\!\!&{\scriptscriptstyle{-\infty}}\!\!\!&{\scriptscriptstyle{1}}\!\!\!&{\scriptscriptstyle{-\infty}}\!\!\!&{\scriptscriptstyle{-1}}\!\!\!&{\scriptscriptstyle{-\infty}}\!\!\!&{\scriptscriptstyle{0}}\!\!\!&{\scriptscriptstyle{-\infty}}\!\!\!&{\scriptscriptstyle{1}}\!\!\!&{\scriptscriptstyle{-\infty}}\!\!\!&{\scriptscriptstyle{-1}}\!\!\!&{\scriptscriptstyle{-\infty}}\!\!\!&{\scriptscriptstyle{-1}}\!\!\!&{\scriptscriptstyle{-2}}\\
{\scriptscriptstyle{-\infty}}\!\!\!&{\scriptscriptstyle{0}}\!\!\!&{\scriptscriptstyle{-\infty}}\!\!\!&{\scriptscriptstyle{-1}}\!\!\!&{\scriptscriptstyle{-\infty}}\!\!\!&{\scriptscriptstyle{1}}\!\!\!&{\scriptscriptstyle{-\infty}}\!\!\!&{\scriptscriptstyle{1}}\!\!\!&{\scriptscriptstyle{2}}\!\!\!&{\scriptscriptstyle{-1}}\!\!\!&{\scriptscriptstyle{-4}}\!\!\!&{\scriptscriptstyle{0}}\!\!\!&{\scriptscriptstyle{4}}\!\!\!&{\scriptscriptstyle{1}}\!\!\!&{\scriptscriptstyle{-2}}\!\!\!&{\scriptscriptstyle{-1}}\!\!\!&{\scriptscriptstyle{-\infty}}\!\!\!&{\scriptscriptstyle{-1}}\!\!\!&{\scriptscriptstyle{-\infty}}\!\!\!&{\scriptscriptstyle{1}}\\
{\scriptscriptstyle{0}}\!\!\!&{\scriptscriptstyle{-4}}\!\!\!&{\scriptscriptstyle{-1}}\!\!\!&{\scriptscriptstyle{2}}\!\!\!&{\scriptscriptstyle{1}}\!\!\!&{\scriptscriptstyle{-\infty}}\!\!\!&{\scriptscriptstyle{1}}\!\!\!&{\scriptscriptstyle{-\infty}}\!\!\!&{\scriptscriptstyle{-1}}\!\!\!&{\scriptscriptstyle{-\infty}}\!\!\!&{\scriptscriptstyle{0}}\!\!\!&{\scriptscriptstyle{-\infty}}\!\!\!&{\scriptscriptstyle{1}}\!\!\!&{\scriptscriptstyle{-\infty}}\!\!\!&{\scriptscriptstyle{-1}}\!\!\!&{\scriptscriptstyle{-\infty}}\!\!\!&{\scriptscriptstyle{-1}}\!\!\!&{\scriptscriptstyle{-2}}\!\!\!&{\scriptscriptstyle{1}}\!\!\!&{\scriptscriptstyle{4}}\\ 
{\scriptscriptstyle{-\infty}}\!\!\!&{\scriptscriptstyle{-1}}\!\!\!&{\scriptscriptstyle{-\infty}}\!\!\!&{\scriptscriptstyle{1}}\!\!\!&{\scriptscriptstyle{-\infty}}\!\!\!&{\scriptscriptstyle{1}}\!\!\!&{\scriptscriptstyle{2}}\!\!\!&{\scriptscriptstyle{-1}}\!\!\!&{\scriptscriptstyle{-4}}\!\!\!&{\scriptscriptstyle{0}}\!\!\!&{\scriptscriptstyle{4}}\!\!\!&{\scriptscriptstyle{1}}\!\!\!&{\scriptscriptstyle{-2}}\!\!\!&{\scriptscriptstyle{-1}}\!\!\!&{\scriptscriptstyle{-\infty}}\!\!\!&{\scriptscriptstyle{-1}}\!\!\!&{\scriptscriptstyle{-\infty}}\!\!\!&{\scriptscriptstyle{1}}\!\!\!&{\scriptscriptstyle{-\infty}}\!\!\!&{\scriptscriptstyle{0}}\\ 
{\scriptscriptstyle{-1}}\!\!\!&{\scriptscriptstyle{2}}\!\!\!&{\scriptscriptstyle{1}}\!\!\!&{\scriptscriptstyle{-\infty}}\!\!\!&{\scriptscriptstyle{1}}\!\!\!&{\scriptscriptstyle{-\infty}}\!\!\!&{\scriptscriptstyle{-1}}\!\!\!&{\scriptscriptstyle{-\infty}}\!\!\!&{\scriptscriptstyle{0}}\!\!\!&{\scriptscriptstyle{-\infty}}\!\!\!&{\scriptscriptstyle{1}}\!\!\!&{\scriptscriptstyle{-\infty}}\!\!\!&{\scriptscriptstyle{-1}}\!\!\!&{\scriptscriptstyle{-\infty}}\!\!\!&{\scriptscriptstyle{-1}}\!\!\!&{\scriptscriptstyle{-2}}\!\!\!&{\scriptscriptstyle{1}}\!\!\!&{\scriptscriptstyle{4}}\!\!\!&{\scriptscriptstyle{0}}\!\!\!&{\scriptscriptstyle{-4}}\\
{\scriptscriptstyle{-\infty}}\!\!\!&{\scriptscriptstyle{1}}\!\!\!&{\scriptscriptstyle{-\infty}}\!\!\!&{\scriptscriptstyle{1}}\!\!\!&{\scriptscriptstyle{2}}\!\!\!&{\scriptscriptstyle{-1}}\!\!\!&{\scriptscriptstyle{-4}}\!\!\!&{\scriptscriptstyle{0}}\!\!\!&{\scriptscriptstyle{4}}\!\!\!&{\scriptscriptstyle{1}}\!\!\!&{\scriptscriptstyle{-2}}\!\!\!&{\scriptscriptstyle{-1}}\!\!\!&{\scriptscriptstyle{-\infty}}\!\!\!&{\scriptscriptstyle{-1}}\!\!\!&{\scriptscriptstyle{-\infty}}\!\!\!&{\scriptscriptstyle{1}}\!\!\!&{\scriptscriptstyle{-\infty}}\!\!\!&{\scriptscriptstyle{0}}\!\!\!&{\scriptscriptstyle{-\infty}}\!\!\!&{\scriptscriptstyle{-1}}\\
{\scriptscriptstyle{1}}\!\!\!&{\scriptscriptstyle{-\infty}}\!\!\!&{\scriptscriptstyle{1}}\!\!\!&{\scriptscriptstyle{-\infty}}\!\!\!&{\scriptscriptstyle{-1}}\!\!\!&{\scriptscriptstyle{-\infty}}\!\!\!&{\scriptscriptstyle{0}}\!\!\!&{\scriptscriptstyle{-\infty}}\!\!\!&{\scriptscriptstyle{1}}\!\!\!&{\scriptscriptstyle{-\infty}}\!\!\!&{\scriptscriptstyle{-1}}\!\!\!&{\scriptscriptstyle{-\infty}}\!\!\!&{\scriptscriptstyle{-1}}\!\!\!&{\scriptscriptstyle{-2}}\!\!\!&{\scriptscriptstyle{1}}\!\!\!&{\scriptscriptstyle{4}}\!\!\!&{\scriptscriptstyle{0}}\!\!\!&{\scriptscriptstyle{-4}}\!\!\!&{\scriptscriptstyle{-1}}\!\!\!&{\scriptscriptstyle{2}}\\ 
{\scriptscriptstyle{-\infty}}\!\!\!&{\scriptscriptstyle{1}}\!\!\!&{\scriptscriptstyle{2}}\!\!\!&{\scriptscriptstyle{-1}}\!\!\!&{\scriptscriptstyle{-4}}\!\!\!&{\scriptscriptstyle{0}}\!\!\!&{\scriptscriptstyle{4}}\!\!\!&{\scriptscriptstyle{1}}\!\!\!&{\scriptscriptstyle{-2}}\!\!\!&{\scriptscriptstyle{-1}}\!\!\!&{\scriptscriptstyle{-\infty}}\!\!\!&{\scriptscriptstyle{-1}}\!\!\!&{\scriptscriptstyle{-\infty}}\!\!\!&{\scriptscriptstyle{1}}\!\!\!&{\scriptscriptstyle{-\infty}}\!\!\!&{\scriptscriptstyle{0}}\!\!\!&{\scriptscriptstyle{-\infty}}\!\!\!&{\scriptscriptstyle{-1}}\!\!\!&{\scriptscriptstyle{-\infty}}\!\!\!&{\scriptscriptstyle{1}}\\ 
{\scriptscriptstyle{1}}\!\!\!&{\scriptscriptstyle{-\infty}}\!\!\!&{\scriptscriptstyle{-1}}\!\!\!&{\scriptscriptstyle{-\infty}}\!\!\!&{\scriptscriptstyle{0}}\!\!\!&{\scriptscriptstyle{-\infty}}\!\!\!&{\scriptscriptstyle{1}}\!\!\!&{\scriptscriptstyle{-\infty}}\!\!\!&{\scriptscriptstyle{-1}}\!\!\!&{\scriptscriptstyle{-\infty}}\!\!\!&{\scriptscriptstyle{-1}}\!\!\!&{\scriptscriptstyle{-2}}\!\!\!&{\scriptscriptstyle{1}}\!\!\!&{\scriptscriptstyle{4}}\!\!\!&{\scriptscriptstyle{0}}\!\!\!&{\scriptscriptstyle{-4}}\!\!\!&{\scriptscriptstyle{-1}}\!\!\!&{\scriptscriptstyle{2}}\!\!\!&{\scriptscriptstyle{1}}\!\!\!&{\scriptscriptstyle{-\infty}}\\ 
{\scriptscriptstyle{2}}\!\!\!&{\scriptscriptstyle{-1}}\!\!\!&{\scriptscriptstyle{-4}}\!\!\!&{\scriptscriptstyle{0}}\!\!\!&{\scriptscriptstyle{4}}\!\!\!&{\scriptscriptstyle{1}}\!\!\!&{\scriptscriptstyle{-2}}\!\!\!&{\scriptscriptstyle{-1}}\!\!\!&{\scriptscriptstyle{-\infty}}\!\!\!&{\scriptscriptstyle{-1}}\!\!\!&{\scriptscriptstyle{-\infty}}\!\!\!&{\scriptscriptstyle{1}}\!\!\!&{\scriptscriptstyle{-\infty}}\!\!\!&{\scriptscriptstyle{0}}\!\!\!&{\scriptscriptstyle{-\infty}}\!\!\!&{\scriptscriptstyle{-1}}\!\!\!&{\scriptscriptstyle{-\infty}}\!\!\!&{\scriptscriptstyle{1}}\!\!\!&{\scriptscriptstyle{-\infty}}\!\!\!&{\scriptscriptstyle{1}}\\ 
{\scriptscriptstyle{-1}}\!\!\!&{\scriptscriptstyle{-\infty}}\!\!\!&{\scriptscriptstyle{0}}\!\!\!&{\scriptscriptstyle{-\infty}}\!\!\!&{\scriptscriptstyle{1}}\!\!\!&{\scriptscriptstyle{-\infty}}\!\!\!&{\scriptscriptstyle{-1}}\!\!\!&{\scriptscriptstyle{-\infty}}\!\!\!&{\scriptscriptstyle{-1}}\!\!\!&{\scriptscriptstyle{-2}}\!\!\!&{\scriptscriptstyle{1}}\!\!\!&{\scriptscriptstyle{4}}\!\!\!&{\scriptscriptstyle{0}}\!\!\!&{\scriptscriptstyle{-4}}\!\!\!&{\scriptscriptstyle{-1}}\!\!\!&{\scriptscriptstyle{2}}\!\!\!&{\scriptscriptstyle{1}}\!\!\!&{\scriptscriptstyle{-\infty}}\!\!\!&{\scriptscriptstyle{1}}\!\!\!&{\scriptscriptstyle{-\infty}}\\ 
{\scriptscriptstyle{-4}}\!\!\!&{\scriptscriptstyle{0}}\!\!\!&{\scriptscriptstyle{4}}\!\!\!&{\scriptscriptstyle{1}}\!\!\!&{\scriptscriptstyle{-2}}\!\!\!&{\scriptscriptstyle{-1}}\!\!\!&{\scriptscriptstyle{-\infty}}\!\!\!&{\scriptscriptstyle{-1}}\!\!\!&{\scriptscriptstyle{-\infty}}\!\!\!&{\scriptscriptstyle{1}}\!\!\!&{\scriptscriptstyle{-\infty}}\!\!\!&{\scriptscriptstyle{0}}\!\!\!&{\scriptscriptstyle{-\infty}}\!\!\!&{\scriptscriptstyle{-1}}\!\!\!&{\scriptscriptstyle{-\infty}}\!\!\!&{\scriptscriptstyle{1}}\!\!\!&{\scriptscriptstyle{-\infty}}\!\!\!&{\scriptscriptstyle{1}}\!\!\!&{\scriptscriptstyle{2}}\!\!\!&{\scriptscriptstyle{-1}}\\
{\scriptscriptstyle{0}}\!\!\!&{\scriptscriptstyle{-\infty}}\!\!\!&{\scriptscriptstyle{1}}\!\!\!&{\scriptscriptstyle{-\infty}}\!\!\!&{\scriptscriptstyle{-1}}\!\!\!&{\scriptscriptstyle{-\infty}}\!\!\!&{\scriptscriptstyle{-1}}\!\!\!&{\scriptscriptstyle{-2}}\!\!\!&{\scriptscriptstyle{1}}\!\!\!&{\scriptscriptstyle{4}}\!\!\!&{\scriptscriptstyle{0}}\!\!\!&{\scriptscriptstyle{-4}}\!\!\!&{\scriptscriptstyle{-1}}\!\!\!&{\scriptscriptstyle{2}}\!\!\!&{\scriptscriptstyle{1}}\!\!\!&{\scriptscriptstyle{-\infty}}\!\!\!&{\scriptscriptstyle{1}}\!\!\!&{\scriptscriptstyle{-\infty}}\!\!\!&{\scriptscriptstyle{-1}}\!\!\!&{\scriptscriptstyle{-\infty}}
\endmatrix
\right].\tag8.2
$$
Denote the matrix (8.2) by $B=(b_{ij})_{1\leq i,j\leq20}$ (the position 
corresponding to a zero of $A$ is recorded in $B$ by $-\infty$).

The form of $B$ suggests the following generalization of our original weight period
$A$. Consider the matrix $A(q)=(a_{ij}(q))_{1\leq i,j\leq20}$ given by
$$
a_{ij}(q)=q^{b_{ij}}, \ \ \ \ \ 1\leq i,j\leq 20,
$$
with the convention that $q^{-\infty}=0$. By construction, $A(1)$ is our original weight period $A$. 

Applying the operator $\de$ successively thirty times to $A(q)$ (a simple task with the aid of Maple), we obtain 
$$
\de^{(30)}(A(q))=A(9q).\tag8.3
$$
Therefore, by the Reduction Theorem and Lemma 3.5 we obtain that $\M(AD_n;\wt_{A(q)})$ is equal to  
$\M(AD_{n-30};\wt_{A(9q)})$ multiplied by a product of cell-factors. The latter product is obtained
multiplying together all cell-factors from the (thirtyfold!) application of the
Reduction Theorem. With the assistance of Maple it is not hard to find this product explicitly. 
One obtains
$$
\M(AD_n;\wt_{A(q)})=3^{4y_n}\M(AD_{n-30};\wt_{A(9q)}), \ \ \ \ n\geq31,\tag8.4
$$
where $y_{10k+1}=8k-13$, $y_{10k+2}=y_{10k+6}=8k-9$, $y_{10k+3}=y_{10k+7}=8k-7$, $y_{10k+4}=y_{10k+5}=8k-8$, 
$y_{10k+8}=8k-5$ and $y_{10k+9}=y_{10k+10}=8k-4$, for $k\geq3$.

Recurrence (8.4) determines $\M(AD_n;\wt_{A(q)})$ for all $n$ provided we know it for $n\leq30$. 
This can be easily found using Maple to carry out successive applications of the Reduction Theorem.
One obtains the values in Table 8.2. These clearly specialize to the values in Table 8.1 when $q=1$.
It is straightforward to check that for $q=1$ recurrence (8.4) and these initial values imply (8.1).
\sendpf

\topinsert \eightpoint
$$\hbox{\vbox{\offinterlineskip
  \def\spacing{\omit\vrule height4pt&&&&&&\cr}
  \halign{\vrule\strut#&&\ \hfil$#$\hfil\ &\vrule#\cr
  \noalign{\hrule}\spacing
  &n&&\M(AD_n;\wt_{A(q)})&&\text{factorization}&\cr
  \spacing\noalign{\hrule}\spacing
  &1&&\,q^2&&\,q^2&\cr\spacing
  &2&&2\,q^{-2}&&2\,q^{-2}&\cr\spacing
  &3&&6\,q^{-2}&&2\cdot3\,q^{-2}&\cr\spacing
  &4&&6&&2\cdot3&\cr\spacing
  &5&&6&&2\cdot3&\cr\spacing
  &6&&6\,q^2&&2\cdot3\,q^2&\cr\spacing
  &7&&27\,q^2&&3^3\,q^2&\cr\spacing
  &8&&486\,q^2&&2\cdot3^5\,q^2&\cr\spacing
  &9&&486&&2\cdot3^5&\cr\spacing
  &10&&486&&2\cdot3^5&\cr\spacing
  &11&&486\,q^2&&2\cdot3^5\,q^2&\cr\spacing
  &12&&6561\,q^{-2}&&3^8\,q^{-2}&\cr\spacing
  &13&&531441\,q^{-2}&&3^{12}\,q^{-2}&\cr\spacing
  &14&&531441&&3^{12}&\cr\spacing
  &15&&531441&&3^{12}&\cr\spacing  
  \noalign{\hrule}}}}
\qquad
\hbox{\vbox{\offinterlineskip
  \def\spacing{\omit\vrule height4pt&&&&&&\cr}
  \halign{\vrule\strut#&&\ \hfil$#$\hfil\ &\vrule#\cr
  \noalign{\hrule}\spacing
  &n&&\M(AD_n;\wt_{A(q)})&&\text{factorization}&\cr
  \spacing\noalign{\hrule}\spacing
  &16&&531441\,q^2&&3^{12}\,q^2&\cr\spacing  
  &17&&43046721\,q^2&&3^{16}\,q^2&\cr\spacing  
  &18&&10460353203\,q^2&&3^{21}\,q^2&\cr\spacing  
  &19&&10460353203&&3^{21}&\cr\spacing  
  &20&&10460353203&&3^{21}&\cr\spacing  
  &21&&10460353203\,q^2&&3^{21}\,q^2&\cr\spacing  
  &22&&7625597484987\,q^{-2}&&3^{27}\,q^{-2}&\cr\spacing  
  &23&&5559060566555523\,q^{-2}&&3^{33}\,q^{-2}&\cr\spacing  
  &24&&5559060566555523&&3^{33}&\cr\spacing  
  &25&&5559060566555523&&3^{33}&\cr\spacing  
  &26&&5559060566555523\,q^2&&3^{33}\,q^2&\cr\spacing  
  &27&&24315330918113857602\,q^2&&2\cdot3^{40}\,q^2&\cr\spacing  
  &28&&79766443076872509863361\,q^2&&3^{48}\,q^2&\cr\spacing  
  &29&&79766443076872509863361&&3^{48}&\cr\spacing  
  &30&&79766443076872509863361&&3^{48}&\cr\spacing  
  \noalign{\hrule}}}}
$$
\tenpoint\centerline{{\smc Table~8.2.} The first thirty values of $\M(AD_n;\wt_{A(q)})$.}
\endinsert

\mysec{9. Concluding remarks}

As mentioned in the proof of Theorem 6.1, a result analogous to Lemma 4.4 holds for 
the columns of the edge-weight array. It is easy to see that together with Lemma 6.2
this analog of Lemma 4.4 implies that {\it single} columns of the edge-weight array can
be scaled with a simple multiplicative effect on the matching generating function 
of a periodically weighted Aztec diamond. The same is
true for the rows of the edge-weight array. 

A natural equivalence relation arises this way on the set of period matrices.  
In order to obtain a recurrence
for the values $\M(AD_n;\wt_A)$, it is enough to find an iterate $\de^{(k)}(A)$
in the equivalence class of $A$. 

Sometimes it is useful to have an alternative approach, as illustrated in the proof of Theorem 8.1, 
where it was crucial for the proof that we first generalized the original weight---by introduction of the parameter $q$.
This larger class of matrices had the desired property that a convenient iterate (the thirtieth) was
in the same class, yielding thus a recurrence. 

All the problems treated in this paper led, in their unweighted form, to periodic weightings
of Aztec diamonds with periods given by certain 0-1 even by even matrices (except for the case of 
Aztec dungeons and fortresses, when entries equal to 1/2 also arise).
A natural problem to consider would be to classify the 0-1 matrices that give rise
to nice enumeration results when regarded from the point of view of the approach described 
in this paper. Two ways for this to come about are mentioned in the previous two paragraphs.
Given the similarity of the examples considered in Sections 6 and 8, one natural starting point 
would concern sublattices of the square grid obtained by deleting parallel edges in some
more general simple pattern.

\mysec{References}
{\openup 1\jot \frenchspacing\raggedbottom
\roster

\myref{1} 
  M. Ciucu, Enumeration of perfect matchings in graphs with reflective
symmetry, {\it J. Combin. Theory Ser. A}, {\bf 77} (1997), 67-97.
\myref{2} 
  M. Ciucu, A complementation theorem for perfect matchings of graphs having a
cellular completion, {\it J. Combin. Theory Ser. A}, {\bf 81} (1998), 34--68.
\myref{3}
  M. Ciucu and C. Krattenthaler, Enumeration of lozenge tilings of hexagons with cut off
corners, {\it J. Combin. Theory Ser. A}, to appear.
\myref{4} 
  C. Douglas, An illustrative study of the enumeration of tilings: conjecture
discovery and proof techniques, Electronic manuscript dated August 1996 (available at
the time of submission at 
www.math.wisc.edu/$\sim$propp/tiling/www/douglas.ps).
\myref{5}
  N. Elkies, G. Kuperberg, M. Larsen, and J. Propp, Alternating-sign matrices
and domino tilings (Part I), {\it J. Algebraic Combin.} {\bf 1} (1992), 
111--132.
\myref{6} 
  G. Kuperberg, Symmetries of plane partitions and the permanent-de\-ter\-mi\-nant
method, {\it J. Combin. Theory Ser. A} {\bf 68} (1994), 115--151.
\myref{7} 
  P. A. MacMahon, ``Combinatory Analysis,'' vols. 1-2, Cambridge, 1916, reprinted by 
Chelsea, New York, 1960.
\myref{8} 
  J. Propp, Enumeration of matchings: problems and progress, in ``New Perspectives in 
Geometric Combinatorics,'' edited by Billera et al., Cambridge University Press, 1999.
\myref{9} 
  J. Propp, Generalized domino-shuffling, to appear in Theoretical Computer Science 
special issue on tilings.
\myref{10} 
  R. P. Stanley, Private communication.
\myref{11}
  B-Y. Yang, ``Three enumeration problems concerning Aztec diamonds,'' Ph.D. 
thesis, Department of Mathematics, Massachusetts Institute of Technology,
Cambridge, MA, 1991.

\endroster\par}

\enddocument